\documentclass[preprint,12pt]{elsarticle}

\usepackage[pdftex]{hyperref}
\biboptions{sort&compress,numbers,square}

\usepackage{amsfonts}
\usepackage{graphicx}
\usepackage{epsfig}
\usepackage{geometry}
\usepackage{tikz}
\usetikzlibrary{trees,positioning,shapes}
\usepackage{subcaption}
\usepackage{amsmath}
\usepackage{amsthm}

\usepackage{xcolor, colortbl}

\usepackage{algorithm}
\usepackage{algorithmicx, algpseudocode}

\newcommand{\inblue}[1]{\textcolor{blue}{#1}}
\newcommand{\inred}[1]{\textcolor{red}{#1}}
\newcommand{\ingray}[1]{\textcolor{darkgray}{#1}}

\newcommand{\inbrown}[1]{\textcolor{brown}{#1}}

\newcommand{\revc}[1]{\textcolor{black}{#1}}
\newcommand{\revcc}[1]{\textcolor{black}{#1}}

\newcommand{\revMP}[1]{\textcolor{blue}{#1}}

\newtheorem{remark}{Remark}

\newcommand{\Vect}[1]{\mathbf{#1}}

\newcommand{\F}{\Vect{F}}

\newcommand{\mb}[1]{{\mathbb #1}}
\newcommand{\K}{\mathcal{K}}

\begin{document}

\begin{frontmatter}

\title{Isogeometric Residual Minimization Method (iGRM)\\ with Direction Splitting Preconditoner for  Stationary Advection-Dominated Diffusion Problems}

\author{\revcc{V.M.Calo}$^{{\textrm{(3,4)}}}$, M. \L{}o\'{s} $^{{\textrm{(1)}}}$, Q.Deng$^{{\textrm{(3)}}}$, I. Muga$^{{\textrm{(2)}}}$, M. Paszy\'{n}ski $^{{\textrm{(1)}}}$}

\address{$^{\textrm{(1)}}$ Department of Computer Science, \\ AGH University of Science and Technology,
Krakow, Poland \\
e-mail: paszynsk@agh.edu.pl \\
e-mail: marcin.los.91@gmail.com}

\address{$^{\textrm{(2)}}$ Institute of Mathematics, \\ Pontifical Catholic Univeristy of Valpara\'iso, Chile \\
e-mail: ignacio.muga@pucv.cl}

\address{$^{\textrm{(3)}}$ 
Applied Geology, 
Faculty of Science and Engineering, Curtin University, Perth, WA, Australia,\\
e-mail: victor.calo@curtin.edu.au}

\address{$^{\textrm{(4)}}$  Mineral Resources, Commonwealth Scientific and Industrial Research Organisation (CSIRO), \\ Kensington, WA, Australia 6152}

\begin{abstract}
In this paper, we introduce the isoGeometric Residual Minimization (iGRM) method. The method solves stationary advection-dominated diffusion problems. \revMP{We} stabilize the method via residual minimization. We discretize the problem using B-spline basis functions. We then seek to minimize the isogeometric residual over a spline space built on a tensor product mesh. We construct the solution over a smooth subspace of the residual. We can specify the solution subspace by reducing the polynomial order, by increasing the continuity, or by a combination of these. The Gramm matrix for the residual minimization method is \revMP{approximated by a weighted $H^1$ norm,} which we can express as Kronecker products, due to the tensor-product structure of the approximations.  
\revMP{We use the Gramm matrix as a preconditional which can be applied} in a computational cost proportional to the number of degrees of freedom in 2D and 3D. Building on these approximations, we construct an iterative algorithm. We test the residual minimization method on several numerical examples, and we compare it to the Discontinuous Petrov-Galerkin (DPG) and the Streamline Upwind Petrov-Galerkin (SUPG) stabilization methods. \revMP{The iGRM method delivers similar quality solutions as the DPG method, it uses smaller grids, it does not require breaking of the spaces, but it is limited to tensor-product meshes. The computational cost of the iGRM is higher than for SUPG, but it does not require the determination of problem specific parameters.}
\end{abstract}
	
\begin{keyword}
isogeometric analysis \sep residual minimization \sep iteration solvers \sep advection-diffusion simulations \sep linear computational cost \sep preconditioners \end{keyword}

\end{frontmatter}

\section{Introduction}
We introduce the isoGeometric Residual Minimization (iGRM) method, which results in a stable \revMP{discretization} for  advection-dominated diffusion problems. The method applies the residual minimization technique to stabilize the discrete solution. 
The minimum residual methods aim to find $u_h\in U_h$ such that
$$
u_h=\arg\!\!\!\min_{w_h\in U_h}\|b(w_h,\cdot) - \ell(\cdot)\|_{V^*}\,\,,
$$
where $U$ and $V$ are Hilbert spaces, $V^*$ denotes the dual of $V$, $b:U\times V\to \mathbb R$ is a continuous bilinear (weak) form, $U_h\subset U$ is a discrete solution space, and $\ell\in V^{*}$ is a given right-hand side. Several discretization techniques are particular incarnations of this wide-class of residual minimization methods. These include: the \emph{least-squares finite element method}~\cite{ LSFEM}, the \emph{discontinuous Petrov-Galerkin method (DPG) with optimal test functions}~\cite{ t8}, the \emph{variational stabilization method}~\cite{ VSM}, adaptive stabilized finite element method~\cite{ Sergio}, or the \emph{automatic variationally stable finite element method}~\cite{ romkes}.
\revMP{This residual minimization technique exploits the tensor product structure of the discrete space to deliver Uzawa-like iteration schemes to solve the resulting global system. We approach the residual minimization as a saddle point (mixed) formulation, as described in~\cite{Erikkson}. }

The residual minimization method relies on the selection of the residual space and its inner product. In particular, we select a weighted $H^1$ norm as the inner product to obtain a Kronecker product structure of the Gramm matrix. 
\revMP{This structure enables us to obtain a linear computational cost preconditioner for the Uzawa-like iterative solver. The resulting method converges well in the weighted $H^1$ norm.}

For spatial discretization, we use IsoGeometric Analysis (IGA)~\cite{ IGA}. IGA  uses B-splines or NURBS~\cite{ NURBS} basis functions to construct a Galerkin finite element method (FEM). 
In particular, we use a tensor-product B-spline space to approximate the residual minimization method.

This structure enables for fast solution of the projection problem with isogeometric analysis~\cite{ CG1, CG2, Gao2014}.  A similar idea is applied to fast direction splitting solvers of explicit dynamics~\cite{ ED1, ED2, ED3, ED4, ED5}. Splitting methods deliver fast factorizations as discussed, among many other sources, in~\cite{ Samarskii, Quarteroni, ADS1, ADS2, ADS3, ADS4}. A modern version of this method solves different classes of problems~\cite{ Minev1, Minev2}.

\revMP{We exploit the Kronecker product structure of the approximated Gramm matrix to approximate the Schur complement, with iterative corrections in the outer loop. We approximate the inverse of the Schur complement system with a conjugate gradient solver (CG)~\cite{ nocedal2006conjugate} in the inner loop.}

\revcc{In this paper, we blend isogeometric analysis, residual minimization, and the alternating directions splitting solver.  We solve the isogeometric Residual Minimization (iGRM) with a direction splitting preconditioner. We solve stationary advection-diffusion problems.}
We analyze four stationary computational problems, including a problem with the manufactured solution, the Eriksson-Johnson model problem, and a vortical wind problem. We use the Streamline Upwind Petrov-Galerkin (SUPG) and Discontinuous Petrov-Galerkin methods as references to compare the performance of our method.

\revcc{The structure of the paper is the following. We start from the introduction of the iGRM method in Section 2. Next, in Section 3, we describe the general idea and the details of the iterative solver. In Section 4, we present the numerical results for the manufactured solution problem, Erikkson-Johnson model problem, and the vortical wind problem. We conclude with Section 5.}

\section{The Isogeometric Residual Minimization Method (iGRM)}

To simplify the discussion, we focus on a two-dimensional model problem in space. Nevertheless, the extension of the formulation to three-dimensions is straightforward. 
\revcc{We introduce $b:U\times V\to \mathbb R$ a continuous bilinear (weak) form defined over solution $U$ and residual $V$ Hilbert spaces, and $\ell\in V^{*}$ a given right-hand side, where $V^*$ denotes the dual of $V$.}

\subsection{Residual minimization method for the global problem}

For a general weak problem: Find $u\in U$ such as
\begin{equation}
b(u,v)=l(v) \quad \forall v\in V 
\label{eq:gen_weak}
\end{equation}
we define the operator 
%
\begin{equation}
 B: U \rightarrow V^{*}
\end{equation}
such that
\begin{equation}
 \langle Bu , v \rangle_{V^{*} \times V} = b(u,v)
\end{equation}
to reformulate the problem as
\begin{equation}
 Bu - l = 0
\end{equation}
We minimize the residual
\begin{equation}
 u_h = \textrm{argmin}_{w_h \in U_h} {1 \over 2} \| Bw_h - l \|_{V^{*}}^2
\end{equation}
We define the Riesz representer as
\begin{equation}
 R_V \colon V \ni v \rightarrow (v,.) \in V^{*}
\end{equation}
We then project the problem back to $V$
\begin{equation}
 u_h = \textrm{argmin}_{w_h \in U_h} {1 \over 2} \| R_V^{-1} (Bw_h - l) \|_V^2
\end{equation}
\revcc{We use here the definition of the dual norm 
\begin{equation}
\| f\|_{V^{*}}:=\textrm{sup} \{ \langle f,v\rangle_{V^{*}\times V}: v\in V, \|v\|_V=1 \}. 
\end{equation}}
\revcc{The definition of the dual norm implicitly depends on the $V$ norm. Thus, a vital ingredient of a residual minimization method is the selection of a norm on $V$. The quality of the results strongly depends on the choice for the norm on $V$.} The minimum is attained at $u_h$ when the G\^ateaux derivative is equal to $0$ in all directions:
we define the operator 
\begin{equation}
\left( R_V^{-1} (Bu_h - l), R_V^{-1}(B\, w_h) \right)_V = 0 \quad \forall \, w_h \in U_h
\end{equation}
\revcc{We define the Riesz representation of the residual $r=R_V^{-1}(Bu_h-l)$. Thus, our problem reduces to 
\begin{equation}
 \left( r , R_V^{-1} (B\,  w_h )  \right)_V = 0 \quad \forall \, w_h \in U_h
\end{equation}
The action of the inverse of the Riesz operator on the functional $Bw_h$ is equal to that element itself in $V$. Therefore, we get
\begin{equation}
 \langle Bw_h, r  \rangle_{V^{*}\times V} = 0 \quad \quad \forall w_h \in U_h.
\end{equation}}
From the definition of the error representation in terms of the residual, we have 
\begin{equation}
(r,v)_V=\langle B u_h-l,v \rangle_{V^{*}\times V} \quad \forall v\in V.
\end{equation}
Thus, our problem reduces to the following semi-infinite problem: Find $(r,u_h)\in{V\times U_h}$ such that
\begin{eqnarray}
\begin{aligned}
(r,v)_V - \langle B u_h-l,v \rangle_{V^{*}\times V} &= 0 \quad  \forall v\in V \\
\langle Bw_h,r\rangle_{V^{*}\times V} &= 0 \quad  \forall w_h \in U_h
\end{aligned}
\end{eqnarray}
We discretize the residual space as $V_h \in V$ to \revMP{obtain} the discrete problem: Find $(r_h,u_h)\in{V_h\times U_h}$ such as 
\begin{eqnarray}
\begin{aligned}
(r_h,v_h)_{V} - \langle B u_h-l,v_h \rangle_{V^{*}\times V} &= 0 \quad  \forall v_h\in V_h \\
\langle Bw_h,r_h\rangle_{V^{*}\times V} &= 0 \quad \forall w_h \in U_h
\end{aligned}
\label{eq:resmin}
\end{eqnarray}
where $(\cdot,\cdot)_{V}$ is the inner product in $V_h$ inherited from $V$, $\langle Bu_h,v_h\rangle_{V^{*}\times V} = b\left(u_h,v_h\right)$, $\langle Bw_h,r_h\rangle_{V^{*}\times V} = b\left(w_h,r_h\right)$. 
\begin{remark}
We define the discrete residual space $V_h$ to be sufficiently close to the continuous $V$ space to ensure stability. Indeed, we gain stability by enriching the residual space $V_h$ (while fixing the solution space $U_h$) until we reach the discrete inf-sup stability.
\end{remark}

\subsection{Minimal residual discretization for the global problem with 
B-splines}

We approximate the solution with a tensor product of one-dimensional B-splines. We choose a uniform order $p$ in all directions to simplify the discussion.

We denote the basis functions in the $x$-direction for the discrete solution and residual spaces as $n_a$ and $N_A$, respectively.  Similarly, we denote the basis functions in the $y$-direction for the discrete solution and residual spaces as $m_b$ and $M_B$, respectively.

\revMP{We} approximate the residual space with a \revMP{sufficiently} larger discrete space \revMP{than} the solution space to guarantee discrete stability (cf.~\cite{Sergio}). We achieve this by either increasing the polynomial order or decreasing the continuity of the residual space to define the solution space.  As a result, we now can define the discrete solution $u_h$ and residual $w_h$ functions as:
\begin{equation} 
\begin{aligned}
  u_h &= \sum_{a,b} w_{ab} n_{a}m_{b} \\
  w_h &= \sum_{A,B} u_{AB} N_{A}M_{B}
\end{aligned}
\end{equation}
and the discrete solution $r_h$ and residual $v_h$ functions as
\begin{equation}
\begin{aligned}
r_h &= \sum_{a,b} r_{ab} n_{a}m_{b} \\
v_h &= \sum_{A,B} v_{AB} N_{A}M_{B}
\end{aligned}
\end{equation}


\newcommand{\dK}{\partial K}
\newcommand{\FK}{\F_{\dK}}
\newcommand{\ZF}{Z_{\F,0}}
\newcommand{\ZFd}{Z_{\F,0}^\dagger}
\newcommand{\VKc}{V_\K^{k,\mathrm{c}}}
\newcommand{\eKc}{e_\K^{k,\mathrm{c}}}

\section{Conjugate Gradients method for Schur complement}

In this section, we derive an iterative solution algorithm for the residual minimization problem. We denote by $\Omega$ a bounded open set of $\mb{R}^d$ with Lipschitz continuous boundary $\partial \Omega$, where $d=2,3$.  \revMP{We} assume that the domain $\Omega = \Omega_x \times \Omega_y$ has a tensor product structure. This structure allows us to approximate the Gramm matrix as a Kronecker product of one-dimensional matrices.

The resulting matrix system for the residual minimization method is
%
\begin{equation} \label{eq:mp}
\begin{bmatrix}
G & B \\
B^T & 0 \\
\end{bmatrix}
\begin{bmatrix}
r \\
u \\
\end{bmatrix} 
= 
\begin{bmatrix}
F \\
0 \\
\end{bmatrix},
\end{equation}
where we minimize the isogeometric discretization residual of~\eqref{eq:resmin} in the energy norm $G$ (a weighted $H^1$ norm), where $G$ is defined as
%
\begin{equation} \label{eq:A}
G = M + \eta K 
\end{equation}
\revMP{where we keep the weight $\eta$ as a solver parameter, and }
\revcc{\begin{equation}
\begin{aligned}
M & = \int_{\Omega_x} \int_{\Omega_y}  n_a m_b N_C M_D \,dx \,dy \\
& = \int_{\Omega_x} n_a N_C \,dx \int_{\Omega_y} m_b M_D \,dy = M_x \otimes M_y, \\
K & = K_x \otimes M_y + M_x \otimes K_y,
\end{aligned}
\end{equation}
where $M_x, M_y$, and $K_x, K_y$ denote 1D mass matrices and stiffness matrices, respectively, defined as 
\begin{equation}
\begin{aligned}
\{ M_x \}_{aC} &= m\left( n_a , N_C \right)= \int_{\Omega_x} n_a N_C \, dx \\
\{ M_y \}_{bD} &= m\left( n_b , N_D \right)=\int_{\Omega_y} m_b M_D \, dy \\
\{ K_x \}_{aC} &= a\left( n_a , N_C \right)=\int_{\Omega_x} \frac{d n_a}{dx} \frac{d N_C}{dx} \, dx \\
\{ K_y \}_{bD} &= a\left( m_b , M_D \right)=\int_{\Omega_y} \frac{d m_b}{dy} \frac{d M_D}{dy} \,dy. \\
\end{aligned}
\end{equation}}

\revcc{We approximate the Gramm matrix in the residual minimization as a Kronecker product matrix $\tilde{G} =\left(K_x+\eta M_x\right) \otimes \left(K_y + \eta M_y\right)$, that is inexpensive to solve \revMP{(it requires solutions of two one-dimensional problems with multiple right-hand sides, where the sizes of the problems and the number of right-hand-sides corresponds to the mesh dimensions in both directions)} and only introduces an error of order $\revc{\eta^2}$. From~\eqref{eq:A}, we have}
\begin{equation}
\begin{aligned}
G & =  (M_x +  \eta K_x) \otimes (M_y +  \eta K_y) - \eta^2 K_x \otimes K_y \nonumber\\
&=\tilde{G} - \tilde K.
\end{aligned} \label{eq:Atilde}
\end{equation}
where $\tilde K = \eta^2 K_x \otimes K_y$. 

\revcc{The definition of the $B$ matrix in~\eqref{eq:mp} corresponds to applying an appropriate isogeometric weak form of the scalar field PDE.}


\subsection{Overview of the iterative algoritm}

\revcc{In this section, we describe the iterative algorithm. Given an initial guess $\begin{bmatrix}
r^k  \\
u^k \end{bmatrix}$, we compute the update as
\begin{equation} \label{eq:update}
\begin{bmatrix}
d  \\
c \\
\end{bmatrix} = 
\begin{bmatrix}
r-r^k  \\
u-u^k \\
\end{bmatrix} \nonumber
\end{equation}}
\revcc{Thus
\begin{equation} \label{eq:get_update}
\begin{bmatrix}
G & B \\
B^T & 0 \\
\end{bmatrix}
\begin{bmatrix}
d \\
c \\
\end{bmatrix} 
= 
\begin{bmatrix}
F - G r^k - Bu^k \\
-B^T r^k \\
\end{bmatrix}  \nonumber
\end{equation}
This system of linear equations is expensive to factorize, so we use the approximate $\tilde{G}$ 
\begin{equation} \label{eq:solve_update}
\begin{bmatrix}
\tilde{G} & B \\
B^T & 0 \\
\end{bmatrix}
\begin{bmatrix}
d \\
c \\
\end{bmatrix} 
= 
\begin{bmatrix}
F - G r^k - Bu^k \\
-B^T r^k \\
\end{bmatrix} \nonumber
\end{equation}}
\revcc{We propose the following iterative algorithm}
\begin{algorithm} [H]
\caption{\revcc{Iterative algorithm}} \label{alg:alg1}
\begin{algorithmic}
  \State \revcc{Initialize $\{u^{0} = 0; r^{0} = 0 \}$} \\
  \revcc{\textbf{for }{$k=1 \ldots N$ until convergence}}
  \begin{equation}
    \begin{aligned}
      \revcc{\textrm{Compute Schur }} & \revcc{\textrm {complement with linear } {\cal O}(N) 
      \textrm{ cost}}\\
    \revcc{  \begin{bmatrix}
 & \tilde{G} & B \\
 & B^T       & 0 \\
      \end{bmatrix}
      \begin{bmatrix}
 & d^k \\ c^k
      \end{bmatrix}}
      & \revcc{=
      \begin{bmatrix}
 & F - G r^k - B u^k \\
 & -B^T r^k \\
      \end{bmatrix}} \\
      & \revcc{\textrm{Solve}} \\
 \revcc{ B^T \tilde{G}^{-1} B\, c^k} & \revcc{=
  B^T \tilde{G}^{-1}(F + \tilde{K} r^k - B u^k)} \\
      &  \revcc{\textrm{with PCG or \revMP{Multi-frontal solver}}} \\
      \revcc{r^{k+1}} & \revcc{= d^k+r^k} \\
      \revcc{u^{k+1}} & \revcc{= c^k+u^k} \\
      \revcc{k} & \revcc{\gets k+1} \\
    \end{aligned}
    \label{eq:CGouter}
  \end{equation}
\end{algorithmic}
\end{algorithm}

\subsection{Convergence of the iterative algorithm}

Since the Algorithm~\ref{alg:alg1} uses a preconditioned CG and both the CG and preconditioned CG are convergent (see, for example,~\cite{ nocedal2006conjugate}), we focus on the spectral analysis of~\eqref{eq:CGouter} in Algorithm~\ref{alg:alg1}.

Applying the initialization $u^{(0)} = u^k$\revMP{,} we obtain
\begin{equation} \label{eq:cv1}
u^{k+1} = u^k + c^k
\end{equation}
where $c^k$ is the update of $u^k$.

We now have
\begin{equation} \label{eq:cv2}
r^{k+1} = \tilde{G}^{-1} ( F + \tilde K r^k  - B u^k) - \tilde{G}^{-1} B c^k.
\end{equation}

To simplify the spectral analysis and without loss of generality\revMP{,} we set $F=0$. Thus, combining~\eqref{eq:cv1} and~\eqref{eq:cv2} gives
\begin{equation} \label{eq:cvp}
\begin{bmatrix}
u^{k+1} \\
r^{k+1} \\
\end{bmatrix}
=
\begin{bmatrix}
1 & 0 \\
-\tilde{G}^{-1}B & \tilde{G}^{-1} \tilde K \\
\end{bmatrix}
\begin{bmatrix}
u^k \\
r^k \\
\end{bmatrix} 
+
\begin{bmatrix}
c^k \\
-\tilde{G}^{-1}B c^k \\
\end{bmatrix}.
\end{equation}

Now, we analyze the spectrum of 
\begin{equation} \label{eq:A0}
\tilde{G}^{-1} \tilde K = (M_x +  \revc{\eta} K_x)^{-1} \otimes (M_y +  \revc{\eta} K_y)^{-1} \cdot \left(\revc{\eta} K_x \otimes \revc{\eta} K_y\right).
\end{equation}

We apply a spectral decomposition~\cite{ horn1990matrix} to the matrix $K_\xi, \xi=x,y$ with respect to $M_\xi$ and arrive at
\begin{equation} \label{eq:sd0}
K_\xi = M_\xi P_\xi D_\xi P_\xi^{-1},
\end{equation}
where $D_\xi$ is a diagonal matrix with entries to be the eigenvalues of the generalized eigenvalue problem
\begin{equation} \label{eq:eigKM0}
K_\xi v_\xi = \lambda_\xi M_\xi v_\xi
\end{equation}
and the columns of $P_\xi$ are the eigenvectors. We assume that all the eigenvalues are sorted in ascending order and are listed in $D_\xi$ and the $j$-th column of $P_\xi$ is associated with the eigenvalue $\lambda_{\xi,j} = D_{\xi,jj}$.

Using~\eqref{eq:sd0} and~\eqref{eq:A0}, we now calculate 
\begin{equation}
\begin{aligned}
\tilde{G}^{-1} & = (M_x +  \revc{\eta} K_x)^{-1} \otimes (M_y +  \revc{\eta} K_y)^{-1} \\
& = (M_x + \revc{\eta} M_x P_x D_x P_x^{-1} )^{-1} \otimes (M_y + \revc{\eta} M_y P_y D_y P_y^{-1} )^{-1} \\
& = P_x E_x P_x^{-1} M_x^{-1} \otimes P_y E_y P_y^{-1} M_y^{-1},
\end{aligned}
\end{equation}
where
\begin{equation}
E_\xi = ( I + \revc{\eta} D_\xi )^{-1}, \qquad \xi = x, y.
\end{equation}
We assume here that the mesh is uniform in both directions.

Thus, similarly, we have
\begin{equation}
\begin{aligned}
\tilde{G}^{-1} \tilde K & = \big( P_x E_x  P_x^{-1} M_x^{-1} \otimes P_y E_y P_y^{-1} M_y^{-1} \big)  \\
& \qquad  \cdot \big( \revc{\eta} M_x P_x D_x P_x^{-1} \otimes \revc{\eta} M_y P_y D_y P_y^{-1} \big) \\
& = \big( P_x \otimes P_y \big) \cdot \big( \revc{\eta}E_x D_x \otimes \revc{\eta} E_yD_y  \big) \cdot \big( P_x^{-1} \otimes P_y^{-1} \big) \\
& = \big( P_x \otimes P_y \big) \cdot \big( \revc{\eta}( I + \revc{\eta} D_x )^{-1} D_x \otimes \revc{\eta} ( I + \revc{\eta} D_y )^{-1} D_y  \big) \cdot \big( P_x^{-1} \otimes P_y^{-1} \big). \\
\end{aligned}
\end{equation}

The middle term is a diagonal matrix. Thus, a typical eigenvalue of $\tilde A^{-1} \tilde K$ is
\begin{equation}
\lambda = \frac{\revc{\eta^2} \lambda_{x,i} \lambda_{y,j}}{ (1+\revc{\eta} \lambda_{x,i}) (1+\revc{\eta} \lambda_{y,i})},
\end{equation}
where $i,j$ are indices of eigenvalues in each dimension. The spectral radius of $\tilde{G}^{-1} \tilde K$ is then
\begin{equation}
\rho = \frac{\revc{\eta^2} \lambda_{x,\max} \lambda_{y,\max}}{ (1+\revc{\eta} \lambda_{x,\max}) (1+\revc{\eta} \lambda_{y,\max})},
\end{equation}
where $\lambda_{\xi,\max}, \xi = x,y$ are the maximum eigenvalues in each dimension.  Immediately, we have
\begin{equation}
0 < \lambda \le \rho <1.
\end{equation}

Thus, the eigenvalues of the amplifying block-matrix (the vector in terms of $c^k$ is from inner CG) in~\eqref{eq:cvp} are 1 and $\lambda$. Thus, they are \revMP{bounded} by 1, which implies that the eigenvalues of their powers are also bounded by 1. Hence, the iterative Algorithm~\ref{alg:alg1} is convergent.


\section{Numerical results for stationary problems}

\revcc{In this section, we focus on a model advection-diffusion problem defined as follows. For a unitary square domain $\Omega=(0,1)^2$, the advection vector $\beta$ and the diffusion coefficient $\epsilon$, we seek the solution of the advection-diffusion equation
\begin{equation}
\beta_x\frac{\partial u}{\partial x}+\beta_y\frac{\partial u}{\partial y}-\epsilon \left(\frac{\partial^2 u}{\partial x^2}+
\frac{\partial^2 u}{\partial y^2}\right)=f
\label{eq:ModelProblem}
\end{equation}
with Dirichlet boundary conditions $u=g$ on the whole boundary of $\Gamma=\partial\Omega$.  We partition the boundary $\Gamma=\partial \Omega$ into the inflow $\Gamma^{-}=\{ x\in\Gamma: \beta \cdot n<0 \}=\{(x,y): xy=0\}$ and the outflow $\Gamma^{+}=\{x\in\Gamma: \beta \cdot n \geq 0\}$.}

\revcc{We introduce the discrete weak formulation 
\begin{equation}
b(u_h,v_h)=l(v_h) \quad \forall v\in V 
\label{eq:ModelProblem_weak}
\end{equation}
\begin{eqnarray}
\begin{aligned}
b(u_h,v_h) &=
\beta_x\left(\frac{\partial u_h}{\partial x},v_h\right)_{\Omega}+\beta_y\left(\frac{\partial u_h}{\partial y},v_h\right)_{\Omega}+\epsilon \left(  \frac{\partial u_h}{\partial x}, \frac{\partial v_h}{\partial x}\right)_{\Omega} +\epsilon \left(  \frac{\partial u_h }{\partial y}, \frac{\partial v_h}{\partial y}\right)_{\Omega} \\
&\inbrown{-\left(\epsilon\frac{\partial u_h}{\partial x}n_x,v_h\right)_{\Gamma}
-\left(\epsilon\frac{\partial u_h}{\partial y}n_y,v_h\right)_{\Gamma}}  \nonumber\\
&\inred{-\left(u_h,\epsilon \nabla v_h \cdot n\right)_{\Gamma} +\left(u_h,\beta\cdot n v_h\right)_{\Gamma^{-}}}\ingray{-\sum_K\left(u_h,3p^2 \epsilon/h_K v_h\right)_{\Gamma|_K} } \nonumber 
\label{eq:ModelProblem_b}
\end{aligned}
\end{eqnarray}
where $n=(n_x,n_y)$ is the \revMP{vector} normal to $\Gamma$, and $h_K$ is the \revMP{``normal''} elemental distance (i.e., mesh size in the \revMP{direction} of the normal) (possibly changing with elements $K$ of the mesh),
\begin{equation}
l(v_h)=(f,v_h)_{\Omega}\inred{-\left(g,\epsilon\nabla v_h \cdot n\right)_{\Gamma}
+\left(g,\beta\cdot n v_h\right)_{\Gamma^{-}}}-\sum_K\ingray{\left(g,3p^2 \epsilon/h_K v_h\right)_{\Gamma}}\revMP{,}
\label{eq:ModelProblem_l}
\end{equation}
where the \ingray{gray} and \inred{red} represent the penalty terms, and $f$ corresponds to the forcing term, while the \inbrown{brown} terms result from the integration by parts~\cite{ weak}.}

In our model problem, we seek the solution in space $U = V = H^1\left(\Omega\right)$. The inner product in $V$ is defined as weighted $H^1$ norm (\ref{eq:A}).

\revcc{A popular stabilization technique is the Streamline-upwind Petrov-Galerkin (SUPG) method~\cite{ SUPG2, SUPG}. This method modifies the weak form as follows
\begin{equation}
b(u_h,v_h) + \sum_K \inred{(R(u_h),\tau \beta\cdot \nabla v_h)_K}=l(v_h)
\revMP{+\sum_K \inblue{(f,\tau \beta\cdot \nabla v_h)_K}}
 \quad \forall v\in V 
\label{eq:Erikkson_weak_SUPG}
\end{equation}
where $R(u_h)=\beta \cdot \nabla u_h +\epsilon \Delta u_h$, and $\tau^{-1}= \beta \cdot \left(\frac{1}{h^x_K},\frac{1}{h^y_K} \right) + 3p^2\epsilon \frac{1}{{h^x_K}^2+h{^y_K}^2}$, where $\epsilon$ stands for the diffusion term, and $\beta = (\beta_x,\beta_y)$ for the convection term, and $h^x_K$ and $h^y_K$ are horizontal and vertical dimensions of an element $K$.  Thus, we have
\begin{equation}
b_{\inred{SUPG}}(u_h,v_h)=l_{\inblue{SUPG}}(v_h) \quad \forall v_h\in V_h
\label{eq:SUPG_weak}
\end{equation}
\begin{eqnarray}
\begin{aligned}
b_{\inred{SUPG}}(u_h,v_h)=
& \left(\frac{\partial u_h}{\partial x},v\right)+\epsilon \left(  \frac{\partial u_h}{\partial x}, \frac{\partial v_h}{\partial x}\right) +\epsilon \left(  \frac{\partial u_h }{\partial y}, \frac{\partial v_h}{\partial y}\right)  \\
&-\left(\epsilon\frac{\partial u_h}{\partial x}n_x,v_h\right)_{\Gamma}
-\left(\epsilon\frac{\partial u_h}{\partial y}n_y,v_h\right)_{\Gamma}\nonumber\\
&-\left(u_h,\epsilon \nabla v_h \cdot n\right)_{\Gamma} -\left(u_h,\beta\cdot n v_h\right)_{\Gamma}-\left(u_h,3p^2 \epsilon/h v_h\right)_{\Gamma}  \nonumber 
\\
&+\inred{\left(\frac{\partial u_h}{\partial x}+\epsilon \Delta u_h, 
\left(\frac{1}{h_x}  + 3\epsilon \frac{p^2}{{h^x_K}^2+{h^y_K}^2}\right)^{-1} \frac{\partial v_h}{\partial x}
\right)}
\label{eq:Erikkson_b_SUPG}
\end{aligned}
\end{eqnarray}}
\revMP{\begin{eqnarray}
\begin{aligned}
l_{\inblue{SUPG}}(v_h)&=(f,v_h)-\left(g,\epsilon\nabla v_h \cdot n\right)_{\Gamma} \nonumber \\
&+\left(g,\beta\cdot n v_h\right)_{\Gamma^{-}}-\left(g,3p^2 \epsilon/h_K v_h\right)_{\Gamma}
\revMP{+\inblue{\left(f,\left(\frac{1}{h_x}  + 3\epsilon \frac{p^2}{{h^x_K}^2+{h^y_K}^2}\right)^{-1}\frac{\partial v_h}{\partial x}\right)}}
\revMP{.} \nonumber
\label{eq:ModelProblem_l}
\end{aligned}
\end{eqnarray}}


\subsection{A manufactured solution problem}

First, we select the advection vector $\beta=(1,1)^T$, and $Pe=1/\epsilon=100$ and solve the advection-diffusion equation (\ref{eq:ModelProblem}) with homogeneous Dirichlet boundary conditions.  We utilize a manufactured solution $u(x,y)=\left(x+\frac{e^{Pe*x}-1}{1-e^{Pe}}\right)\left(y+\frac{e^{Pe*y}-1}{1-e^{Pe}}\right)$ enforced by the forcing term $f$.  \revcc{This analytic expression of the solution limits the P\'eclet number to $Pe=100$ due to machine precision.}

We use the weak form (\ref{eq:ModelProblem_weak}) and the inner product (\ref{eq:A}) into the iGRM setup (\ref{eq:resmin}) and use the preconditioned CG solver described in Section 3.
\revcc{In the weighted $H^1$ \revMP{norm} we introduce $\eta=h^2$, where $h$ is the diameter of the element.}

In the following problem, we study the $h$- and $p$-convergence of the iGRM method on uniform grids using different combinations of solution and residual spaces. We do not employ adaptive Shishkin grids here~\cite{ Kopteva}.  We increase the accuracy by increasing the order and continuity of solution spaces ($k$-\revMP{increase}~\cite{IGA, IGA2}) simultaneously. \revMP{Increasing} solution space order and continuity improves the accuracy of the solution \revMP{for this smooth problem}. Similarly, refining the mesh also improves the accuracy of the solution.

Table~\ref{tab:Problem1} illustrates the $h$ and $p$-convergence of the method.  The rows represent $p$-refinement of the solution test, from $(p,p-1)$ (order $p$, continuity $p-1$) to $(p+1,p)$ (order $p+1$, continuity p), and the columns represent $h$-refinement, from $n\times n$ mesh to $2n \times 2n$ mesh.  We report the numbers of degrees of freedom and the error in $L^2$ and $H^1$ norms.

\revcc{The black colors on the numerical results represent unwanted negative values of the solution.  We can read from Table~\ref{tab:Problem1} that $k$-refinement (increasing polynomial order of the solution with maximum continuity) is more attractive than $h$-refinement since the solution improves using fewer degrees of freedom.}

\revcc{We also present in Table~\ref{tab:Problem1SUPG} the results for SUPG method for the corresponding meshes and approximation spaces. \revMP{We conclude that we can achieve less than 1 percent error}, as measured in $L^2$ and $H^1$ norms\revMP{,} for iGRM method for $\#DOF=5594$ for trial (5,4) test (2,0) over $32\times 32$ mesh, and for SUPG method for $\#DOF=4761$ for trial=test=(5,4) over $64\times 64 $ mesh.}

{\tiny\begin{table*}[htp]
\begin{center}
\begin{tabular}{|c|c|c|c|c|}
\hline 
n & trial (2,1)  & trial (3,2)  & trial (4,3) & trial (5,4) \\
 & test (2,0) & test (2,0) & test (2,0)  & test (2,0)  \\
\hline
\#DOF & 389 & 410 & 433 & 458 \\
L2 & 192.47 & 151.23 & 78.69 & 28.11 \\
H1 & 101.14 & 74.54 & 44.33 & 32.05 \\
$8\times 8$ & 
\includegraphics[scale=0.16]{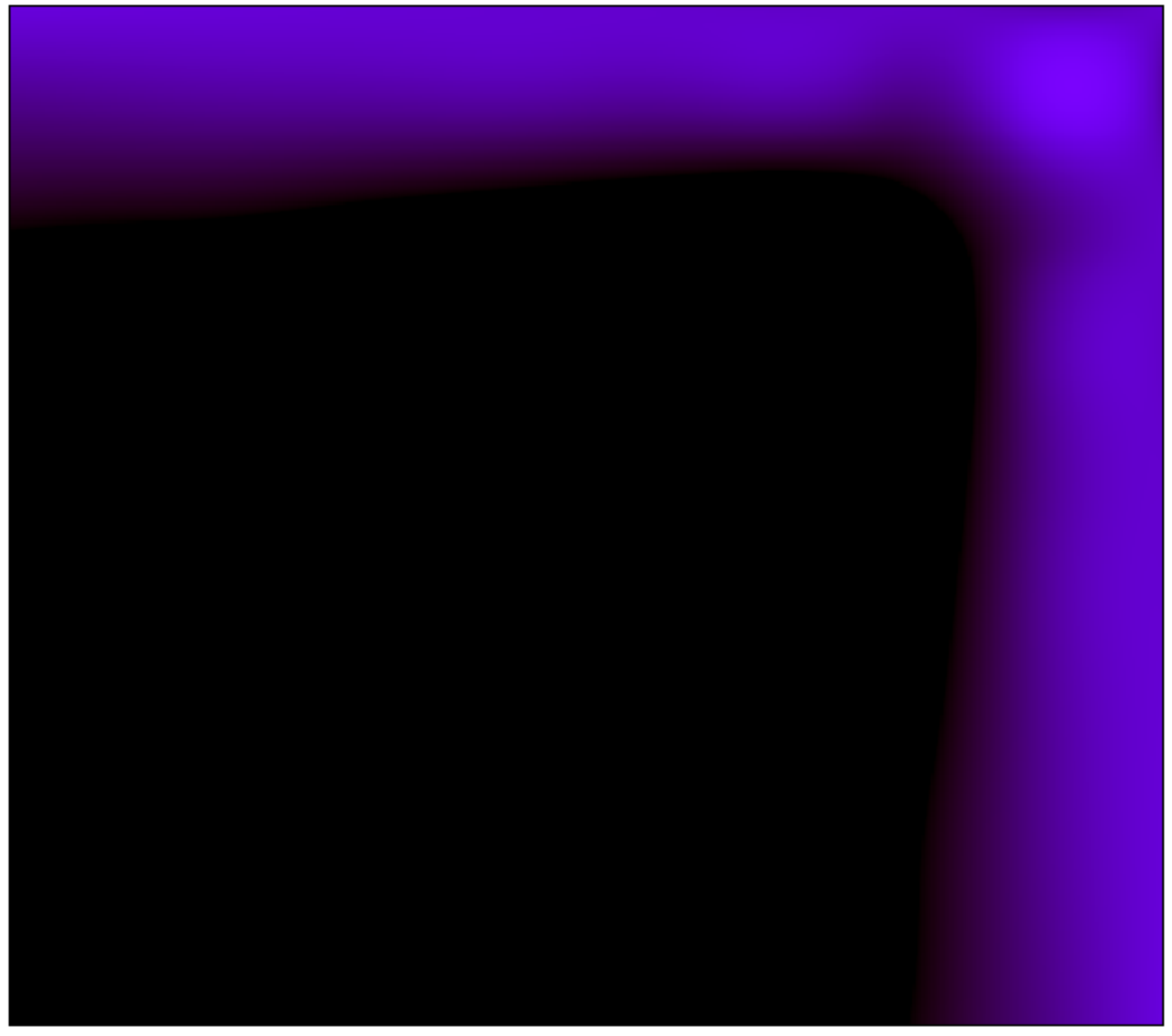}
&
\includegraphics[scale=0.16]{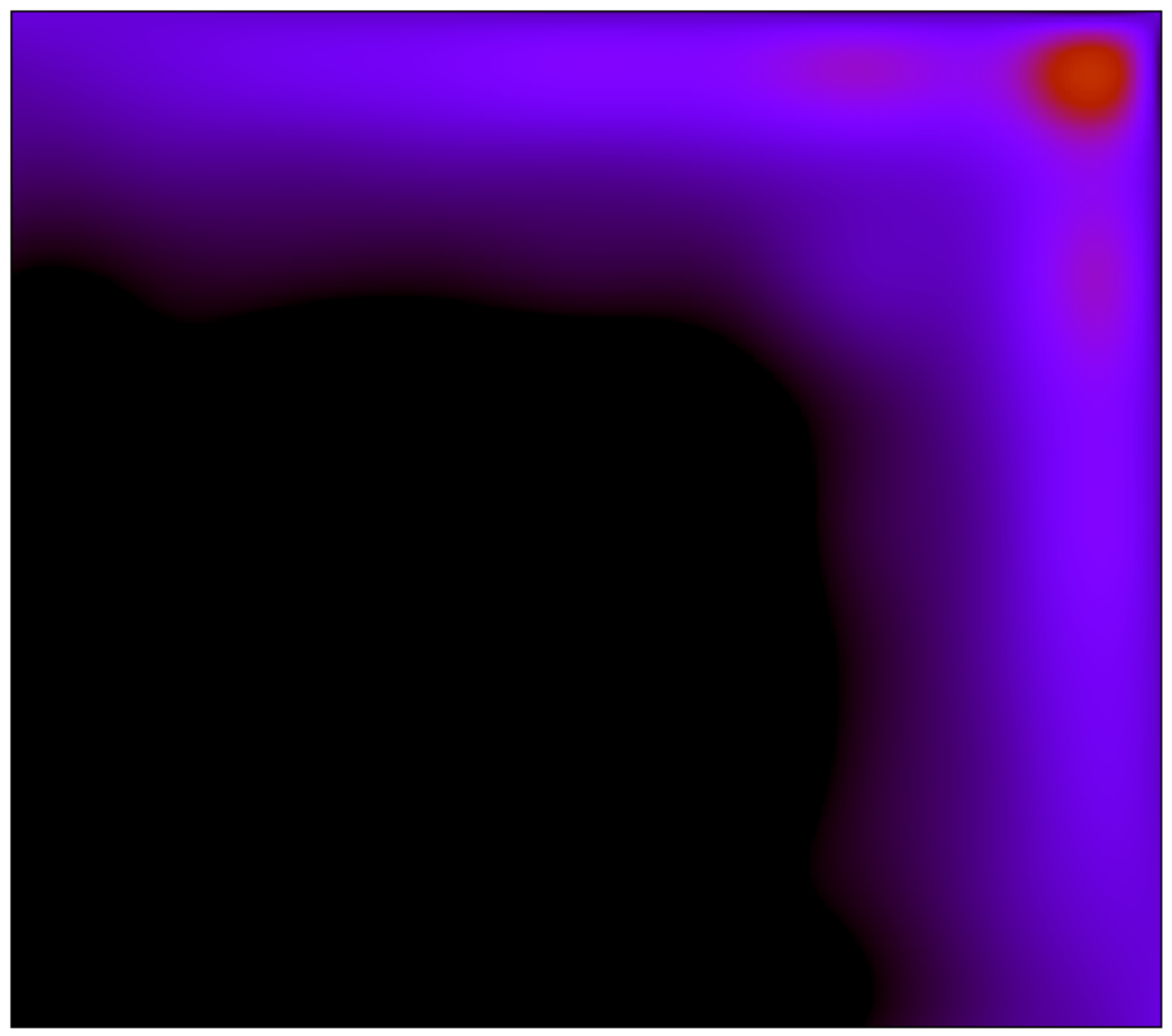}
&
\includegraphics[scale=0.16]{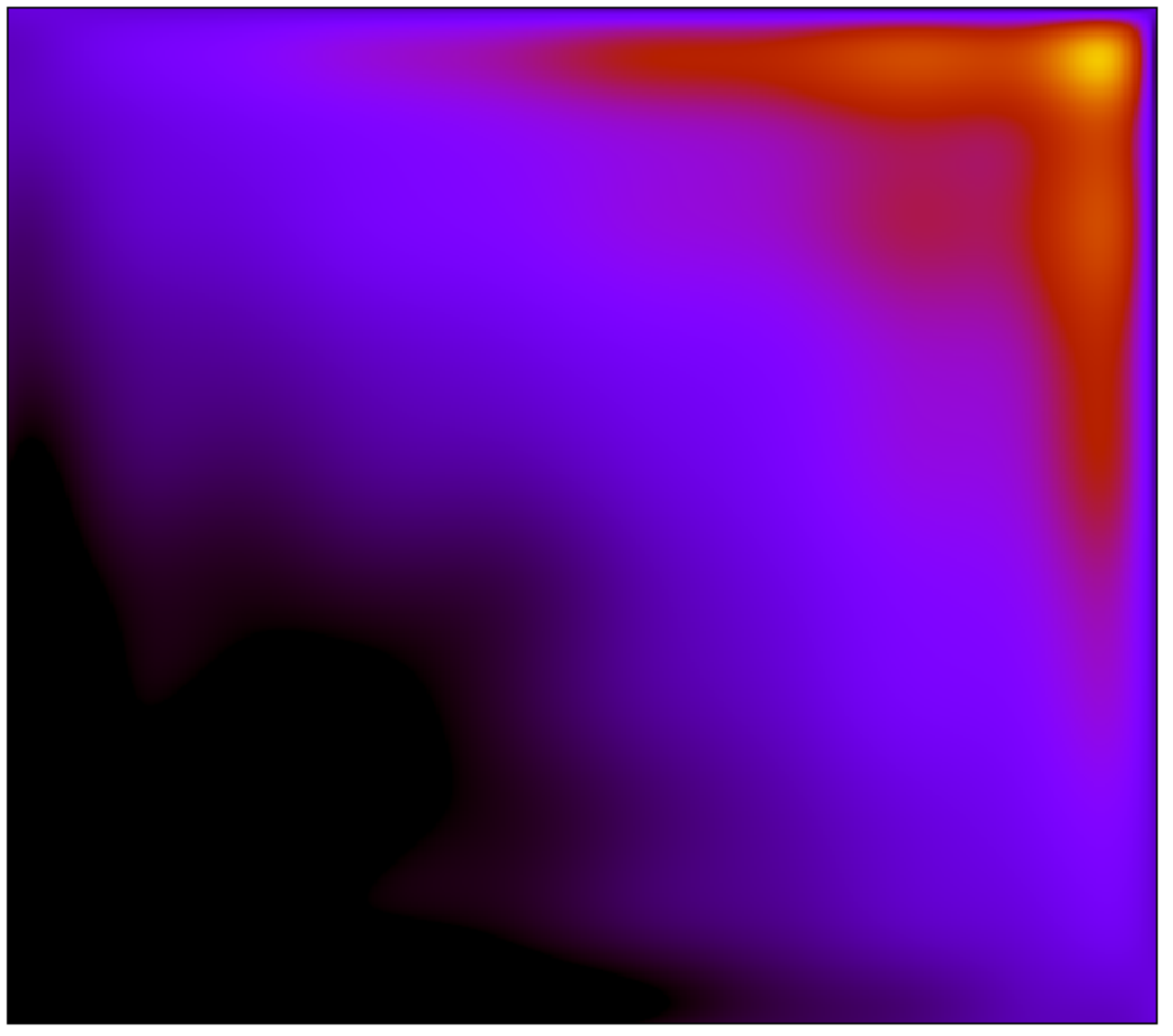}
&
\includegraphics[scale=0.16]{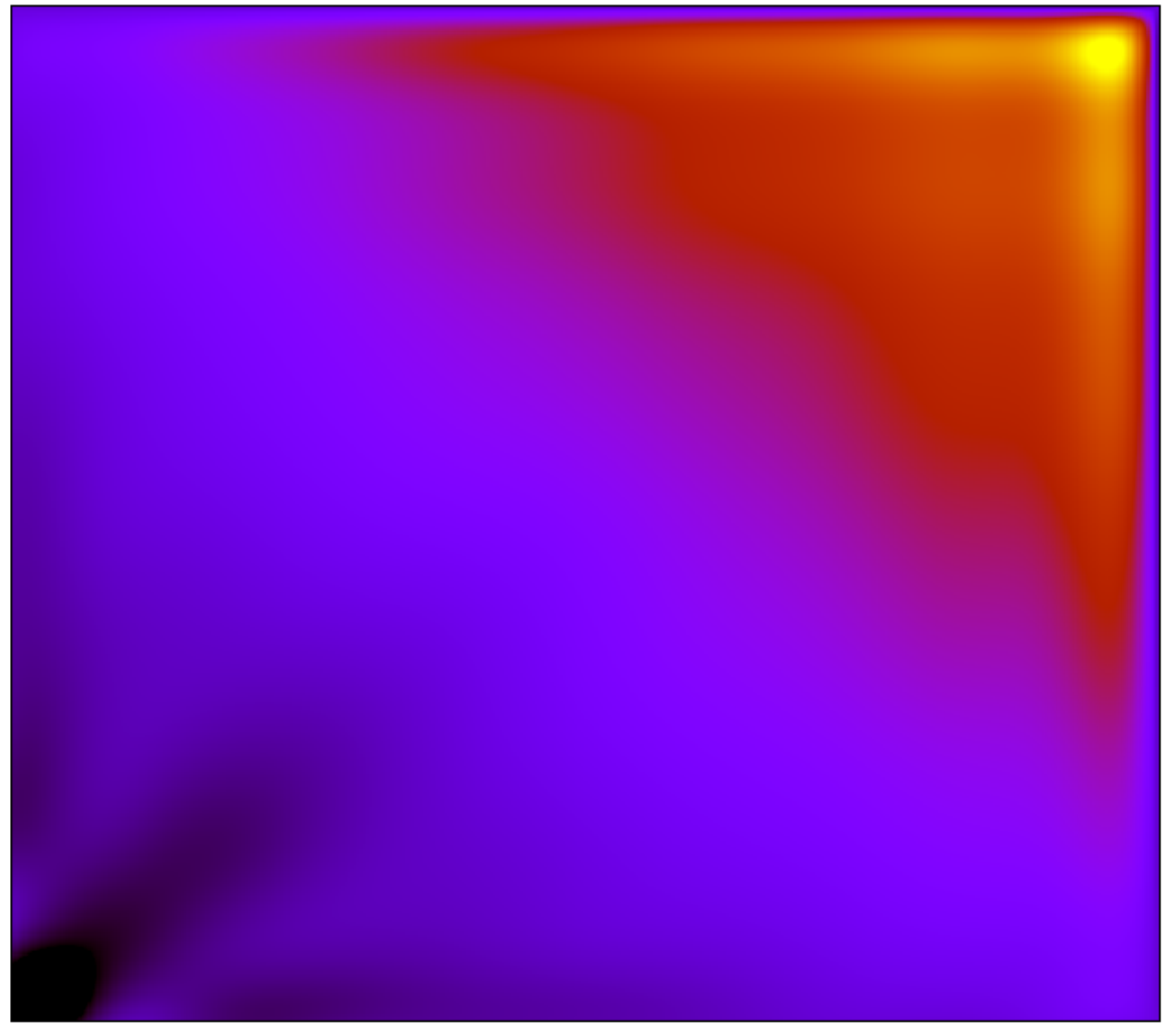}
\\
\hline
\#DOF & 1413 & 1450 & 1489 & 1530 \\
L2 & 80.01 & 16.64 & 3.29 & 1.48 \\
H1 & 59.56 & 29.83 & 18.04 & 10.40
\\
$16\times 16$ & 
\includegraphics[scale=0.16]{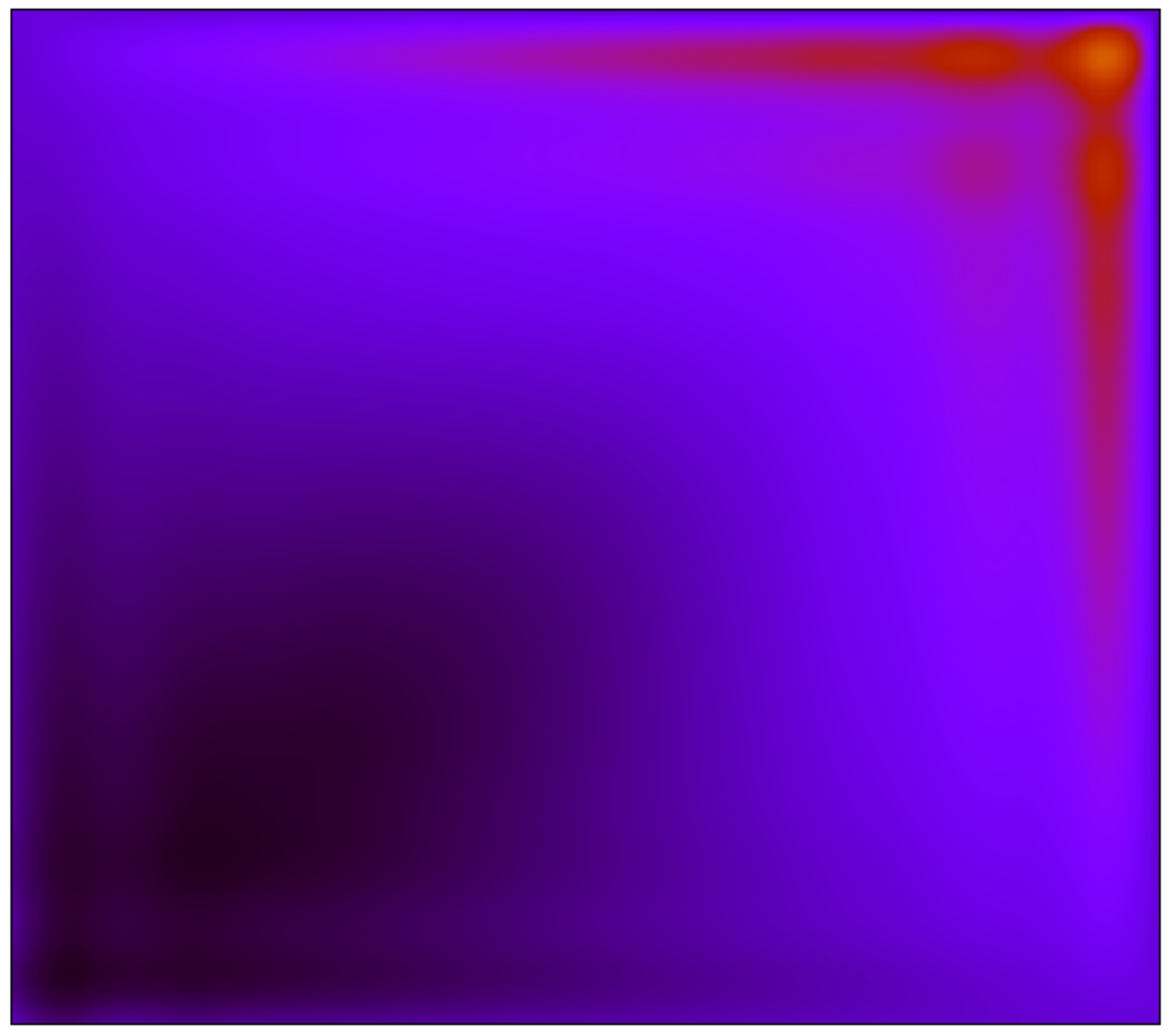}
&
\includegraphics[scale=0.16]{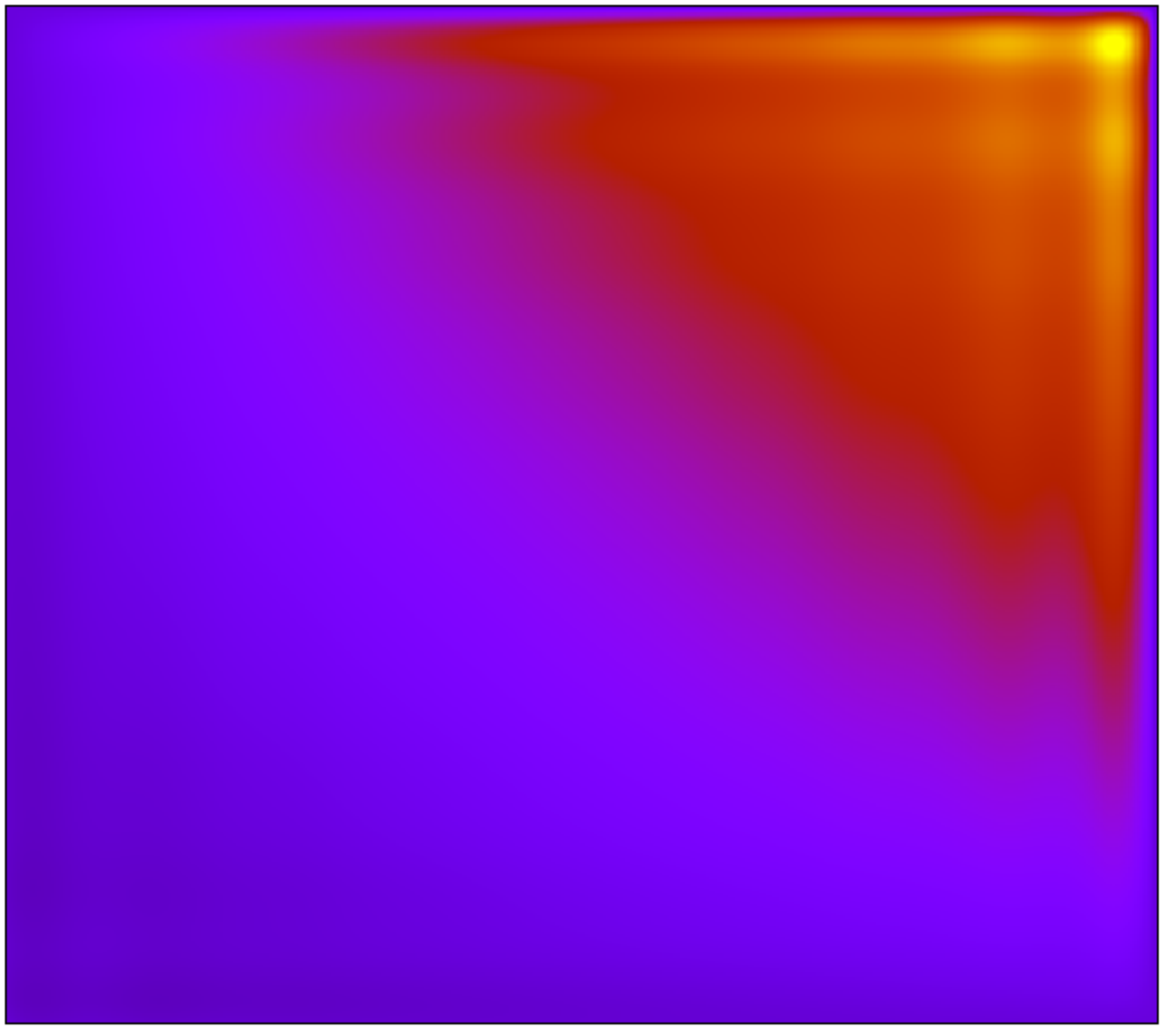}
&
\includegraphics[scale=0.16]{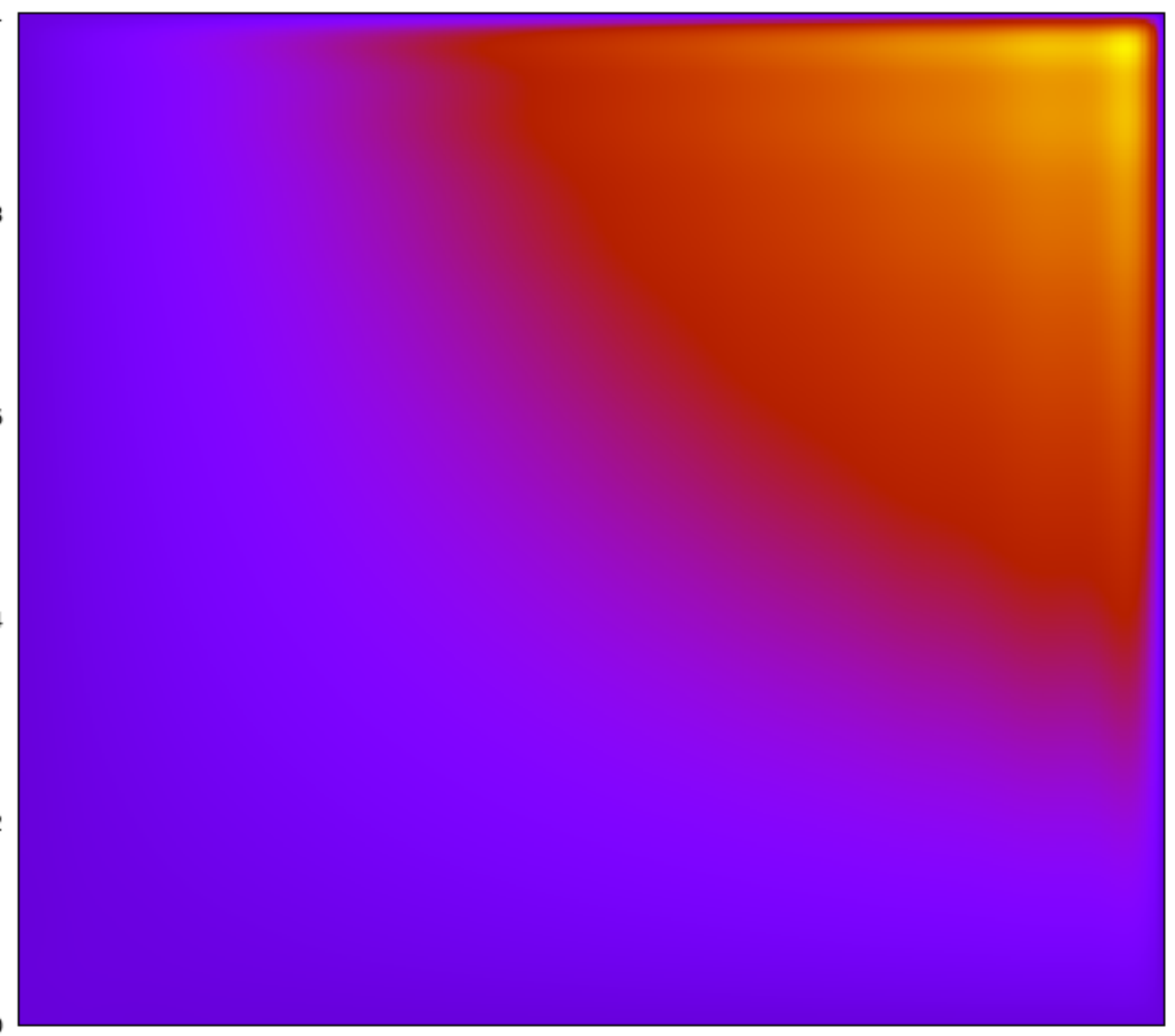}
&
\includegraphics[scale=0.16]{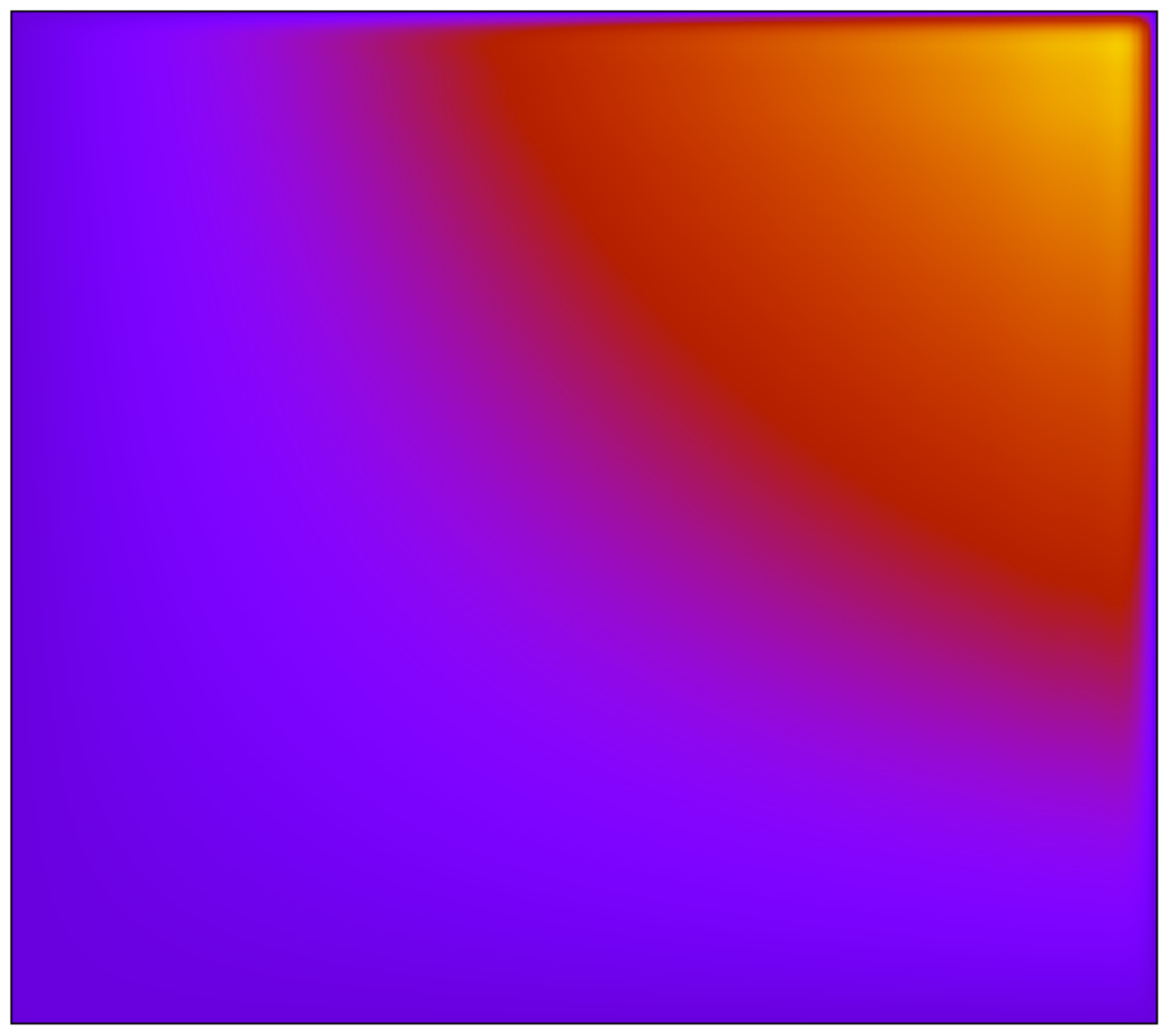}
\\
\hline
\#DOF & 5381 & 5450 & 5521 & 5594 \\
L2 & 32.07 & 1.33 & 0.27 & 0.056 \\
H1 & 31.01 & 9.77 & 3.16 & 0.82
\\
$32\times 32$ & 
\includegraphics[scale=0.16]{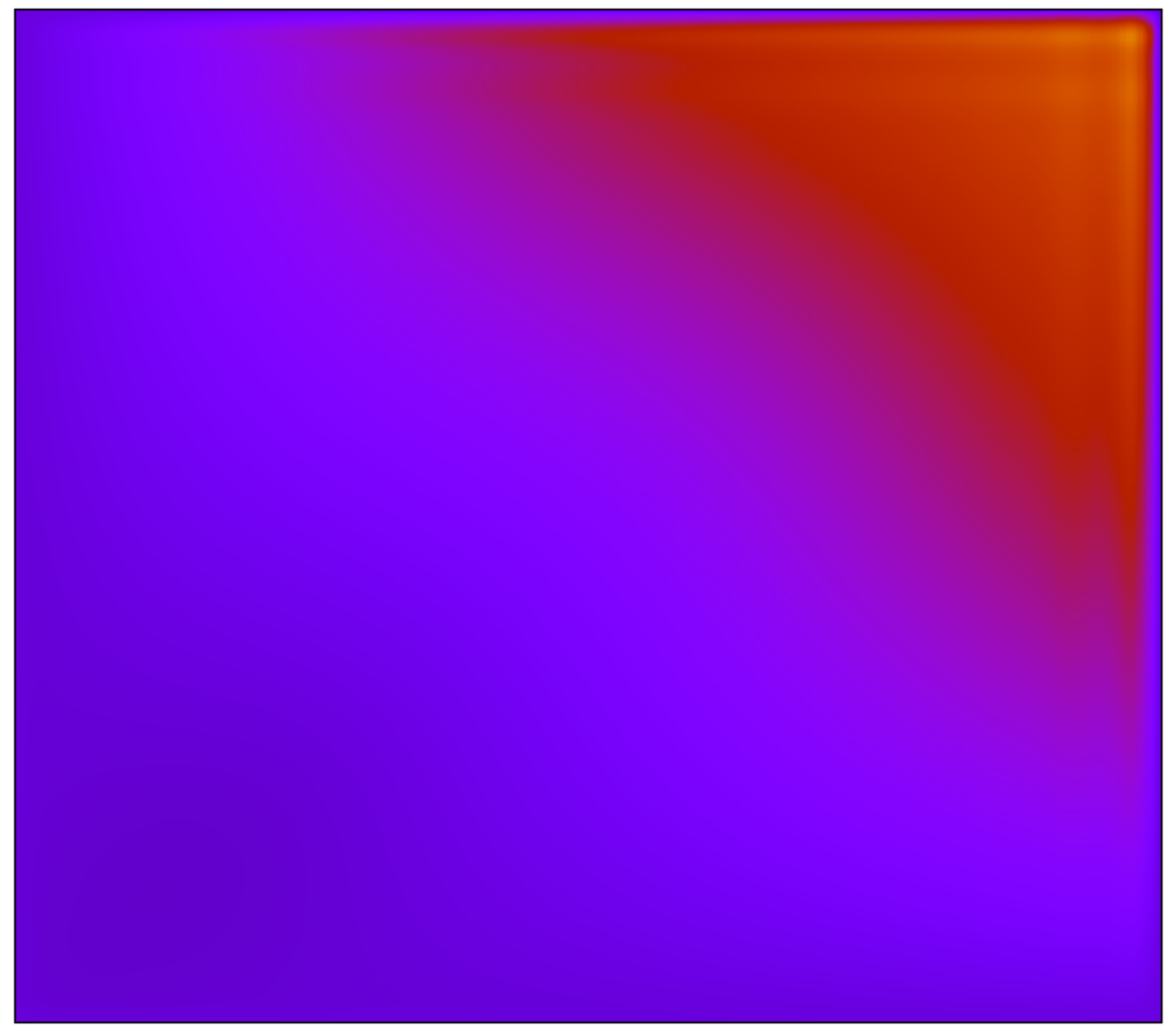}
&
\includegraphics[scale=0.16]{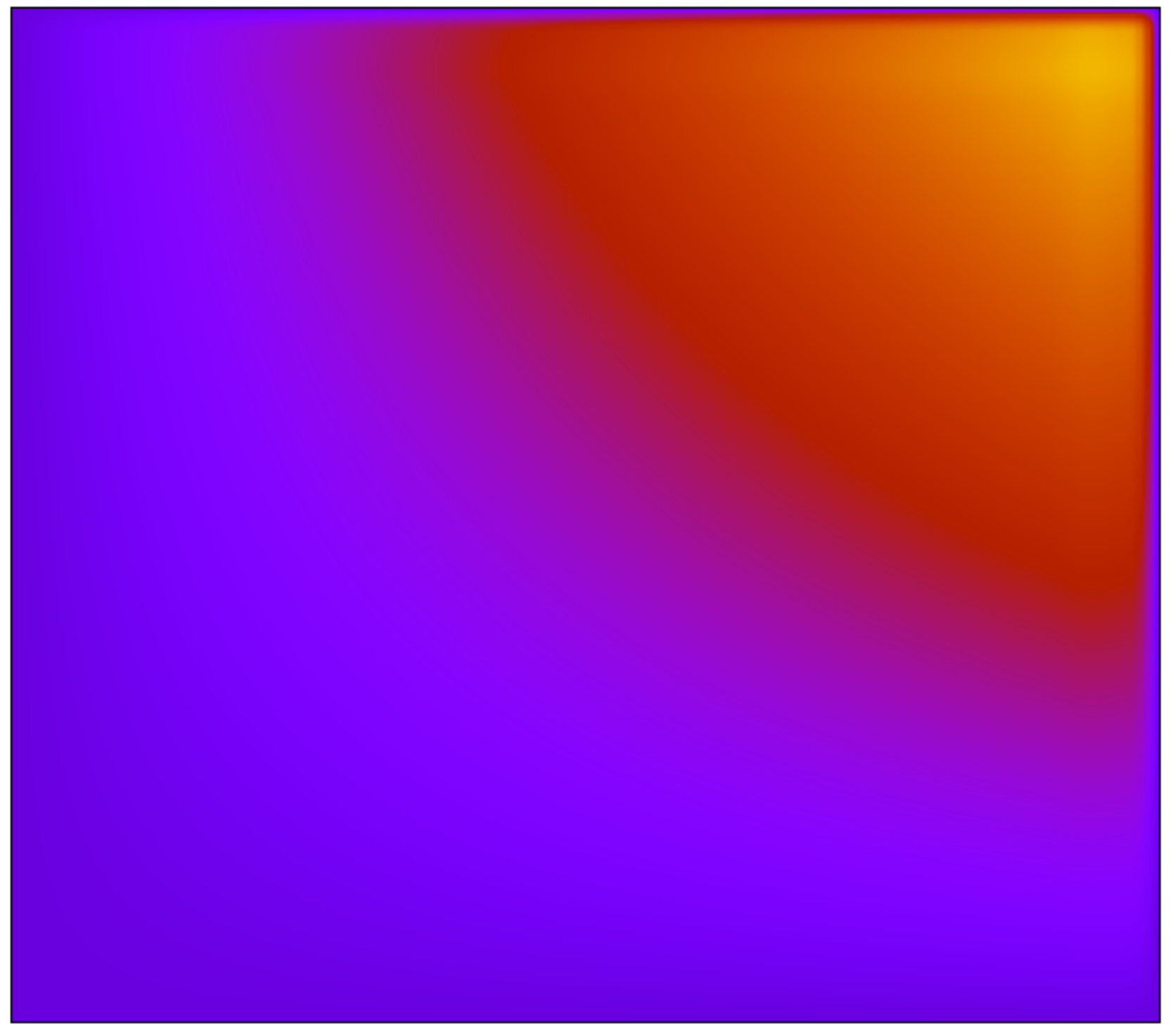}
&
\includegraphics[scale=0.16]{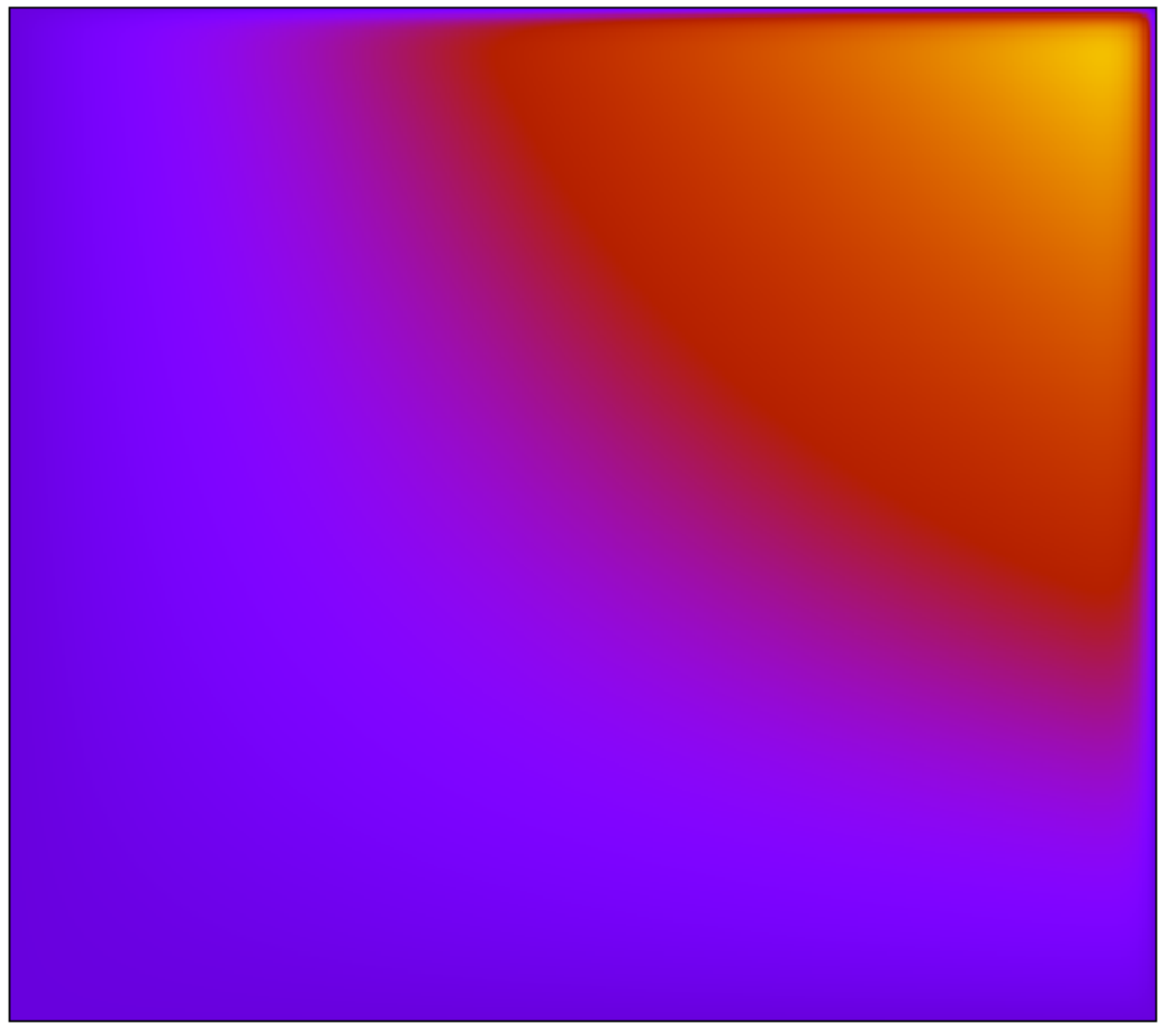}
&
\includegraphics[scale=0.16]{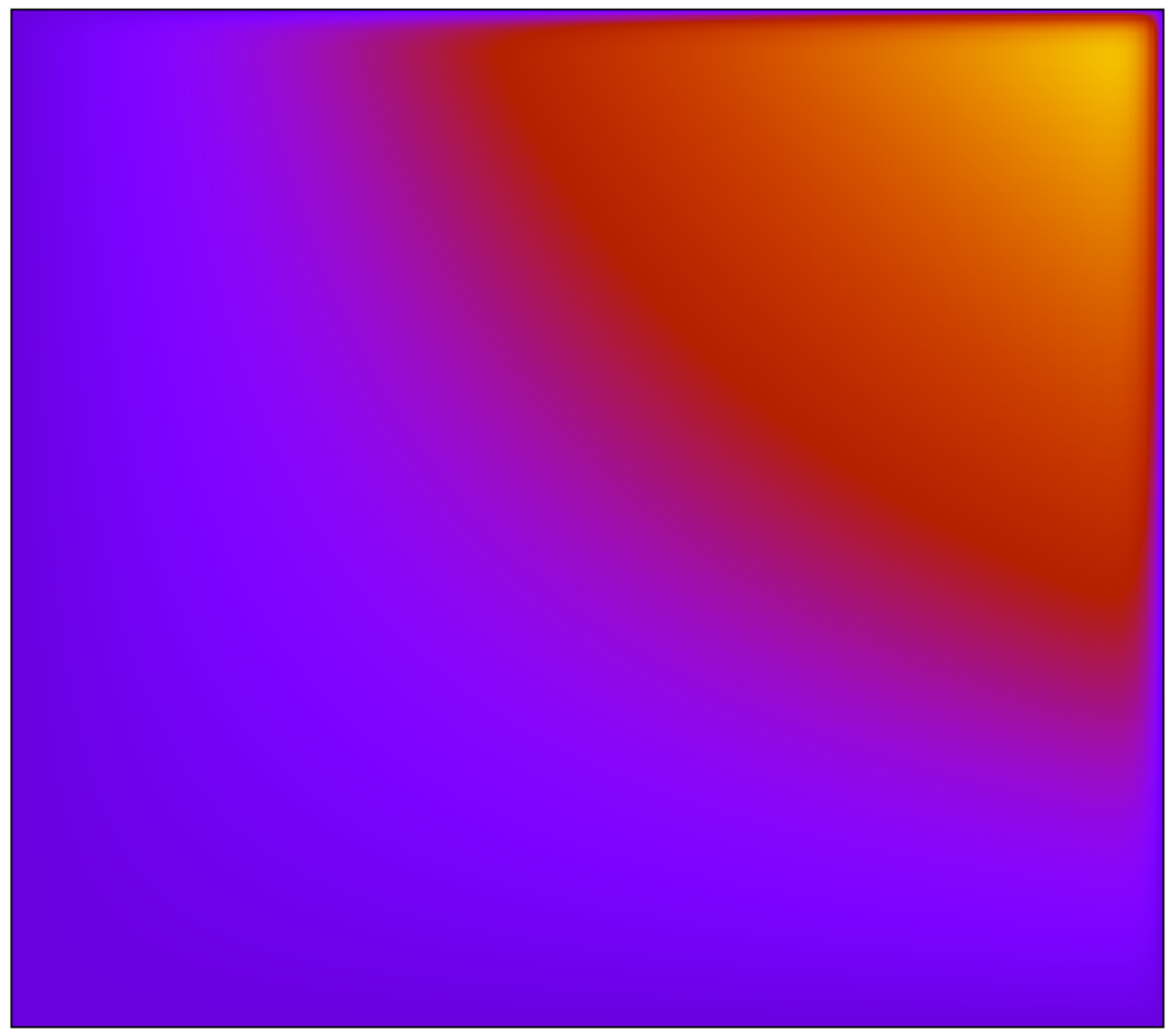}
\\
\hline
\#DOF & 20997& 21130& 21265 & 21402 \\
L2 & 7.66 & 0.07 & 0.01 & 0.003 \\
H1 & 9.86 & 1.67 & 0.26 & 0.068
\\
$64\times 64$ & 
\includegraphics[scale=0.16]{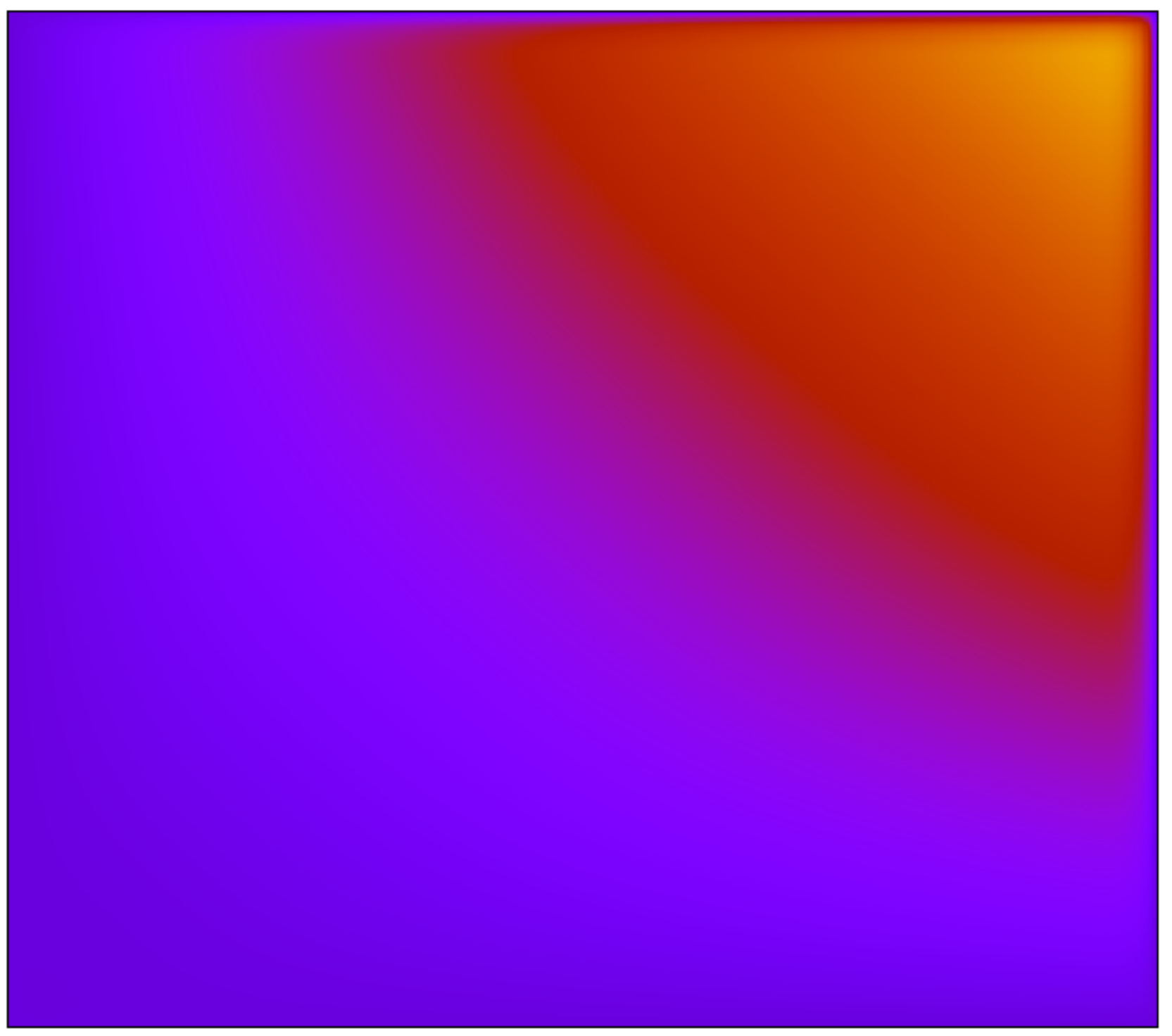}
&
\includegraphics[scale=0.16]{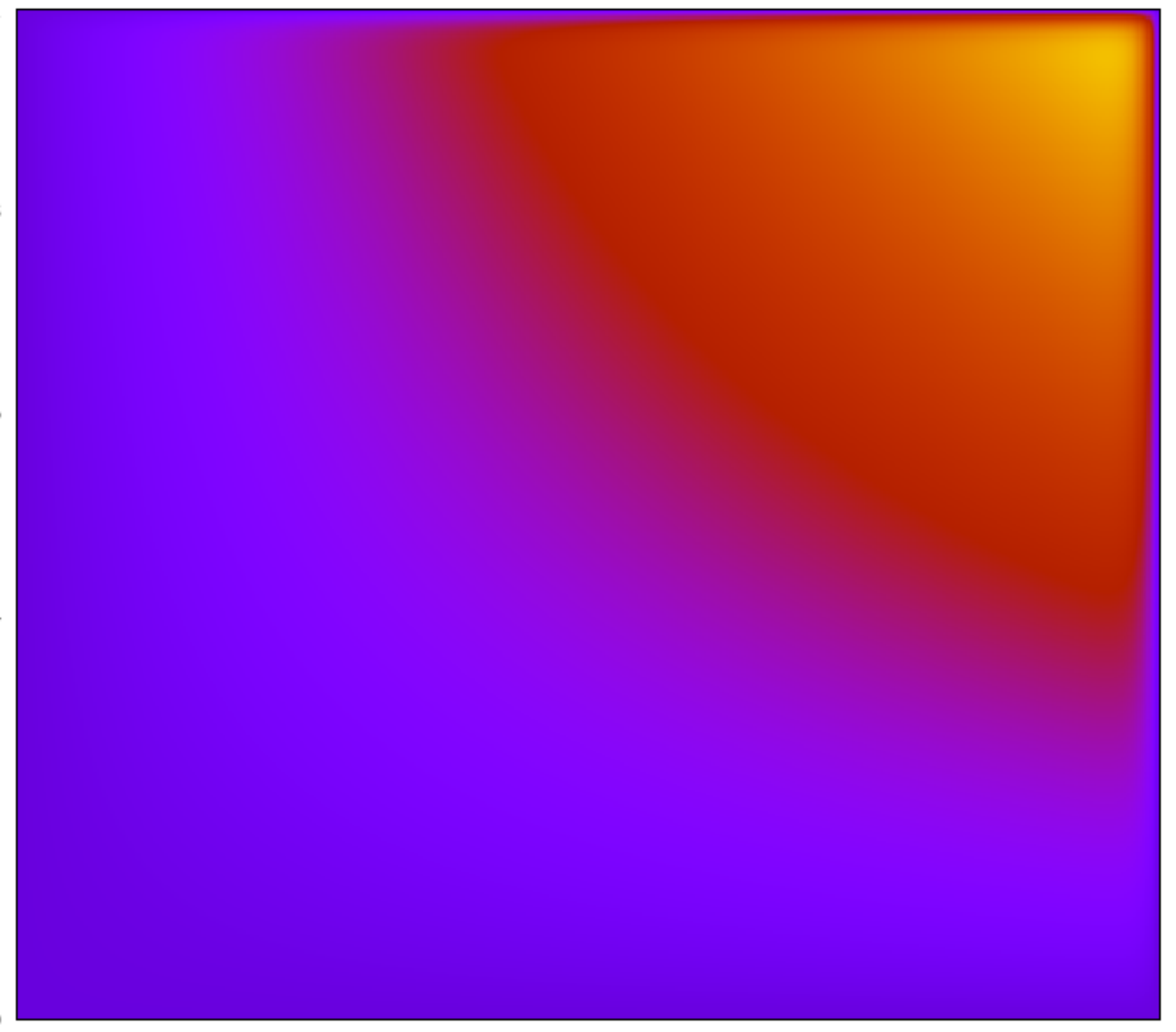}
&
\includegraphics[scale=0.16]{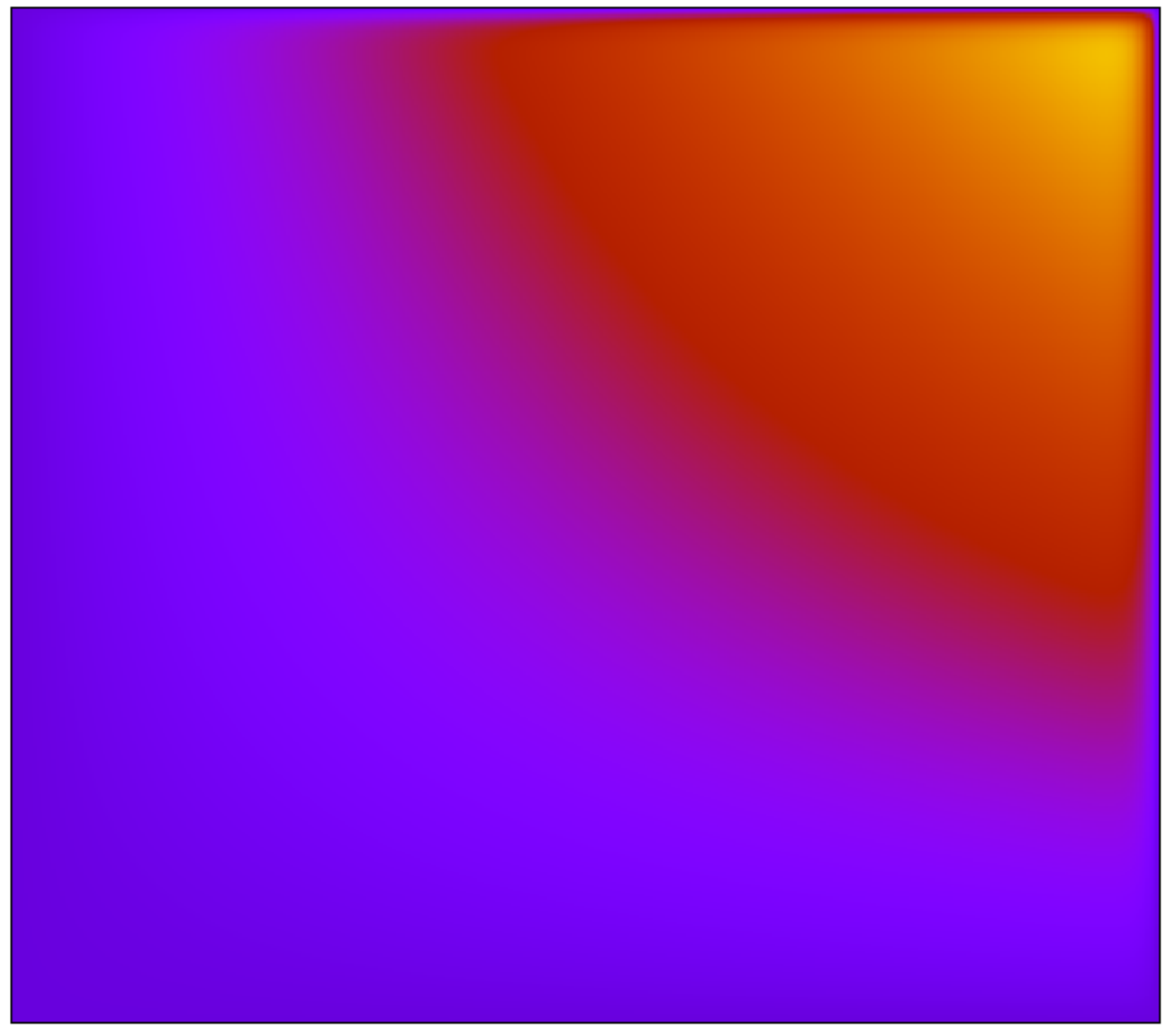}
&
\includegraphics[scale=0.16]{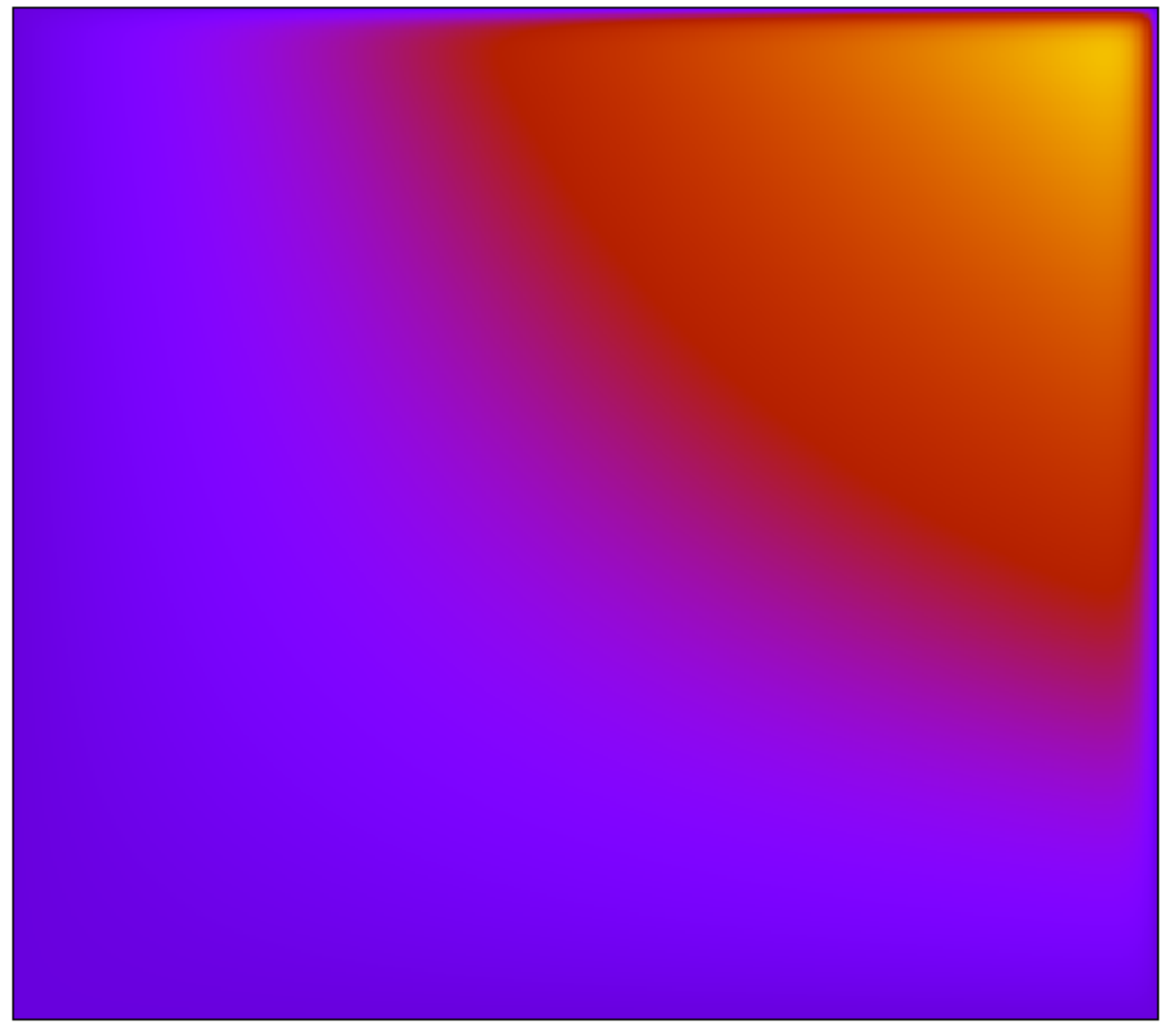}
\\
\hline
\end{tabular}
\includegraphics[scale=0.4]{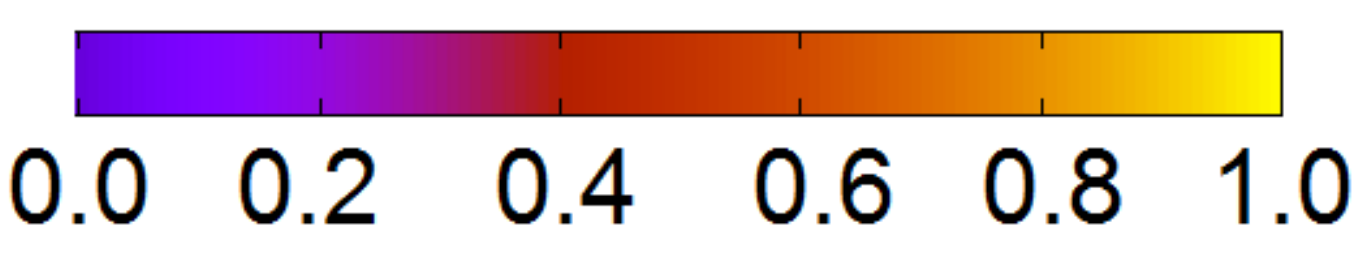}
\end{center}
\caption{Solution of Problem 1 by iGRM method, with different solution and residual spaces, for different mesh dimensions.}
\label{tab:Problem1}
\end{table*}}

{\tiny\begin{table*}[htp]
\begin{center}
\begin{tabular}{|c|c|c|c|c|}
\hline 
n & trial (2,1)  & trial (3,2)  & trial (4,3) & trial (5,4) \\
 & test (2,1) & test (3,2) & test (4,3)  & test (5,4)  \\
\hline
\#DOF & 100 & 121 & 144 & 169 \\
L2 & 40.36 & 38.74 & 38.72 & 38.82 \\
H1 & 73.30 & 80.16 & 82.86 & 83.72 \\
$8\times 8$ & 
\includegraphics[scale=0.15]{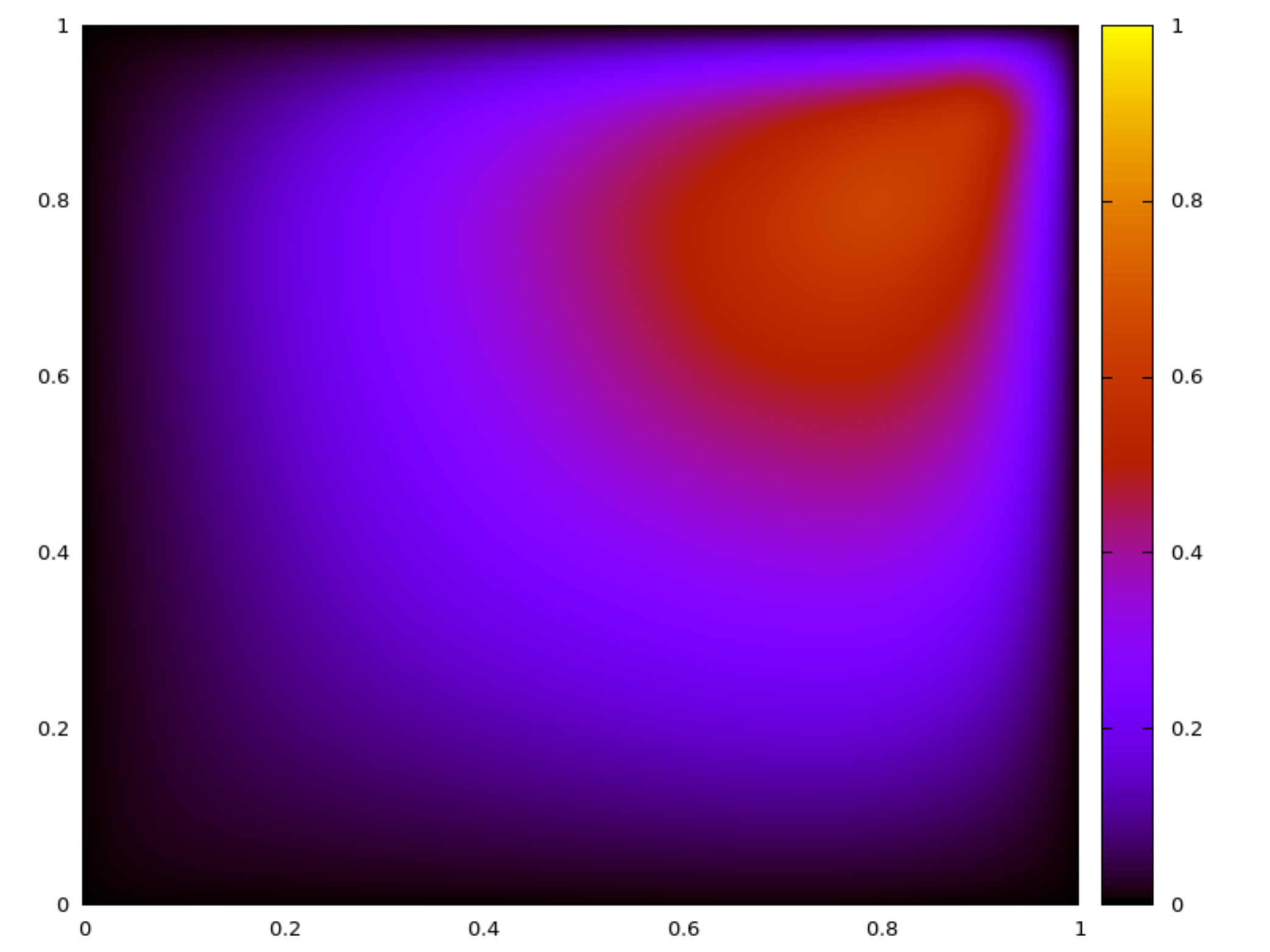}
&
\includegraphics[scale=0.15]{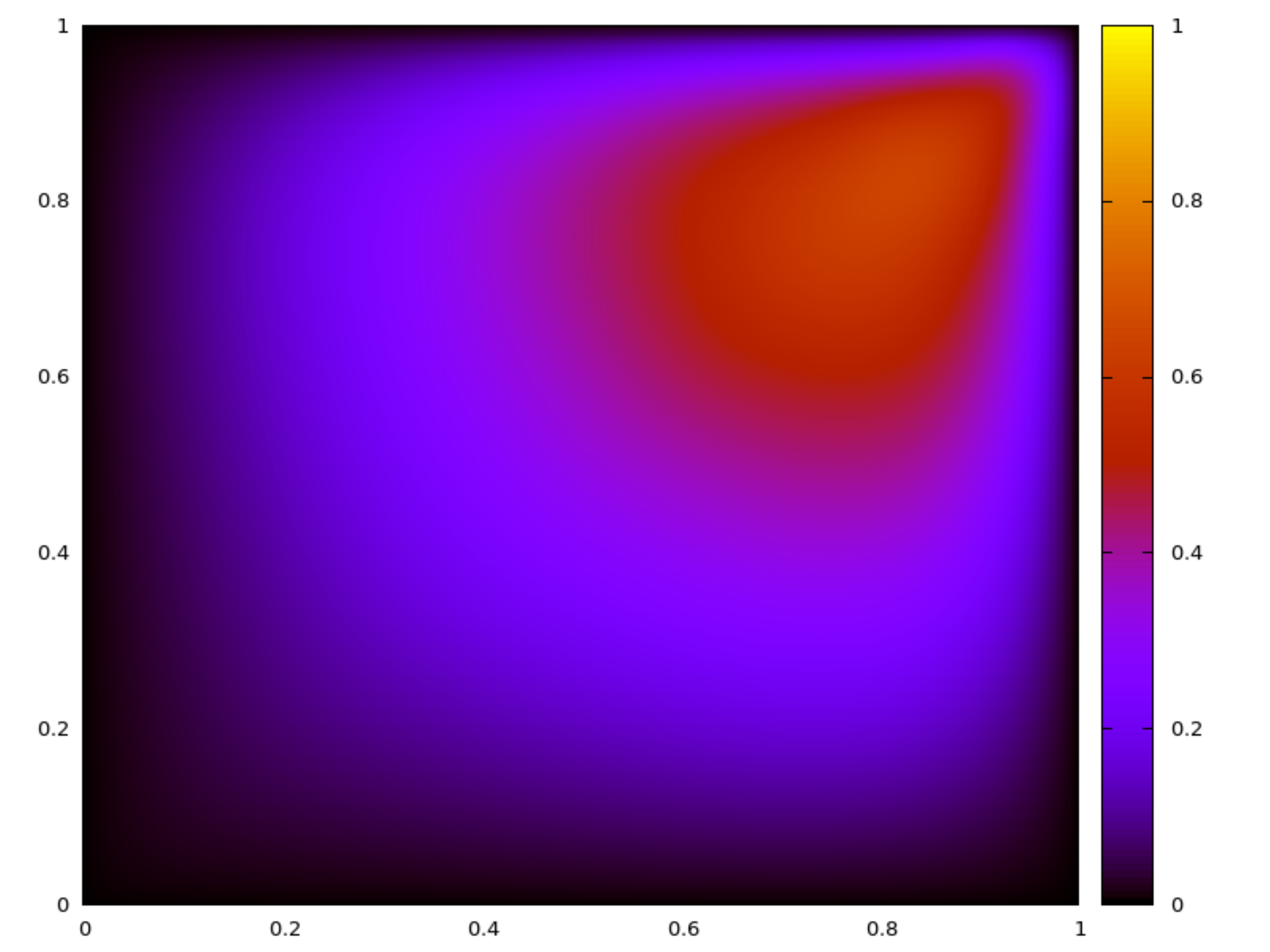}
&
\includegraphics[scale=0.15]{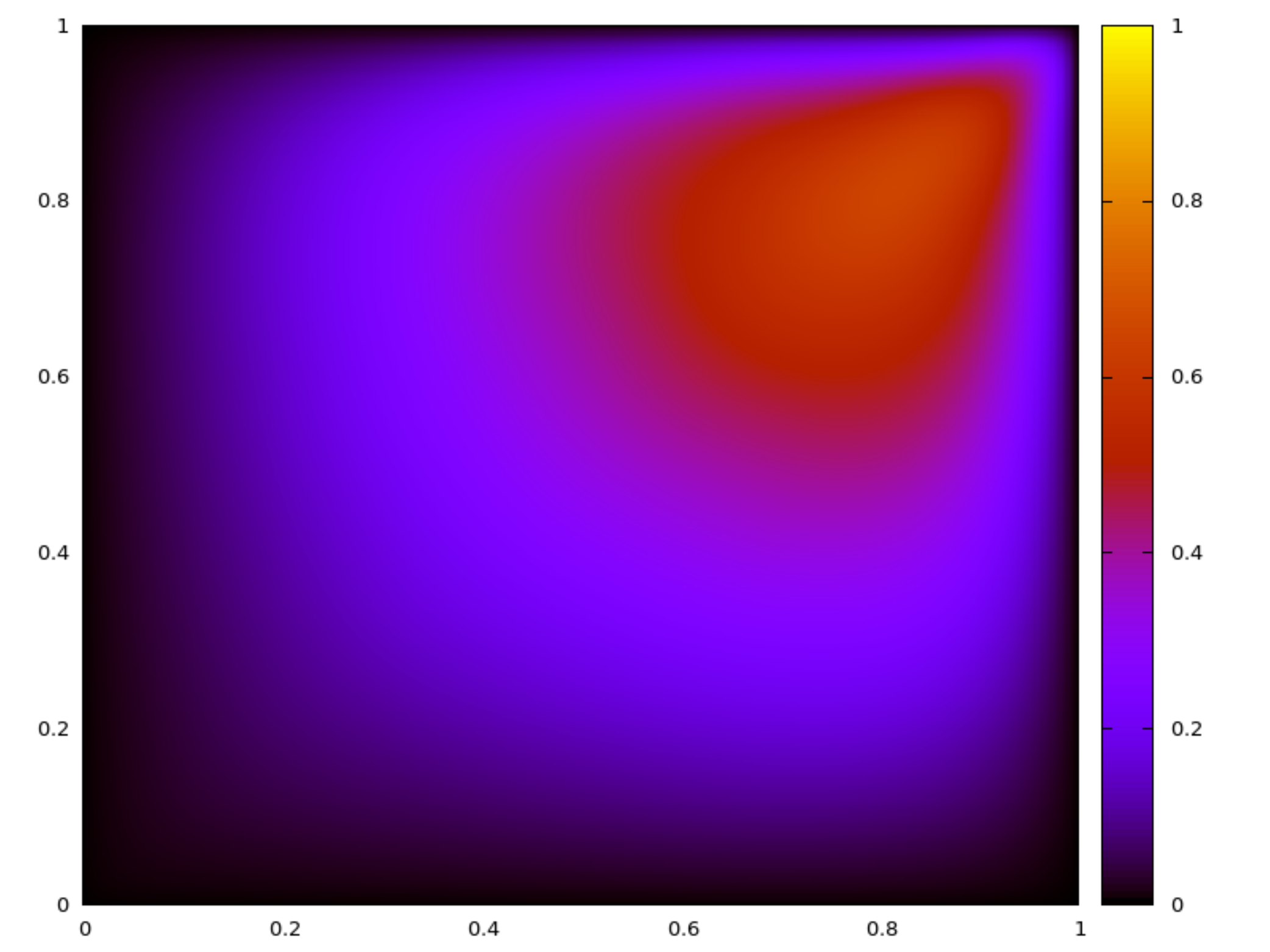}
&
\includegraphics[scale=0.15]{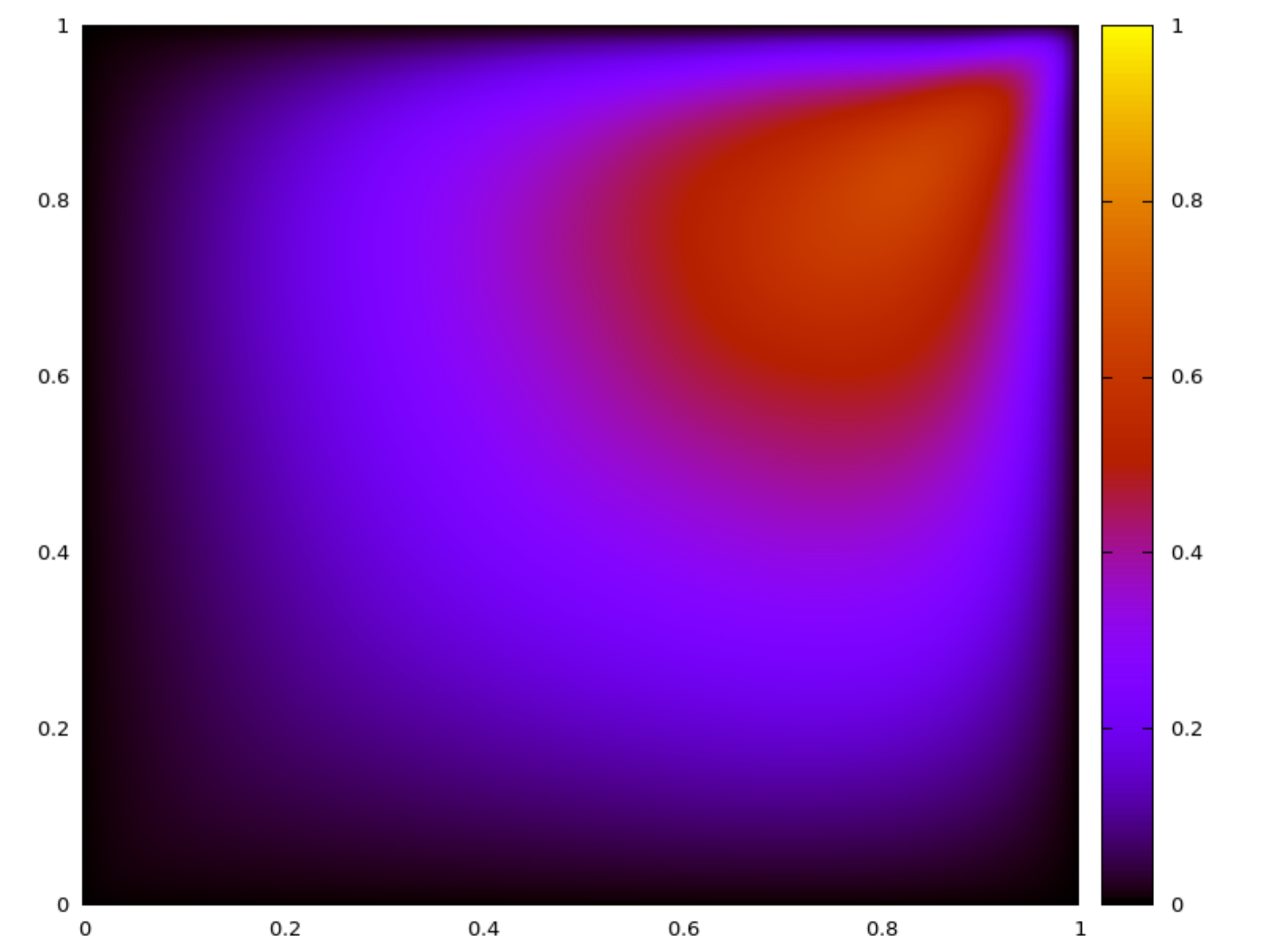}
\\
\hline
\#DOF & 324 & 361 & 400  & 441 \\
L2 & 20.04 & 19.00 & 18.61 & 18.52 \\
H1 & 63.04 & 62.19 & 61.24 & 61.01 \\
$16\times 16$ & 
\includegraphics[scale=0.15]{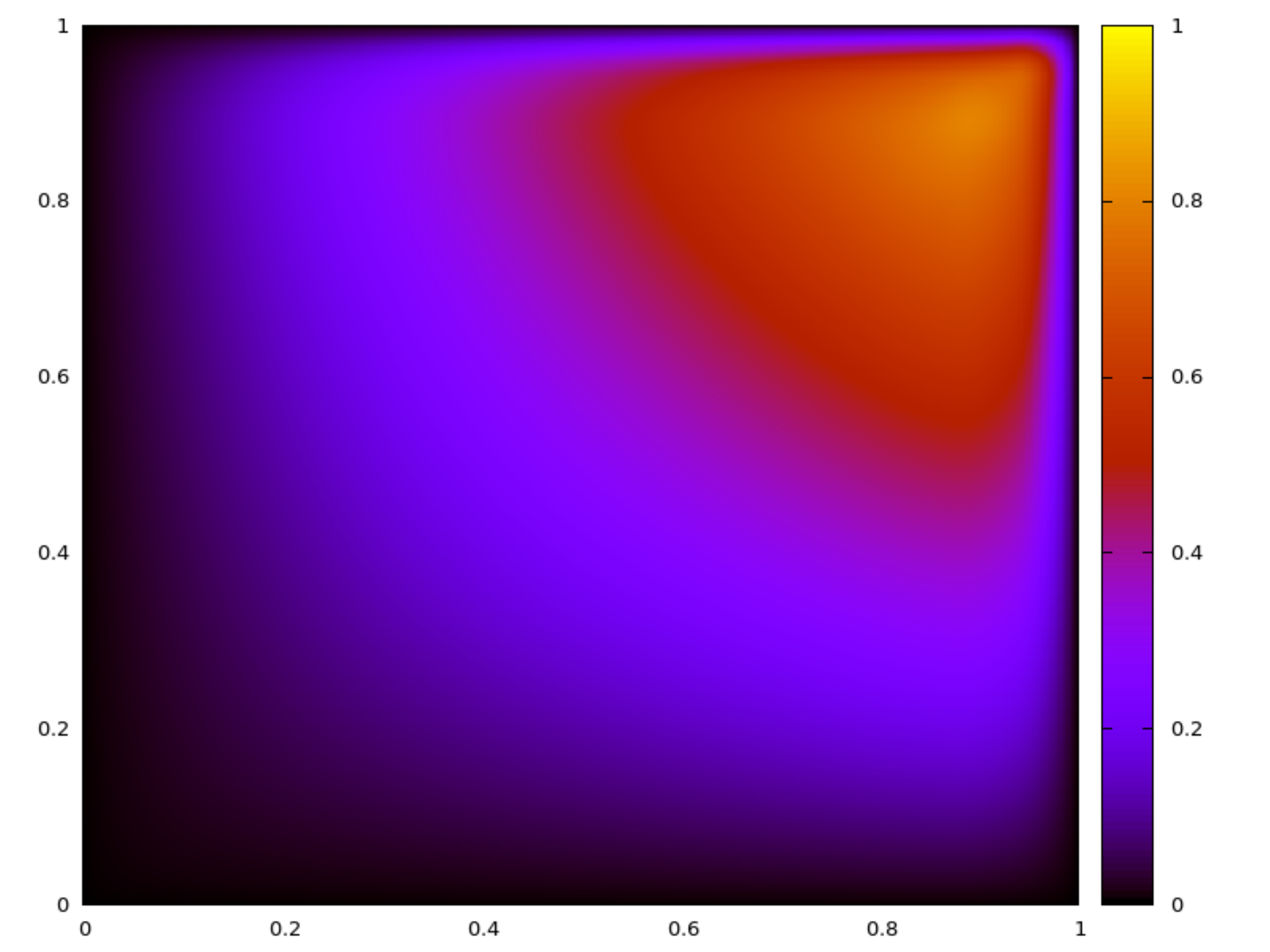}
&
\includegraphics[scale=0.15]{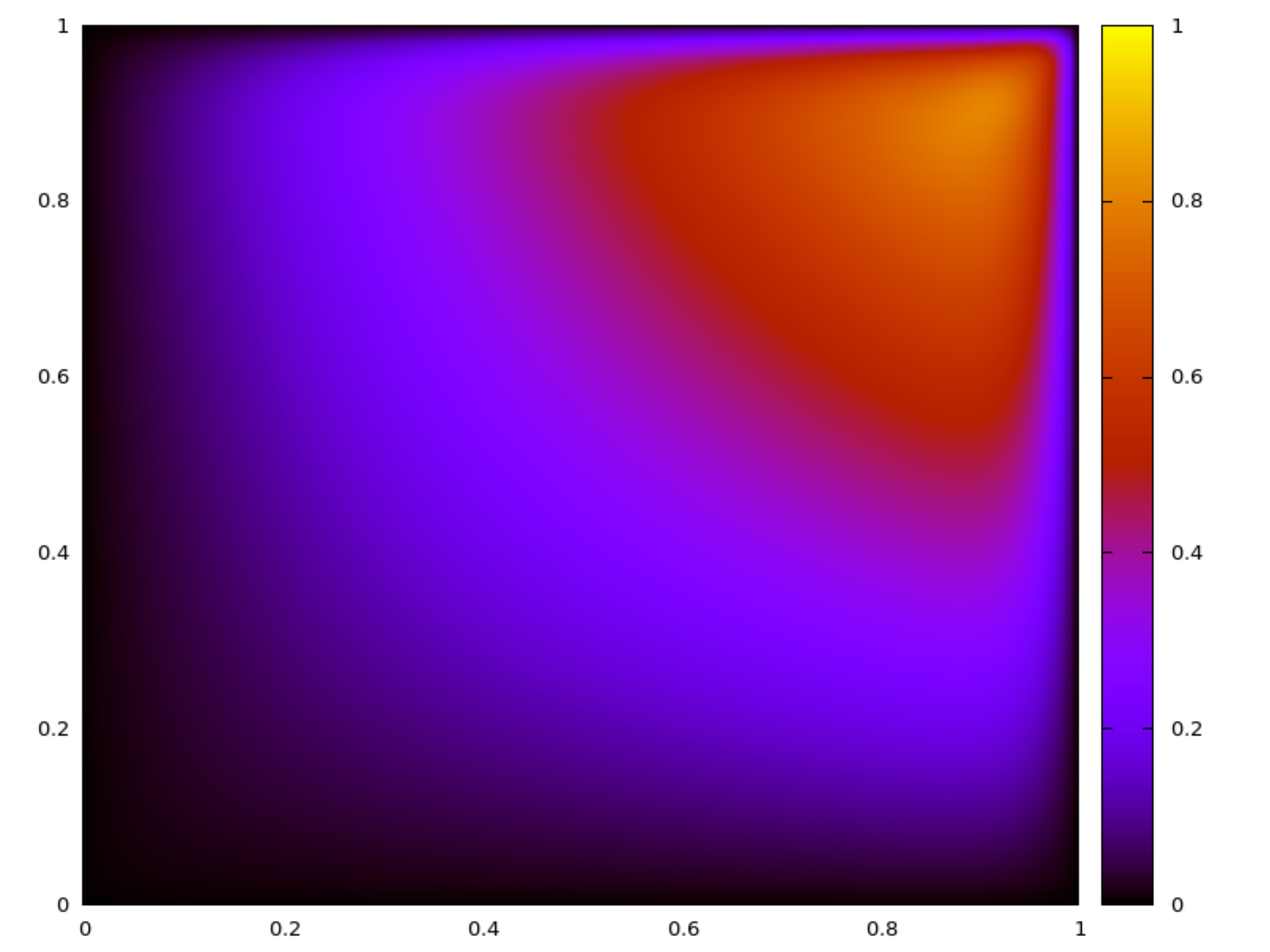}
&
\includegraphics[scale=0.15]{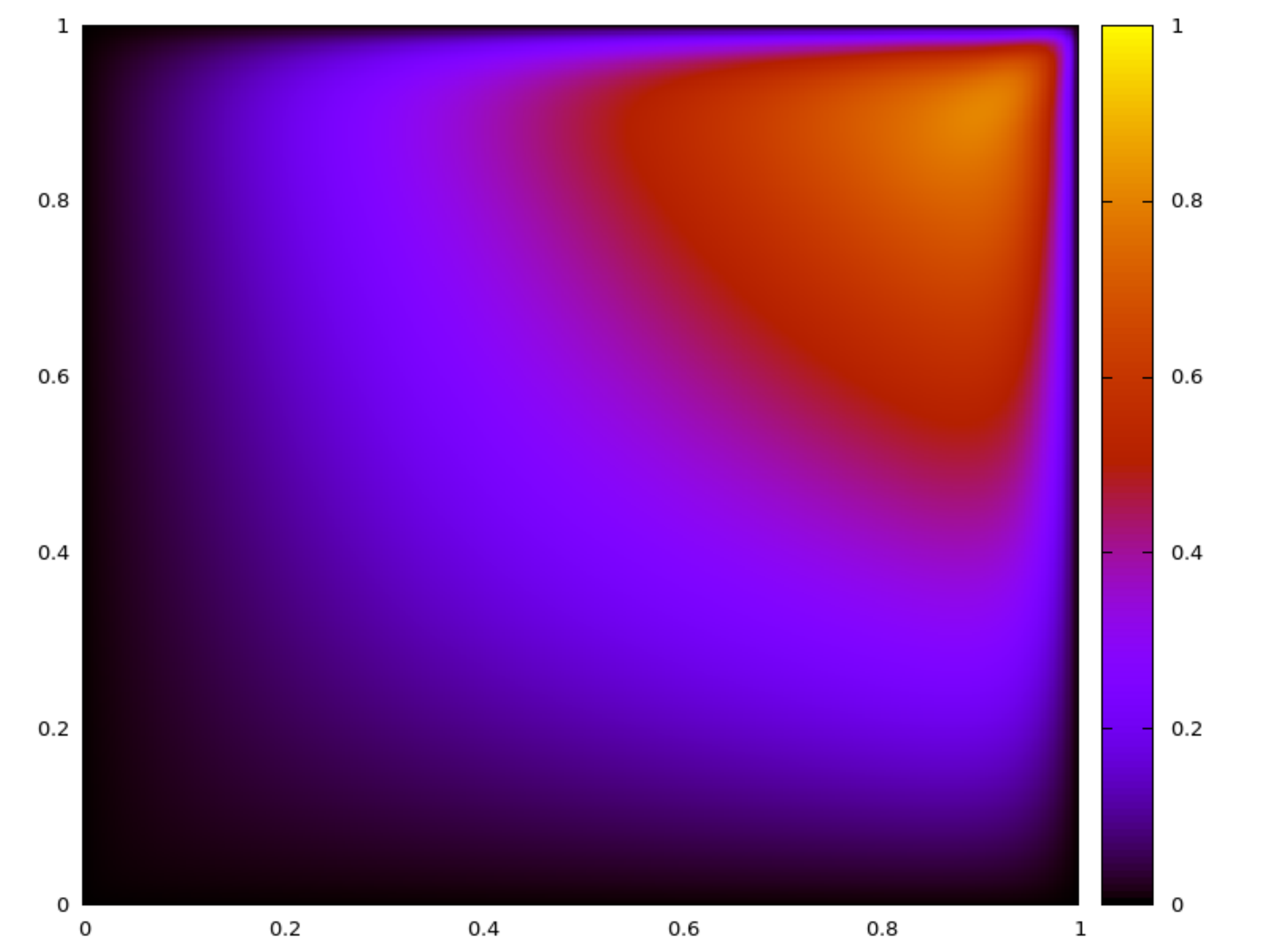}
&
\includegraphics[scale=0.15]{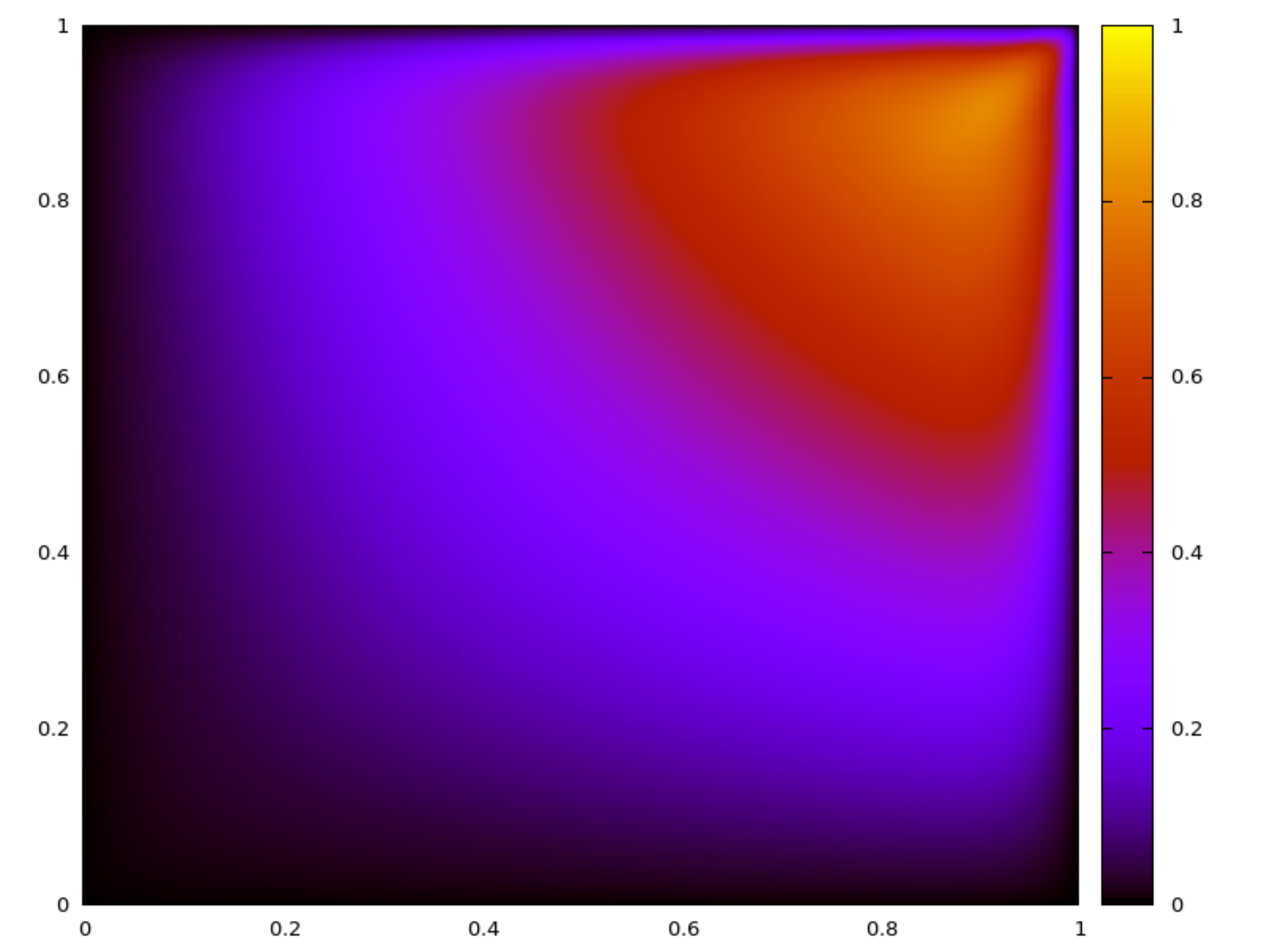}
\\
\hline
\#DOF & 1156 & 1225 & 1296 & 1369 \\
L2 & 6.51 & 5.07 & 4.20 & 3.66 \\
H1 & 36.57 & 26.88 & 21.16 & 17.35 \\
$32\times 32$ & 
\includegraphics[scale=0.15]{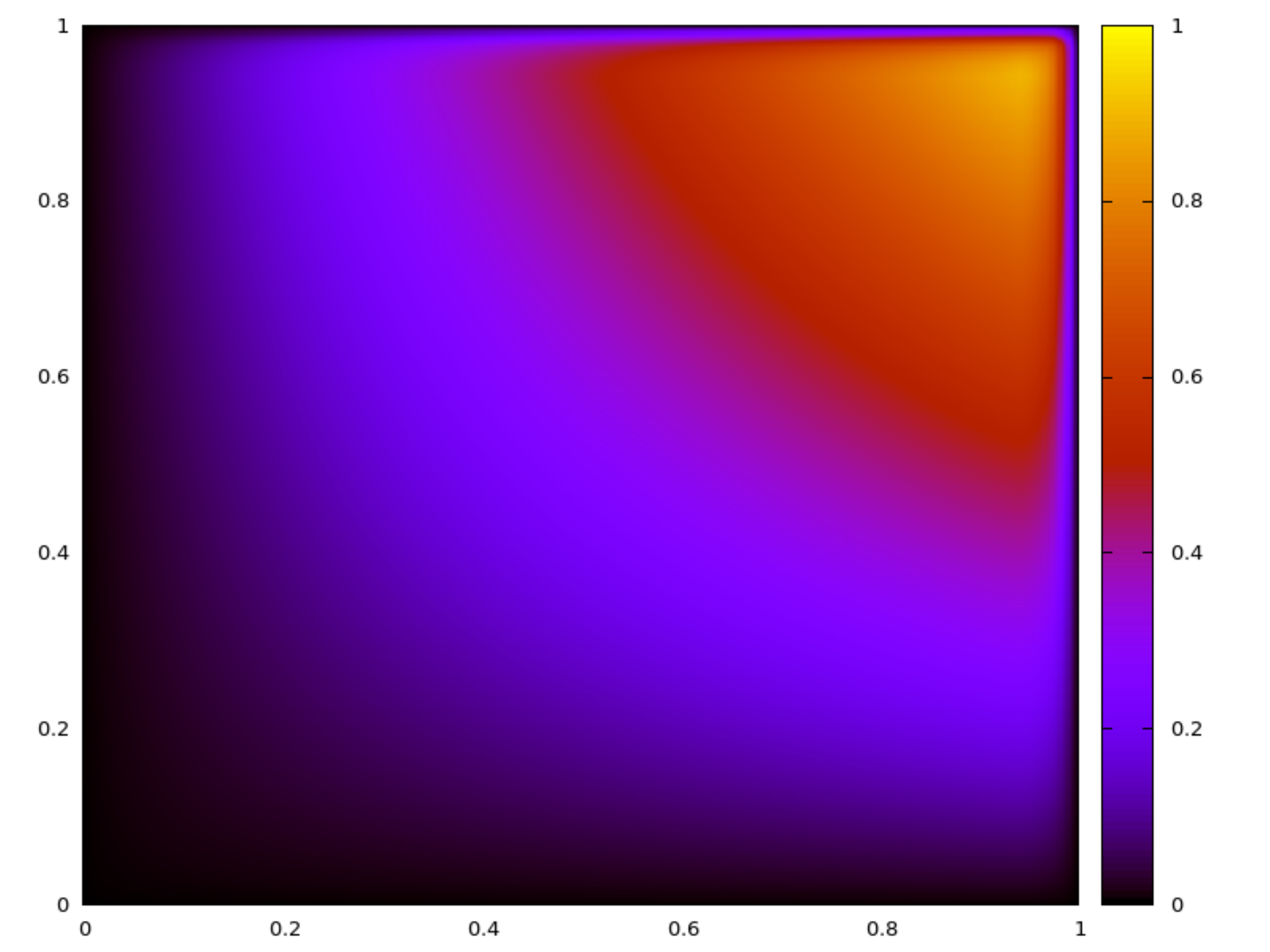}
&
\includegraphics[scale=0.15]{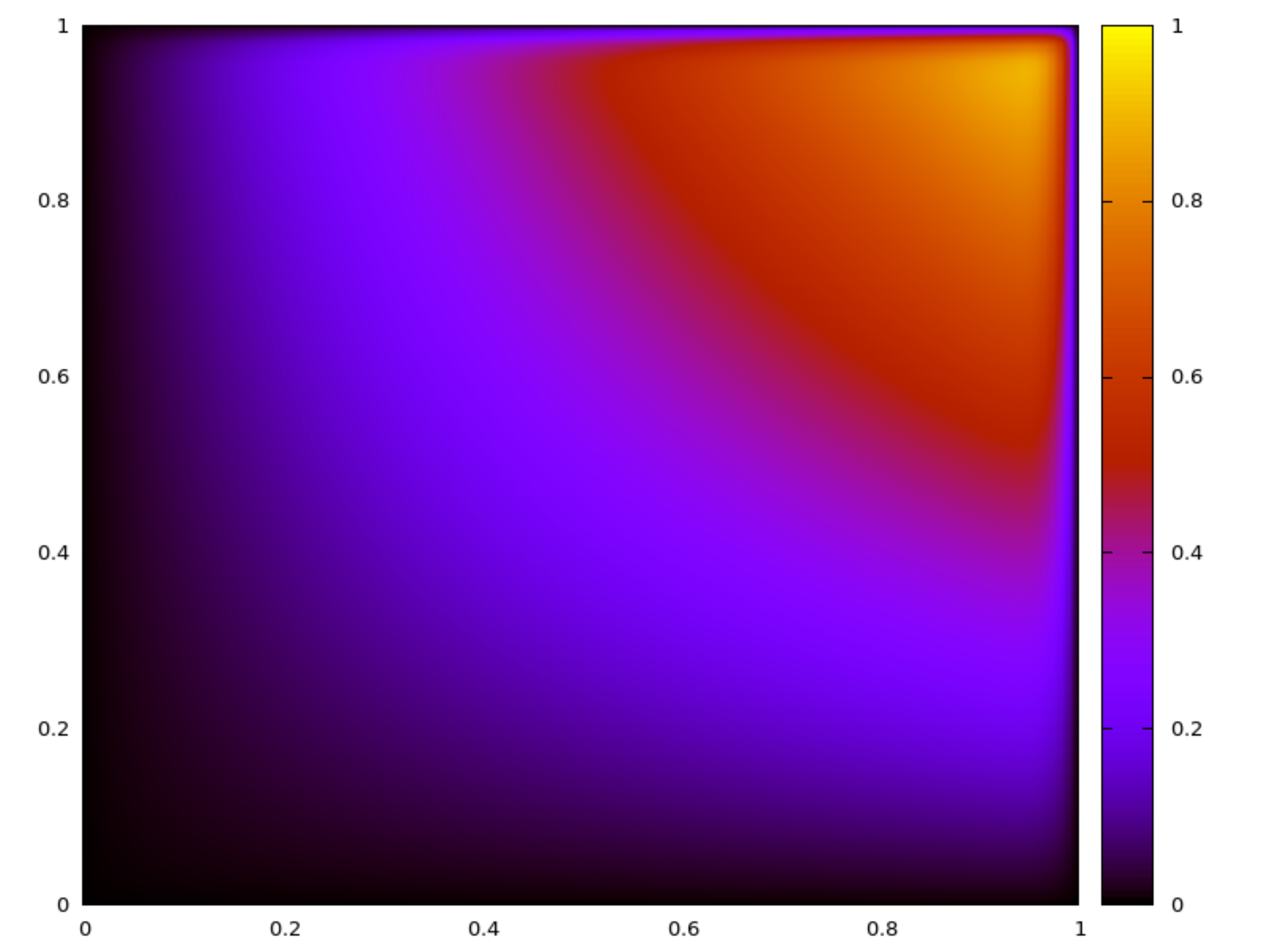}
&
\includegraphics[scale=0.15]{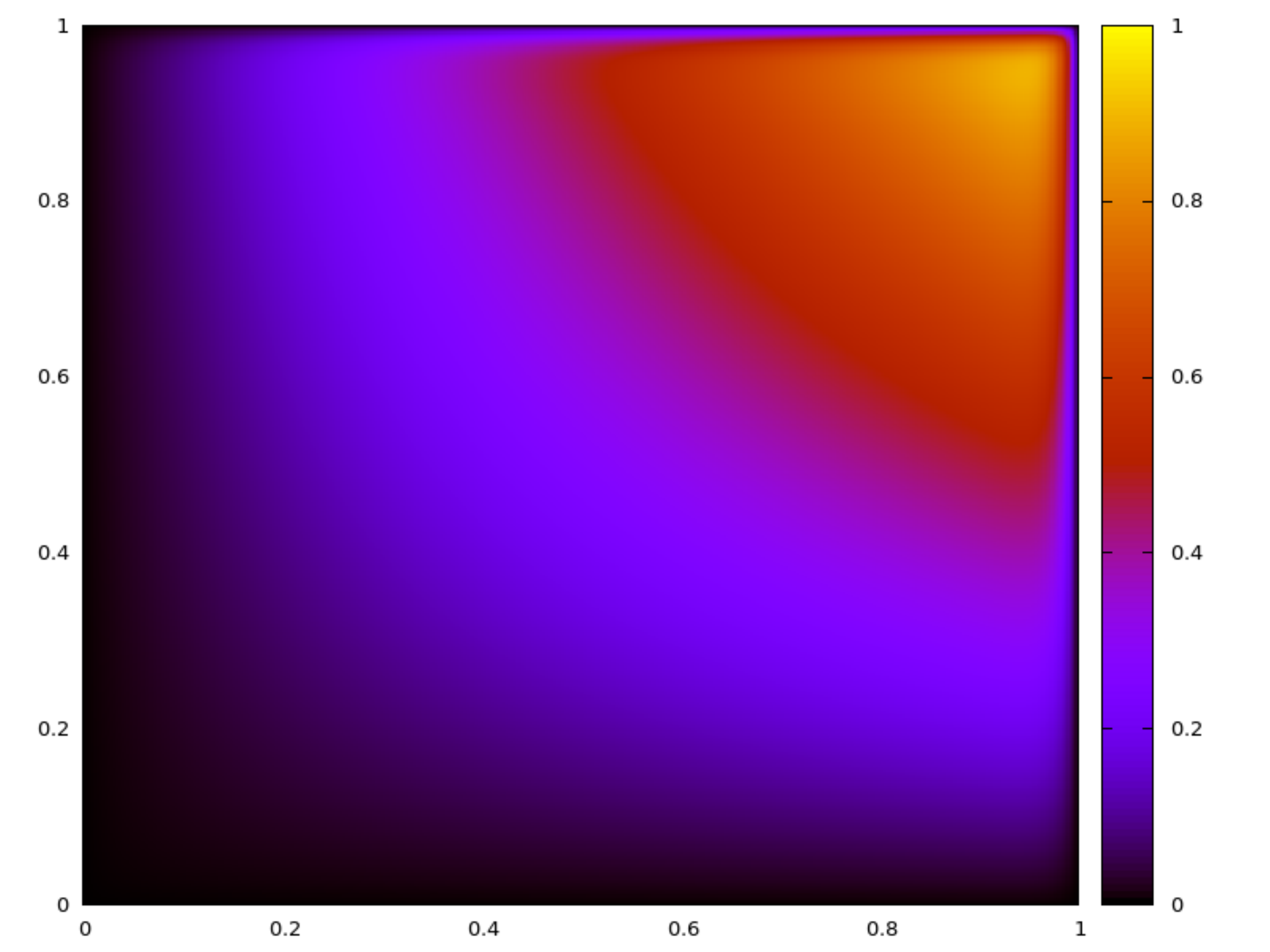}
&
\includegraphics[scale=0.15]{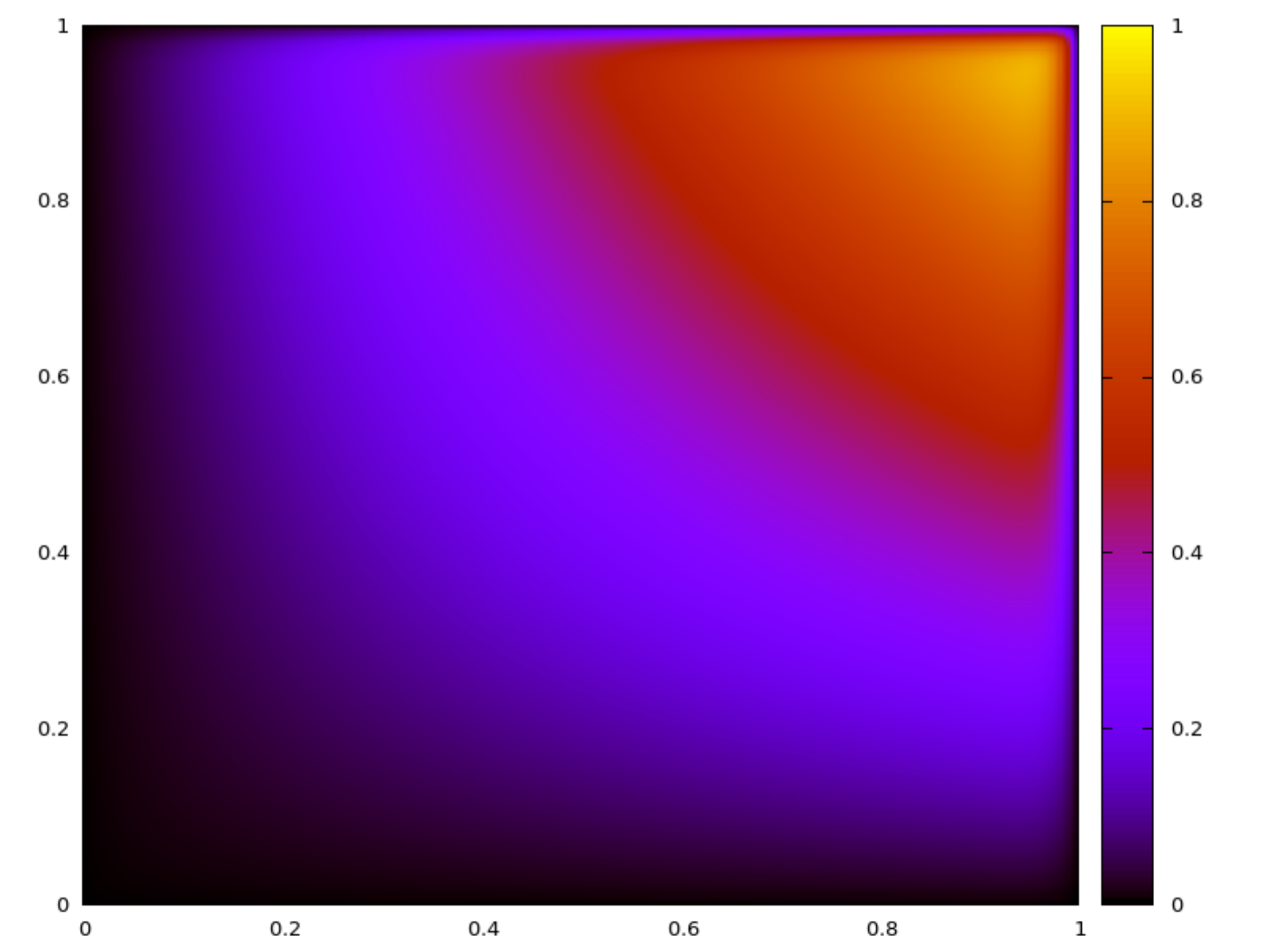}
\\
\hline
\#DOF & 4356 & 4489 & 4624 & 4761 \\
L2 & 1.19 & 0.80 & 0.75 & 0.75 \\
H1 & 11.50 & 3.94 & 1.39 & 0.62 \\
$64\times 64$ & 
\includegraphics[scale=0.15]{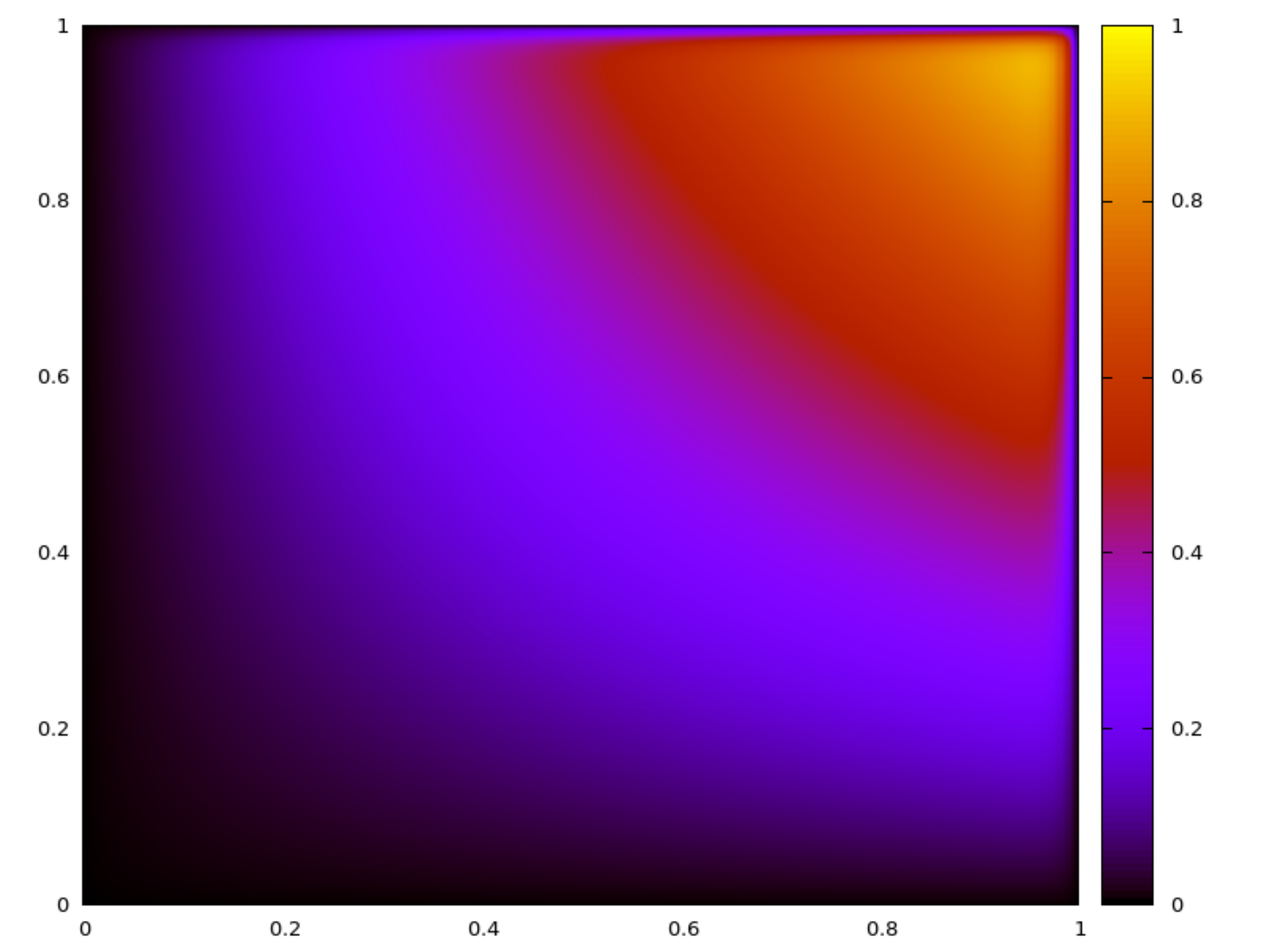}
&
\includegraphics[scale=0.15]{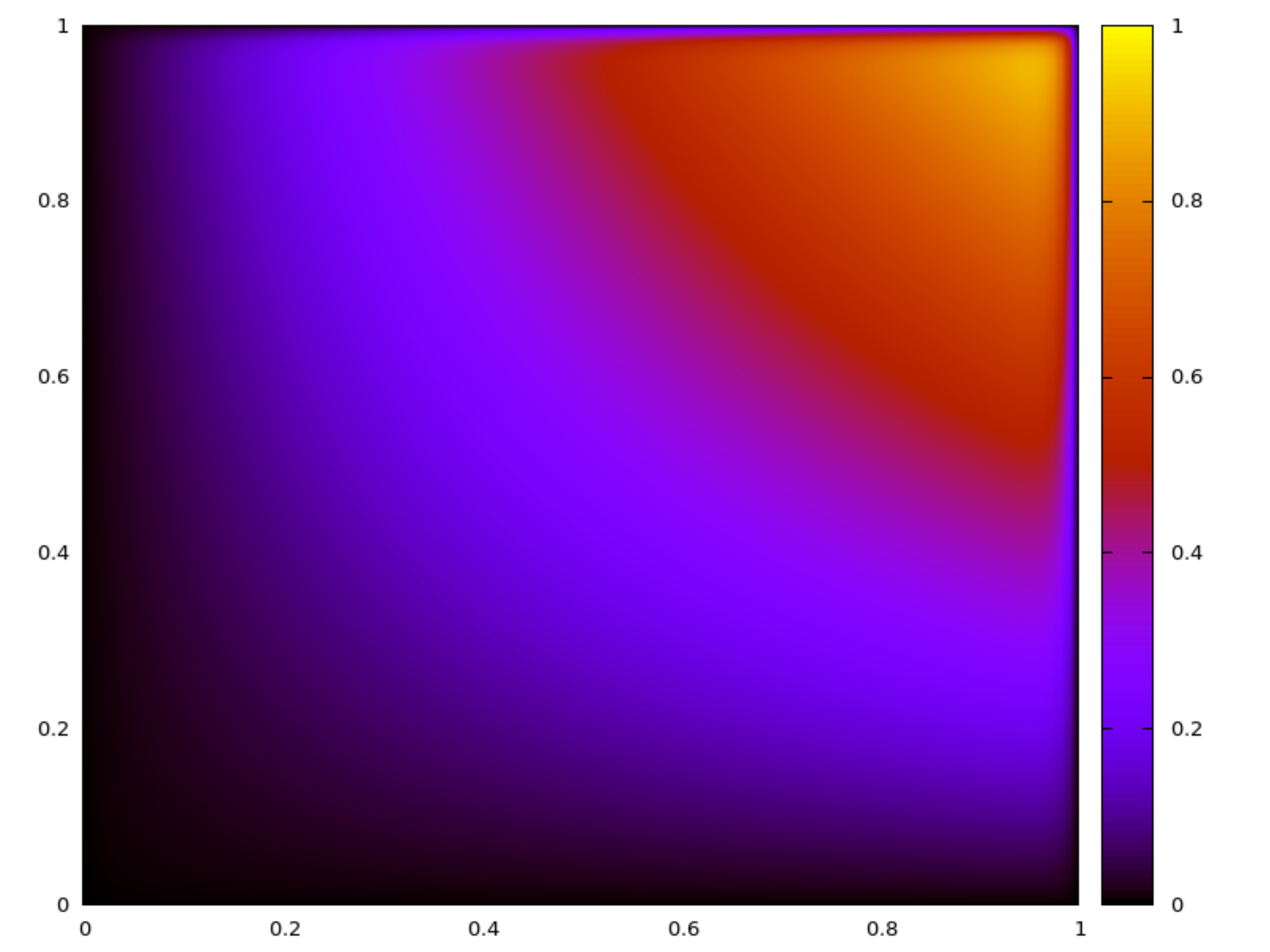}
&
\includegraphics[scale=0.15]{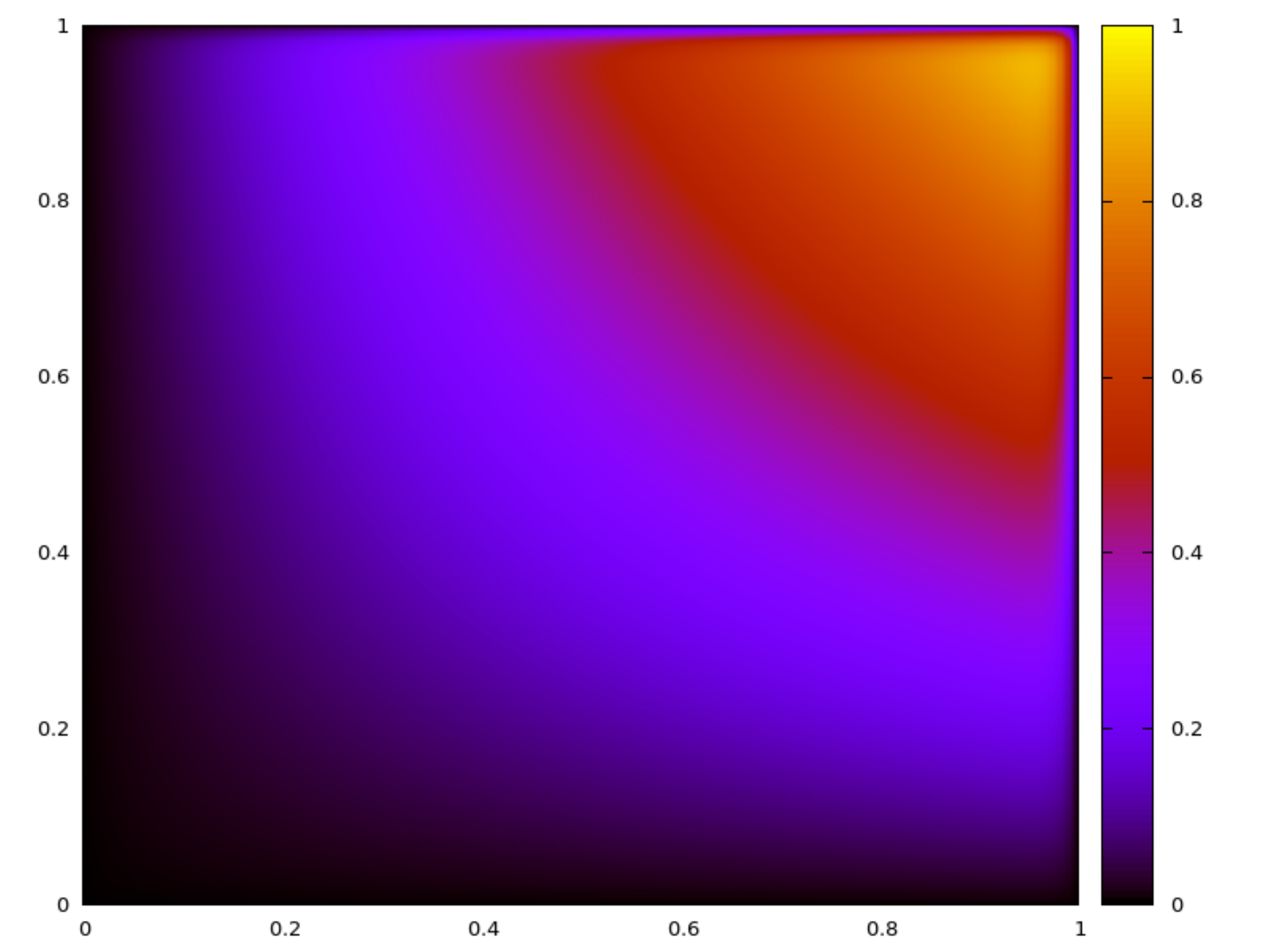}
&
\includegraphics[scale=0.15]{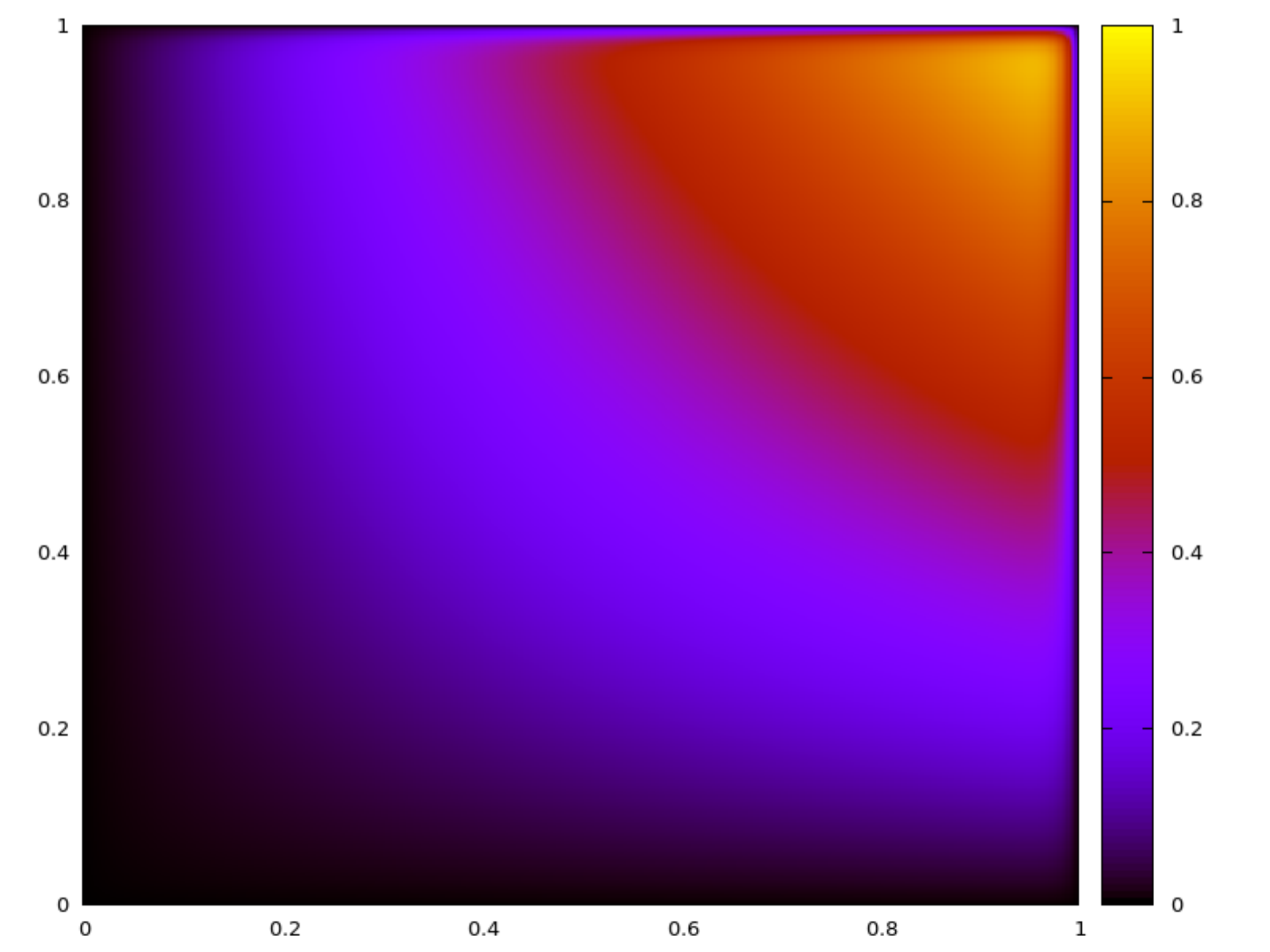}
\\
\hline
\end{tabular}
\includegraphics[scale=0.4]{legend.pdf}
\end{center}
\caption{\revcc{Solution of Problem 1 by SUPG method, with different solution and residual spaces, for different mesh dimensions.}}
\label{tab:Problem1SUPG}
\end{table*}}

\subsection{Erikkson-Johnson model problem}

We focus now on the model Eriksson-Johnson problem with the modifications proposed by~\cite{ Erikkson}. For the square domain $\Omega=(0,1)^2$ and the advection vector $\beta=(1,0)^T$, we seek the solution of the advection-diffusion equation~\eqref{eq:ModelProblem}.  We introduce the Dirichlet boundary conditions
$$
u=g=sin(\pi y) \textrm{ for }x\in \Gamma^{-}
$$
$$
u=0 \textrm{ for } x\in \Gamma^{+}
$$
weakly on the boundary $\Gamma$.  The inflow Dirichlet boundary condition drives the problem and it develops a boundary layer of width $\epsilon$ at the outflow $x = 1$.


\revcc{In the manufactured solution problem, we compare the SUPG and iGRM methods on a sequence of uniform grids. \revMP{Later}, we compare iGRM and SUPG methods on a sequence of adapted grids.  We simulate the Erikkson-Johnson problem with the iGRM and SUPG methods with $Pe=10^6$ on a sequence of grids refined towards the boundary layer. We start from the uniform grid of $2 \times 2$ elements, and we add the knot points to the last interval on the right in the $x$ direction.  We refine the knot vector in the $x$ direction by breaking in half the rightmost element.  We continue to refine in the perpendicular direction after we capture the boundary layer when the element size in one direction is \revMP{below} $\epsilon$.}

\begin{figure}
\centering
\includegraphics[scale=0.45]{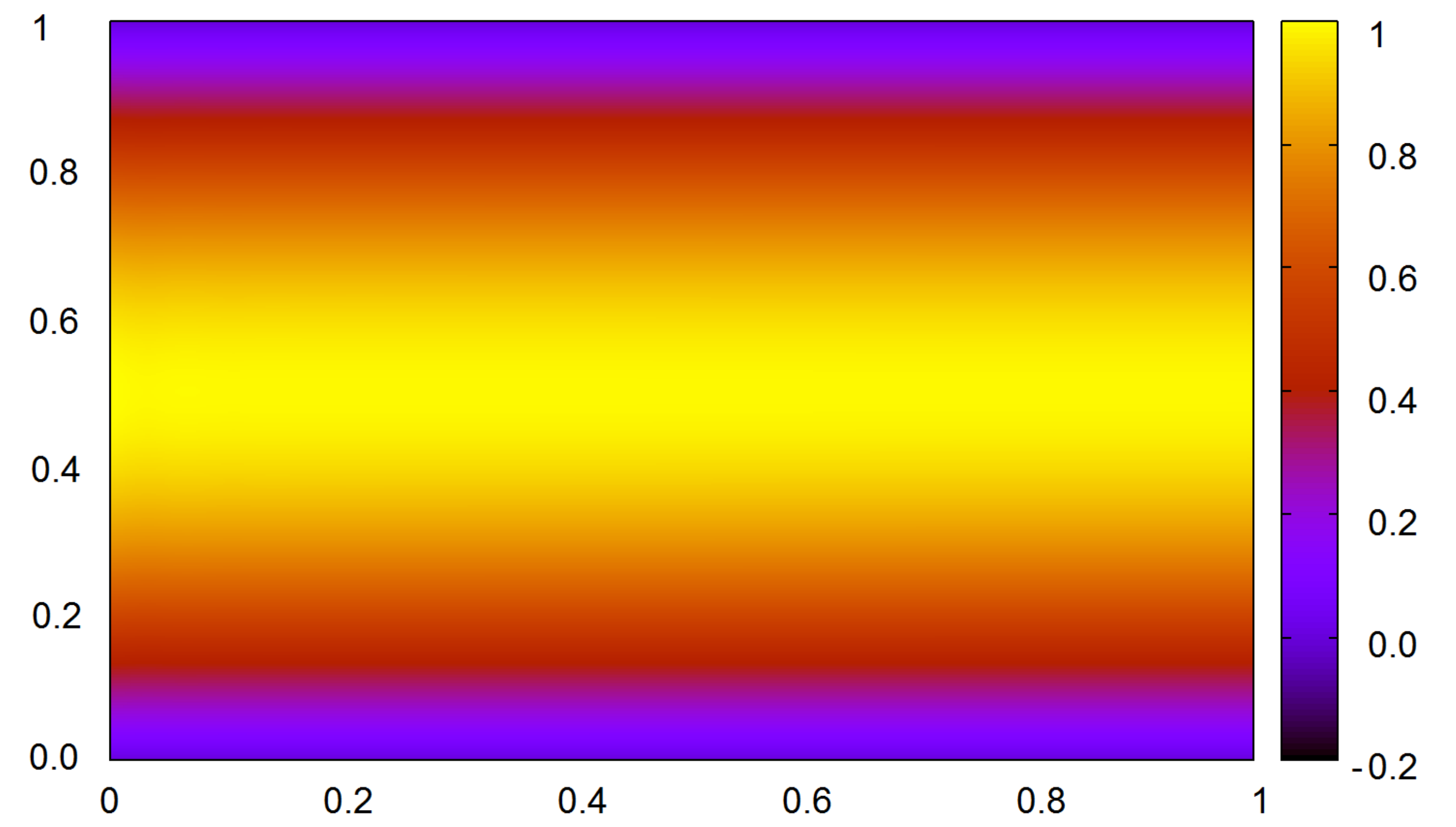}
\caption{\revcc{Solution to the Erikkson-Johnson problem on the mesh from iteration 20.}}
\label{fig:Crosses3}
\end{figure}

\revcc{We illustrate the solutions obtained from the iGRM and SUPG methods in Figures~\ref{fig:Crosses1}-\ref{fig:Crosses3}. We report the convergence in \revMP{$L^2$} and \revMP{$H^1$} norms in Table~\ref{tab:iGRMSUPG_E}.}

We now compare \revMP{iGRM vs} the DPG method~\cite{ Erikkson} using Lagrange (2,0) polynomials for solution and (3,0) polynomials for the residual. In our method, we use a (2,1) space for the solution with either a (3,1) or (3,2) space for the residual. Thus, our solution spaces have identical order and higher continuity, so they are smaller than those used in~\cite{ Erikkson}.

\revcc{We compare the iGRM computations with the DPG results from Figure 5.3 in~\cite{ Erikkson}, which shows the $L^2$ errors. Both iGRM and DPG meshes are adapted.  The DPG method~\cite{ Erikkson} for $Pe=10^6$ delivers \revMP{an} \revMP{$L^2$} norm error of order 0.02 for mesh size of the order 2000. The iGRM method for mesh dimension $20\times 10$ delivers \revMP{$L^2$} norm error 0.02 on mesh size 1312 with a residual space of (3,1).  Summing up, for $Pe=10^6$, iGRM delivers similar quality solutions on smaller grids. We achieve this by using smooth functions with higher continuity approximations. Additionally, the iGRM implementation is simpler since it does not require breaking the spaces. Nevertheless, the present solution strategy for iGRM is limited to tensor-product meshes. }

\revcc{Comparing the iGRM with SUPG, we conclude that iGRM converges to similar quality solutions than the SUPG method. The computational cost is higher for iGRM, but the iGRM method does not require the determination of problem-specific parameters.}  



\begin{table*}[htp]
\begin{center}
\begin{tabular}{|c|c|c|}
\hline 
\includegraphics[scale=0.5]{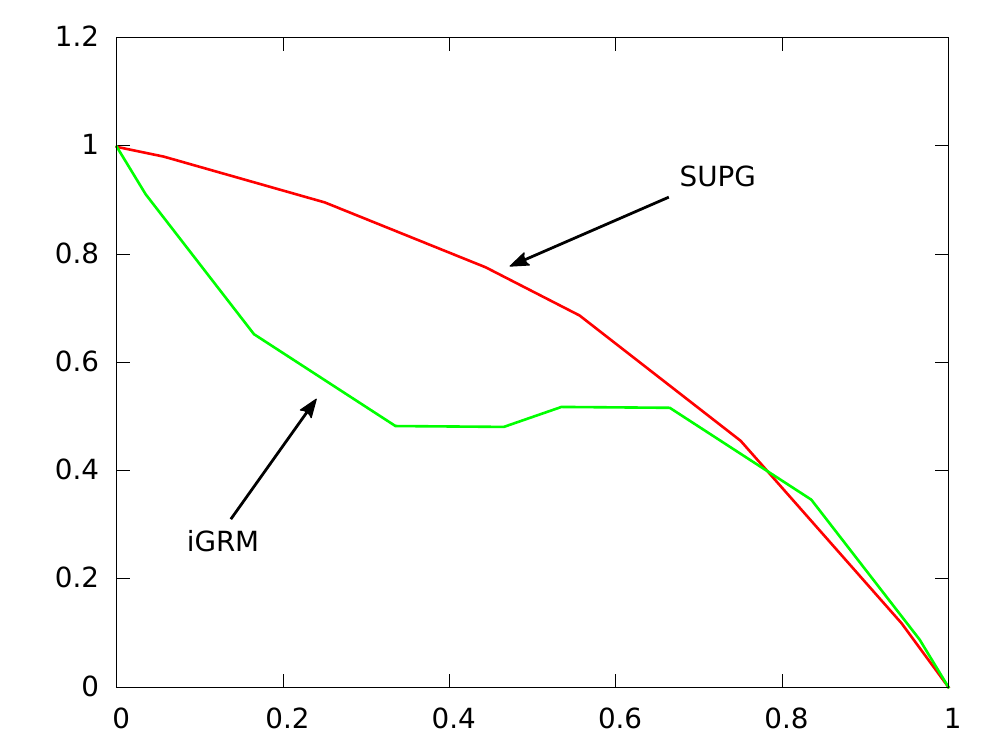} & \includegraphics[scale=0.5]{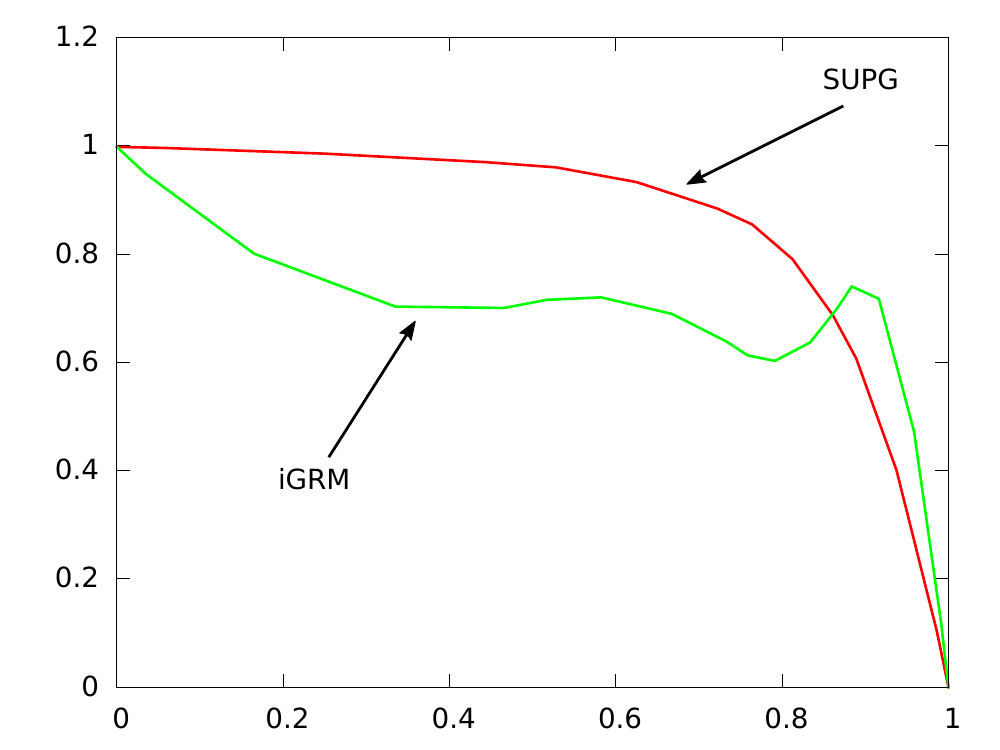} &
\includegraphics[scale=0.5]{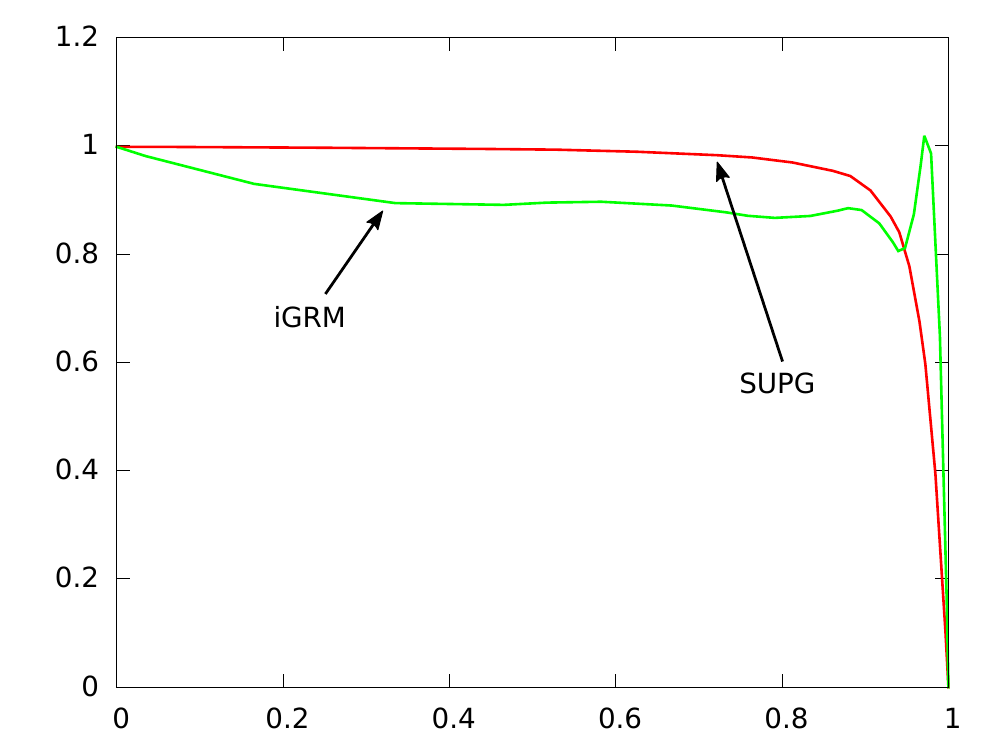} \\
iteration 1 [0  0.5  1]   &  iteration 3  [0  0.5  0.75  0.875  1]  & iteration 5 [0  0.5  0.75  0.875 \\ & &  0.9375  0.96875  1]  \\
\hline
\includegraphics[scale=0.5]{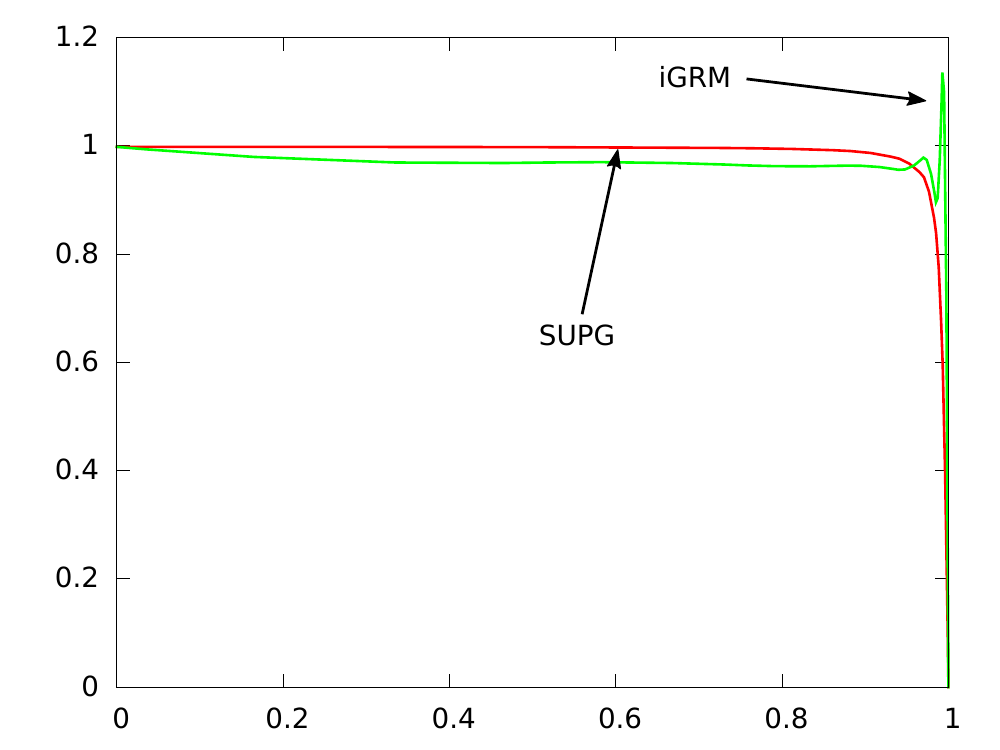} & \includegraphics[scale=0.5]{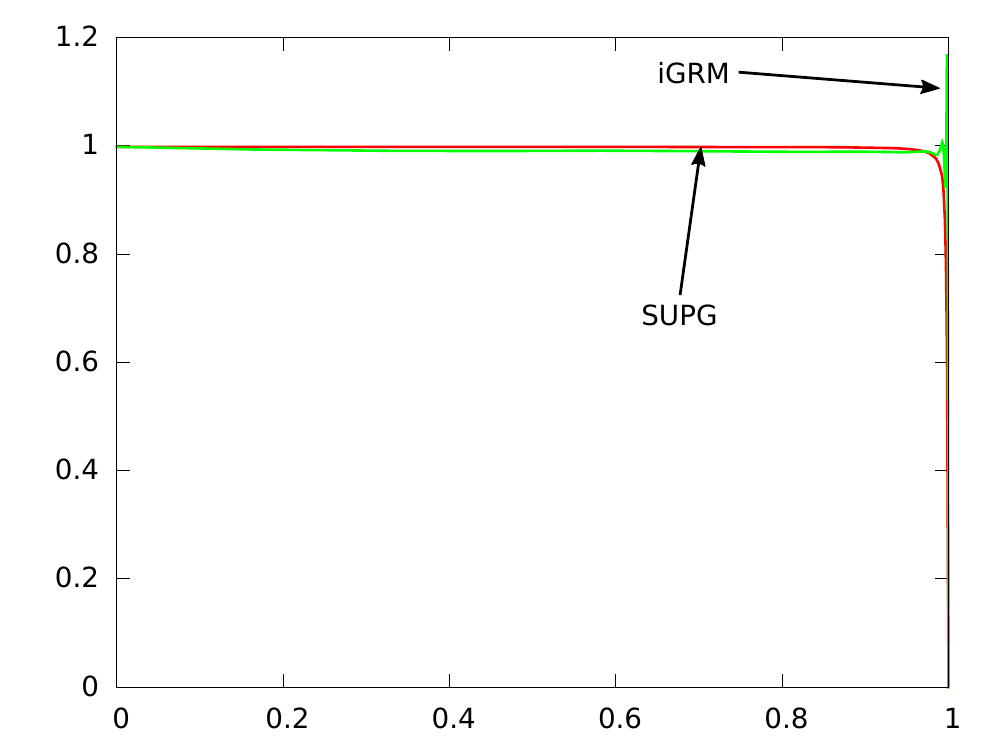} &
\includegraphics[scale=0.5]{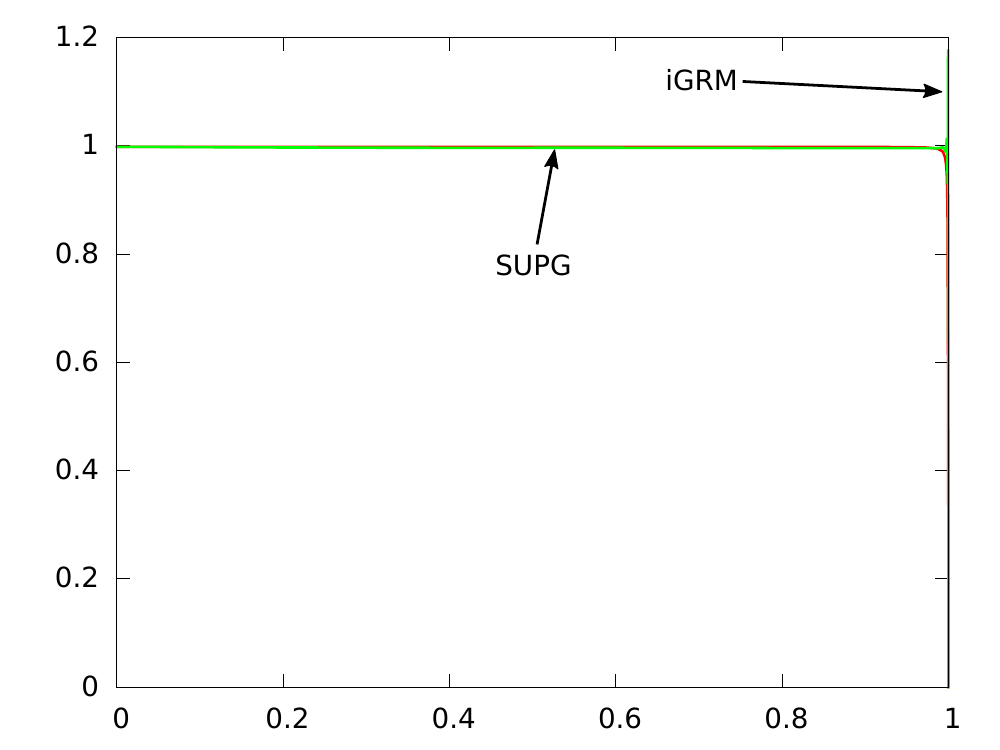} \\
iteration 7  [0  0.5  0.75  0.875  & iteration 9  [0  0.5  0.75  0.875  & iteration 11 [0  0.5  0.75  0.875  \\
 0.9375  0.96875  0.984375  1]  & 0.9375  0.96875  ...  & 0.9375  0.96875  ... \\
                                                & 0.998046875  1]                 &  0.998046875 0.9990234375  \\ 
                                                &                                           & 0.9995117188  1] \\
\hline
\includegraphics[scale=0.5]{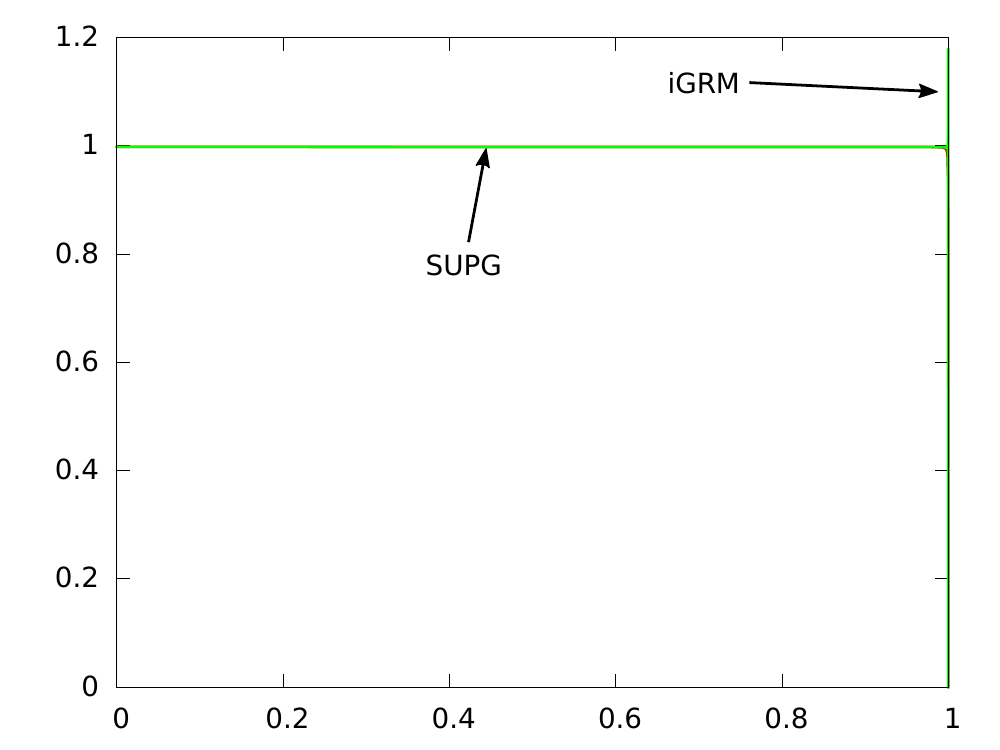} & \includegraphics[scale=0.5]{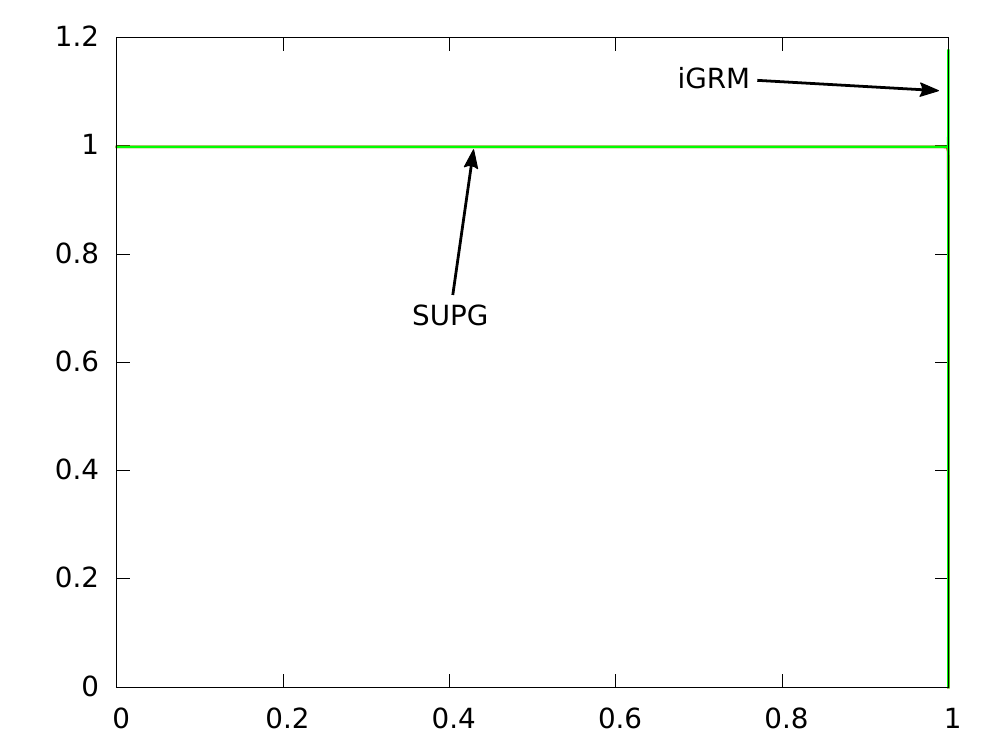} &
\includegraphics[scale=0.5]{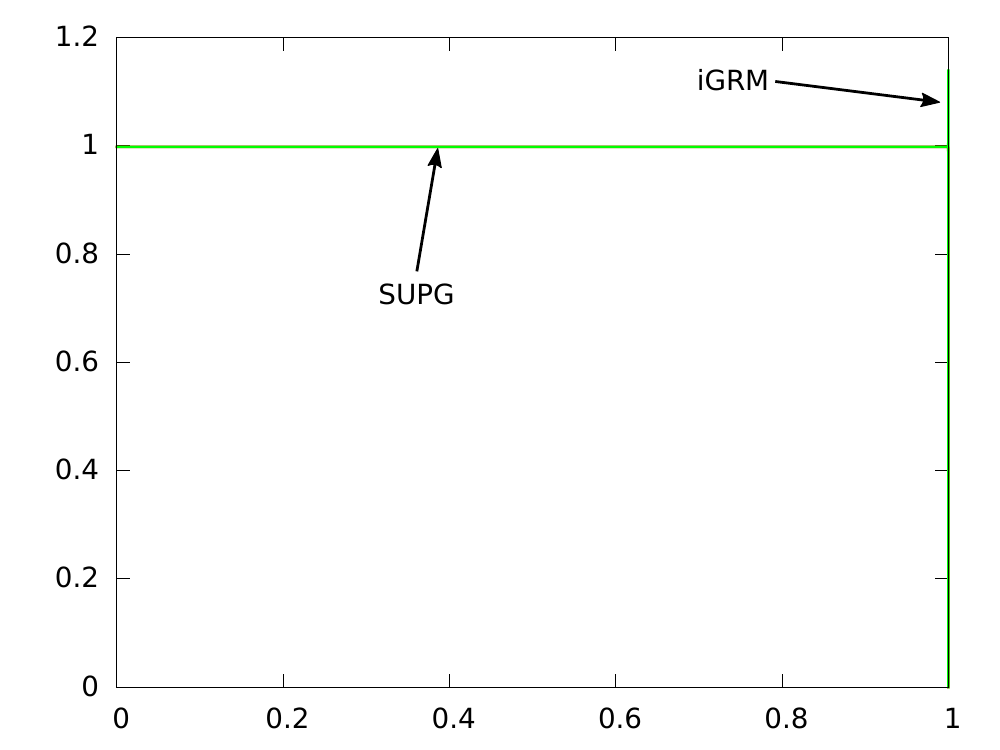} \\
iteration 13 [0  0.5  0.75  0.875   & iteration 15 [0  0.5  0.75  0.875  & iteration 17 [0  0.5  0.75  0.875  \\
0.9375  0.96875  ...  &  0.9375  0.96875  ... &  0.9375  0.96875  ... \\
  0.9998779297  1] &   0.9998779297  0.9999389648 &   0.9998779297  0.9999389648 \\
					  &     0.9999694824  1]  &   0.9999694824  0.9999847412 \\
					  &                                &    0.9999923706  1] \\
\hline
\end{tabular}
\end{center}
\caption{\revMP{Solutions to the Erikkson-Johnsson problem by using iGRM and SUPG methods for Pe=1,000,000 over a sequence of grids. We report the knot vector below each picture.}}
\label{fig:Crosses1}
\end{table*}

\begin{table*}[htp]
\begin{center}
\begin{tabular}{|c|c|c|}
\hline 
\includegraphics[scale=0.5]{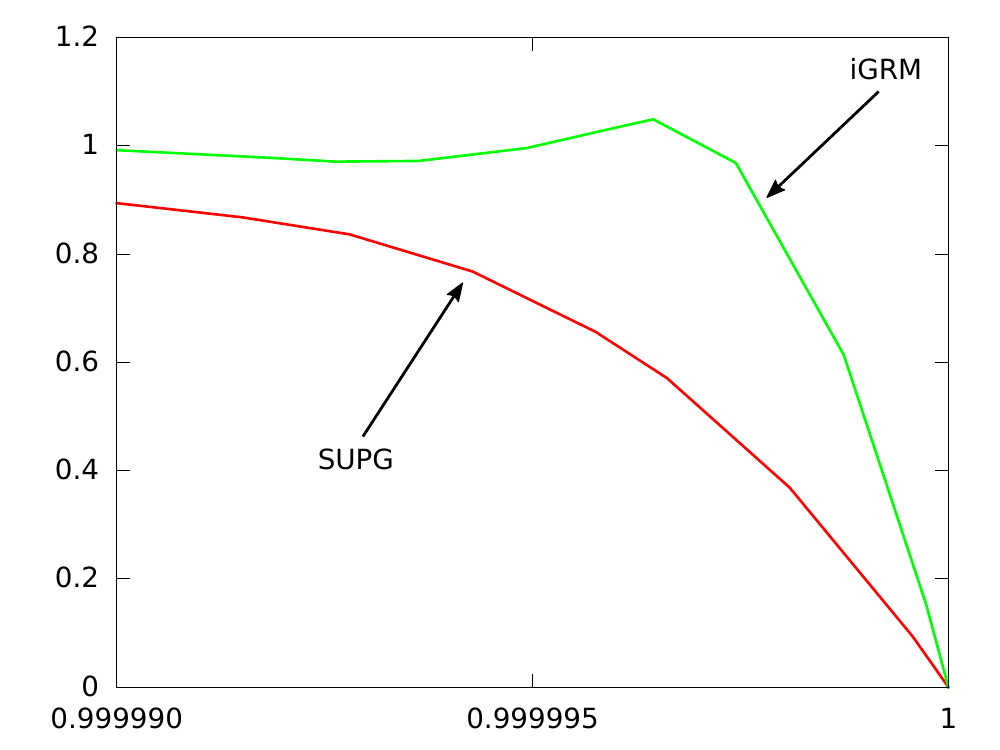} & \includegraphics[scale=0.5]{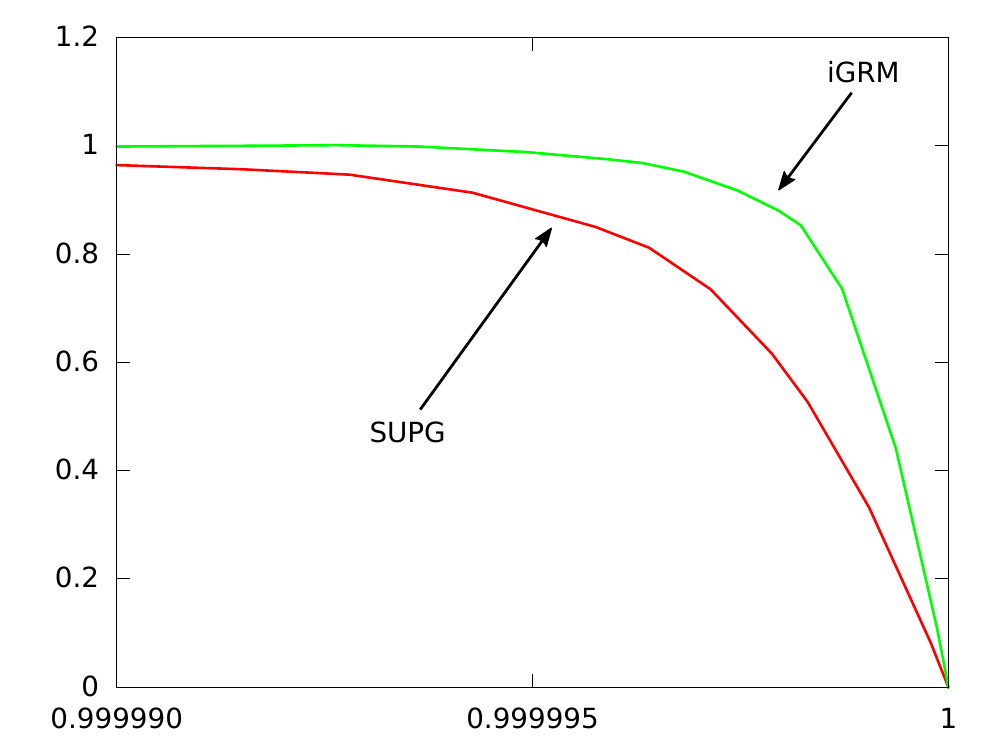} &
\includegraphics[scale=0.5]{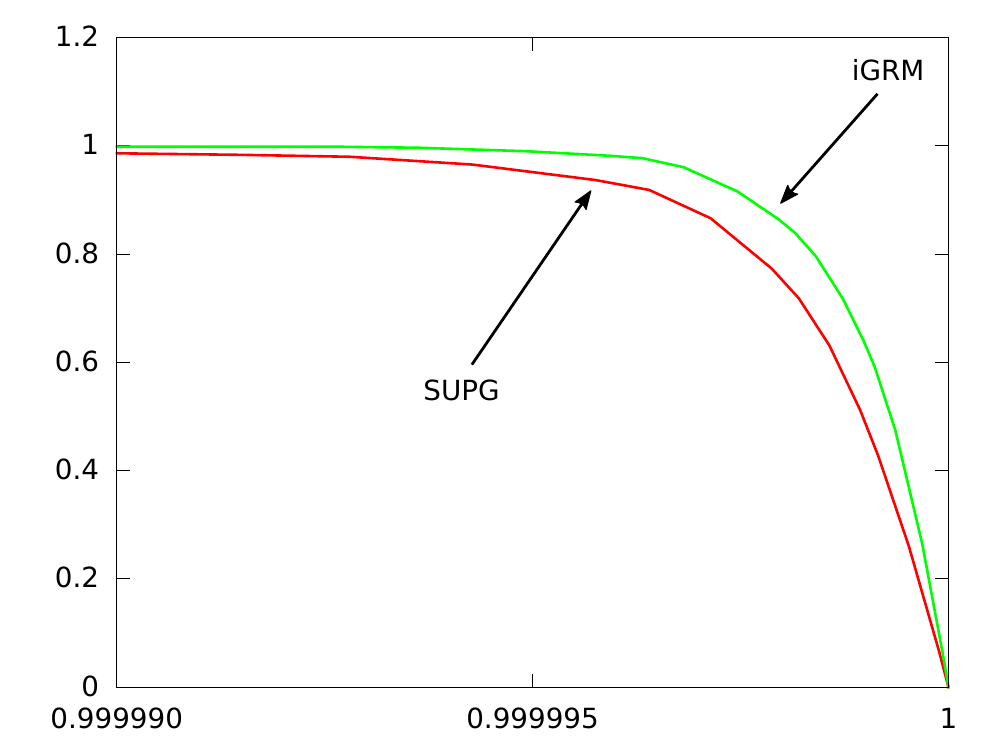} \\
iteration 18 [0  0.5  0.75 &  iteration 19 [0  0.5  0.75 & iteration 20 [0  0.5  0.75 \\
 0.875  0.9375  0.96875 ... &  0.875  0.9375  0.96875 ... &  0.875  0.9375  0.96875 ... \\
0.9999961853  1]  &   0.9999980927  1] & 0.9999980927 0.9999990463  1] \\
\hline
\includegraphics[scale=0.5]{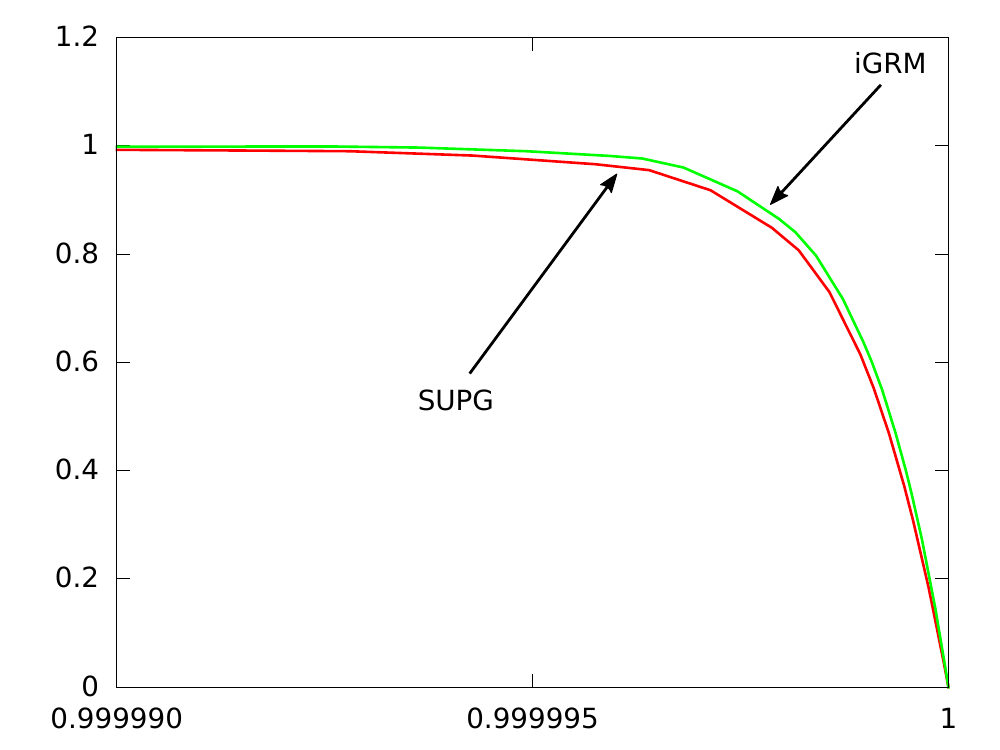} & \includegraphics[scale=0.5]{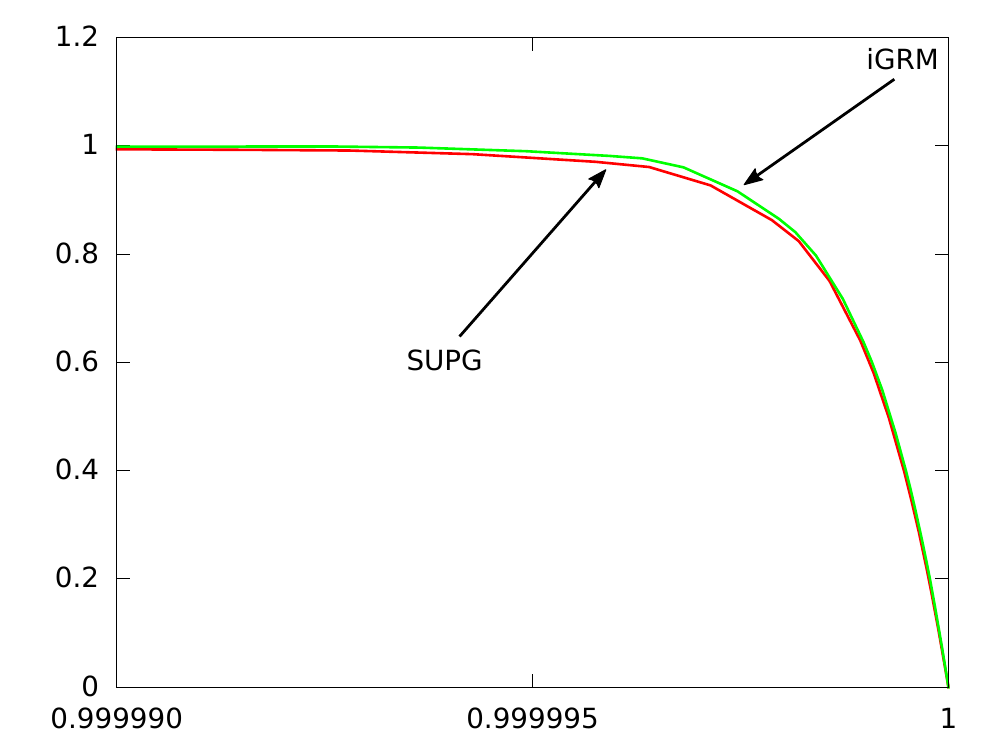} &
\includegraphics[scale=0.5]{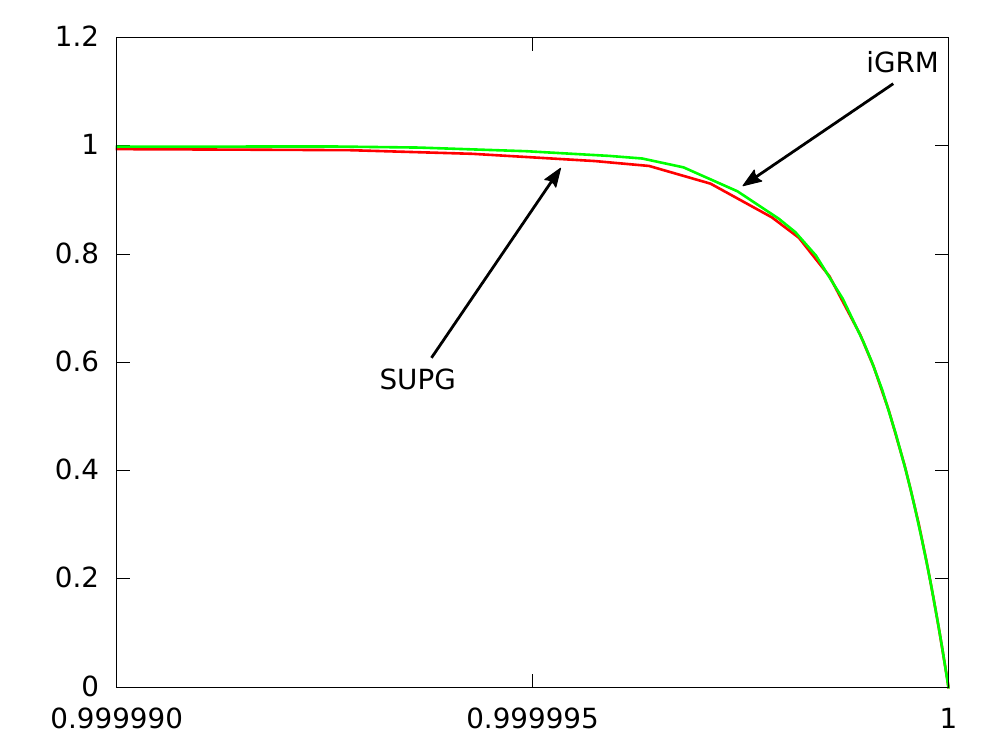} \\
iteration 21  [0  0.5  0.75  0.875  & iteration 22  [0  0.5  0.75  0.875  & iteration 23 [0  0.5  0.75  0.875  \\
0.9375 0.96875    ... & 0.9375 0.96875    ... &  0.9375 0.96875    ... \\
0.9999990463  0.9999995232  1]  & 0.9999990463  0.9999995232  & 0.9999990463  0.9999995232 \\
							& 0.9999997616  1]   & 0.9999997616  0.9999998808  1] \\
\hline
\includegraphics[scale=0.5]{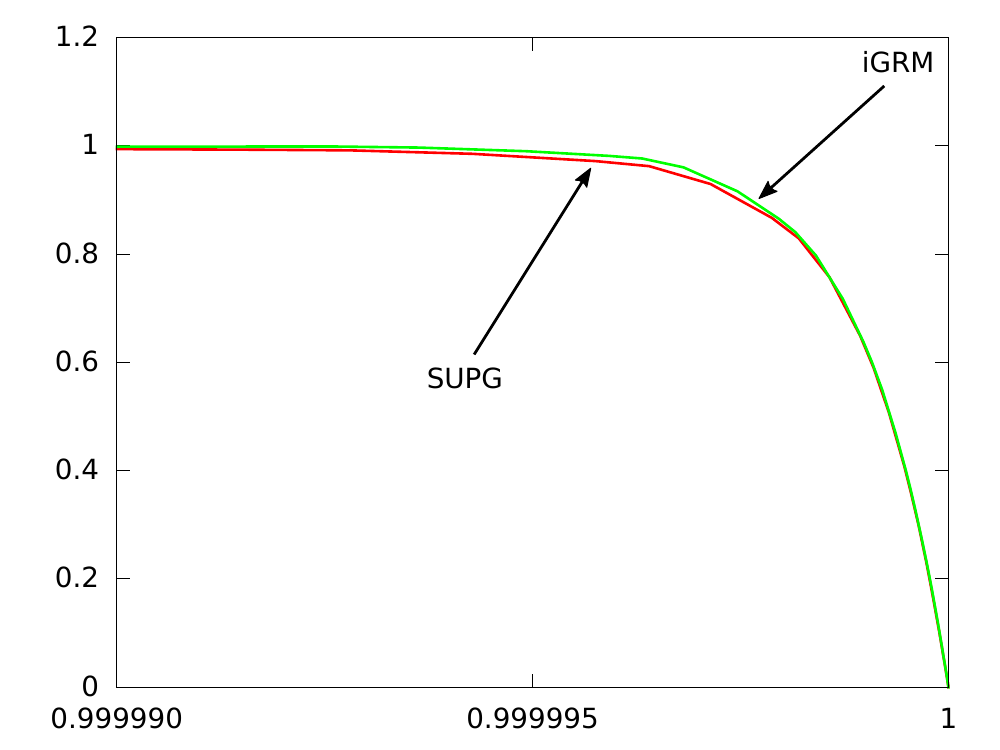} & \includegraphics[scale=0.5]{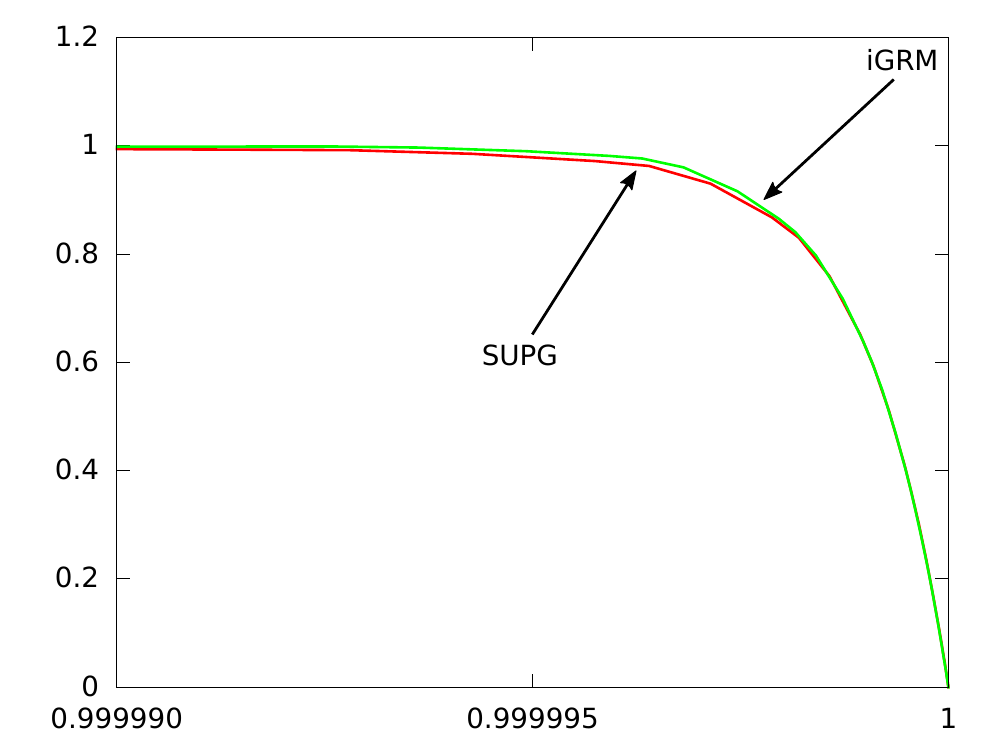} &
\includegraphics[scale=0.5]{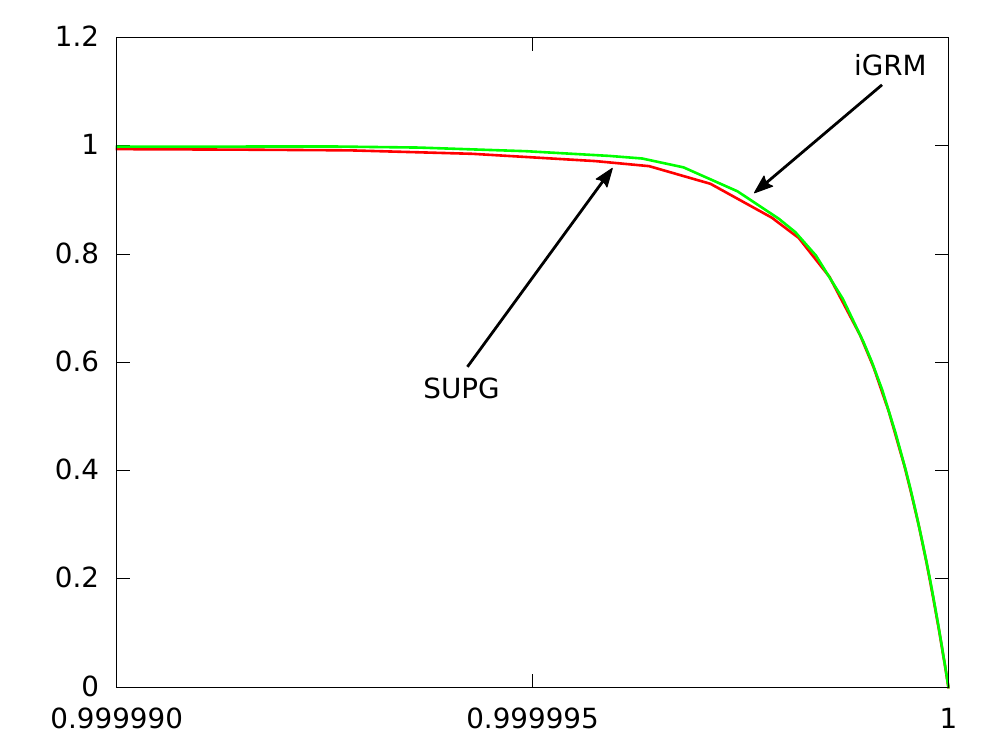} \\
iteration 21  [0  0.5  0.75  0.875  & iteration 22  [0  0.5  0.75  0.875  & iteration 23 [0  0.5  0.75  0.875  \\
0.9375 0.96875    ... & 0.9375 0.96875    ... &  0.9375 0.96875    ... \\
0.9999999404  1]                    & 0.9999999404  0.9999999702  1] & 0.9999999404  0.9999999702 \\								 &                                                            &  0.9999999851  1] \\
\hline
\end{tabular}
\end{center}
\caption{\revMP{Zoom to the right corner of the solutions to the Erikkson-Johnsson problem by using iGRM and SUPG methods for Pe=1,000,000 over a sequence of adapted grids. We report the knot vector below each picture.}}
\label{fig:Crosses2}
\end{table*}

\begin{table*}[htp]
\begin{center}
\begin{tabular}{cccc|ccc}
& \multicolumn{3}{c}{ SUPG (2,1) } & \multicolumn{3}{c}{ iGRM (2,1) (3,1) } \\
iteration & \#NRDOF & L2 & H1  & \#NRDOF & L2 & H1  \\
\hline
1 & 24 & 44.33 & 56.37 & 84 & 53.59 & 68.00  \\
2 & 30 & 33.79 & 56.29 & 110 & 44.00 & 66.03  \\
3 & 36 & 24.74 & 66.30 & 136 & 31.86 & 90.46  \\
4 & 42 & 17.75 & 87.73 & 162 & 20.31 & 143.38  \\
5 & 48 & 12.63 & 121.57 & 188 & 12.42 & 227.77 \\
6 & 54 & 8.96 & 170.93 & 214 & 7.36 & 345.20 \\
7 & 60 & 6.34 & 241.31 & 240 & 4.49 & 506.91 \\
8 & 66 & 4.49 & 341.08 & 266 & 2.85 & 731.00 \\
9 & 72 & 3.18 & 482.28 & 292 & 1.89 & 1044.11 \\
10 & 78 & 2.26 & 681.99 & 318 & 1.31 & 1484.00 \\
11 & 84 & 1.61 & 964.41 & 344 & 0.94 & 2103.96 \\
12 & 90 & 0.70 & 1363.73 & 370 & 0.70 & 2979.06 \\
13 & 96 & 0.84 & 1928.06 & 396 & 0.55 & 4026.87 \\
14 & 102 & 0.63 & 1623.21 & 422 & 0.45 & 353.16 \\
15 & 108 & 0.49 & 104.06 & 448 & 0.39 & 76.60 \\
16 & 114 & 0.39 & 77.93 & 474 & 0.35 & 69.25 \\
17 & 120 & 0.34 & 78.38 & 500 & 0.34 & 58.28 \\
18 & 126 & 0.31 & 68.11 & 526 & 0.34 & 34.03 \\
19 & 132 & 0.30 & 48.56 & 552 & 0.34 & 12.57 \\
20 & 138 & 0.30 & 27.56 & 578 & 0.34 & 3.97 \\
21 & 144 & 0.30 & 9.59 & 604 & 0.34 & 2.45 \\
22 & 150 & 0.30 & 4.40 & 630 & 0.34 & 2.35 \\
23 & 156 & 0.30 & 2.69 & 656 & 0.34 & 2.35 \\
24 & 162 & 0.30 & 2.85 & 682 & 0.34 & 2.35 \\
25 & 168 & 0.30 & 2.85 & 708 & 0.34 & 2.35 \\
\hline
\end{tabular}
\end{center}
\caption{\revcc{Convergence of SUPG and iGRM on grids refined from \revMP{a} $2\times 2$ uniform grid, for the Erikkson-Johnson problem with $Pe=10^6$, for (2,1) B-splines. Here $\#NDOF$ is the total number of degrees of freedom, and \revMP{$L^2$} and \revMP{$H^1$} are the exact norm errors. We use $\eta=0.0001$.}}
\label{tab:iGRMSUPG_E}
\end{table*}


\revcc{We also investigate the convergence and the execution times of the iterative solver. We show in Table \ref{tab:CG} the dependence of the number of iterations of the inner and outer loops on the parameter $\eta$ from the definition of the inner product. We use the Erikkson-Johnson problem as a reference. We \revMP{implement} the iGRM code in the IGA-ADS software \cite{ED4}, which automatically uses all cores of the machine. We run the experiments on Intel i7 processor with 2.7GHz with 8 cores and 16 GB of RAM.}

\begin{table*}[htp]
\begin{center}
\begin{tabular}{ccccccc}
iteration & $\eta$ & \# outer & \# inner &  L2 & H1 & time [ms] \\
\hline
1 & 0.1 & $>100$ & $>100$ & 12.76 & 225.88 & 37 \\
1 & 0.01 & 10 & 32 & 24.95 & 166.275 & 12 \\
1 & 0.001 & 4 & 18 & 49.86 & 73.11 & 6 \\
1 & 0.0001 & 3 & 16 & 53.59 & 68.00 & 5 \\
1 & 0.00001 & 2 & 12 & 58.90 & 66.36 & 4 \\
1 & 0.000001 & 2 & 11 & 59.02 & 66.45 & 4 \\
1 & 0.0000001 & 2 & 11 & 59.02 & 66.46 & 4 \\
1 & 0.00000001 & 2 & 11 & 59.02 & 66.46 & 4 \\
\hline
20 & 0.1 & $>100$ & $>100$ & 0.34 & 42.96 & 132 \\
20 & 0.01 & $>100$ & $>100$ & 0.34 & 8.24 &103 \\
20 & 0.001 & 5 & 69 & 0.34 & 4.07 & 30 \\
20 & 0.0001 & 2 & 51 & 0.34 & 3.97 & 21\\
20 & 0.00001 & 2 & 63 & 0.34 & 3.97 & 26 \\
20 & 0.000001 & 2 & 65 & 0.34 & 3.97 & 23 \\
20 & 0.0000001 & 2 & 81 & 0.59 & 4.05 & 31 \\
20 & 0.00000001 & 2 & 85 & 4.73 & 8.91 & 34 \\
\hline
\end{tabular}
\end{center}
\caption{Convergence of the inner (CG) and outer (corrections) loops of the iterative solver for different $\eta$ parameter from the weighted norm. We use the first and the last (number 20) adaptive grids for the Erikkson-Johnson problem with $Pe=10^6$, for trial (2,1) and test (3,1) B-splines. Here \revMP{$L^2$} and \revMP{$H^1$} are the exact norm errors.}
\label{tab:CG}
\end{table*}


\subsection{Vortical wind problem}

We now solve a vortical flow with the advection-dominated diffusion equations 
over the rectangular domain $\Omega=(0,1)\times(-1,1)$, with zero right-hand side $f=0$, 
the advection vector $\beta(x,y)=\left(\beta_x(x,y),\beta_y(x,y)\right)=\left(-y,x\right)$. This velocity field is modeling a rotating flow.
We introduce 
$\Gamma_{1}=\{(x,y):x=0, 0.5\leq y \leq 1.0\}$, 
$\Gamma_{2}=\{(x,y):x=0, 0.0\leq y \leq 0.5\}$\revMP{.}
The Dirichlet boundary conditions are
\begin{eqnarray}
g = \frac{1}{2}\left(\tanh \left(\left(|y|-0.35\right)\frac{b}{\epsilon}\right)+1\right), \textrm{ for }x\in\Gamma_2 \nonumber \\ 
g = \frac{1}{2}\left(0.65-\tanh \left(\left(|y|\right)\frac{b}{\epsilon}\right)+1\right), \textrm{ for }x\in\Gamma_1 \\ 
g=0, \textrm{ for }x\in\Gamma \setminus \Gamma_1\cup\Gamma_2 
\label{eq:Problem4_bc}
\end{eqnarray}

\begin{figure}
\centering
\includegraphics[scale=0.45]{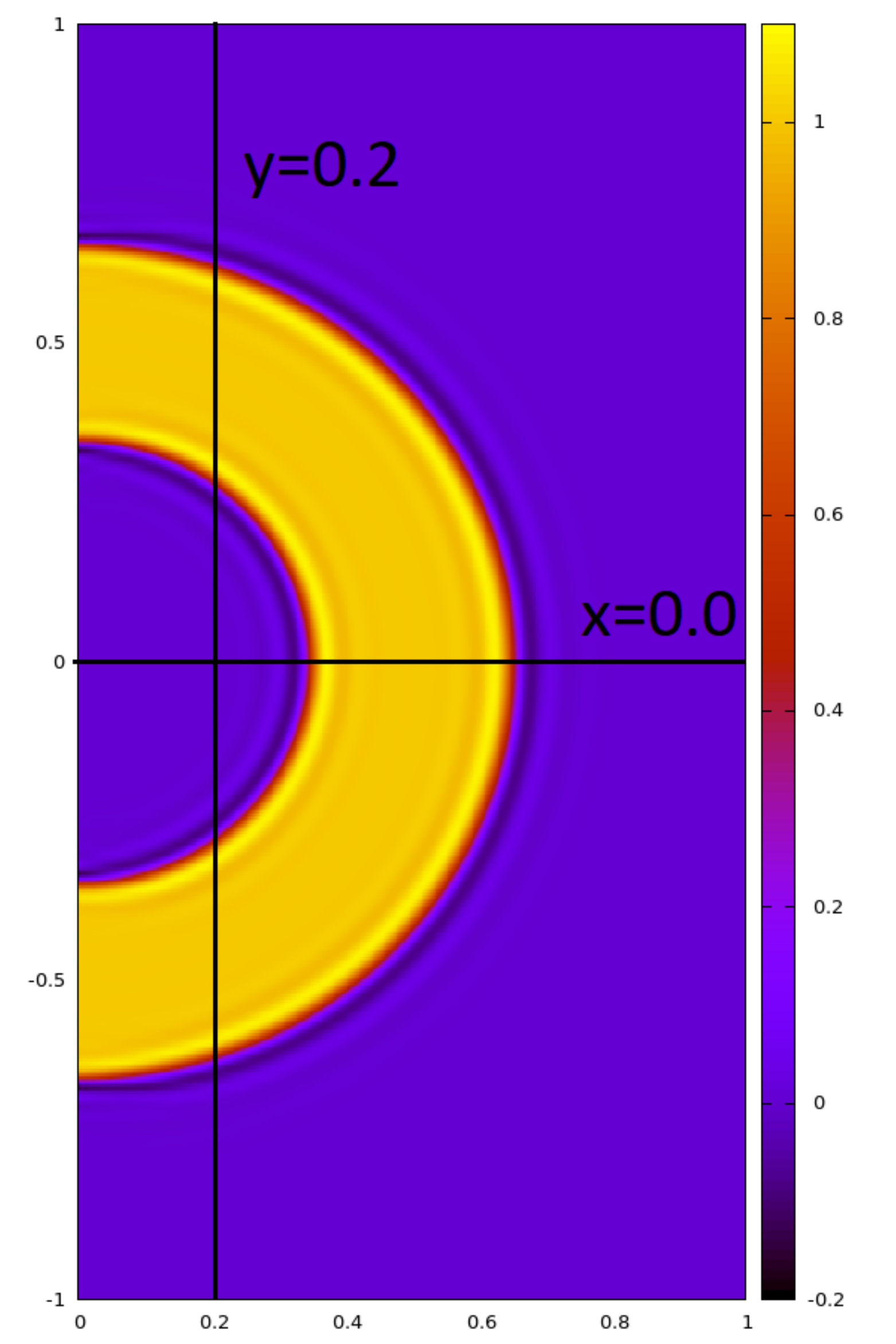}
\caption{Solution to the circular wind problem on a mesh of $128\times
  128$ elements with trial (2,1), test (2,0), for P\'eclet number $Pe=10^6$, wind force $b=1$. We also include two planes (cross-sections) whose solutions are presented in Figures~\ref{fig:Problem4b}-\ref{fig:Problem4c}.}
\label{fig:Problem4a}
\end{figure}

\begin{figure}
\centering
\includegraphics[scale=0.5]{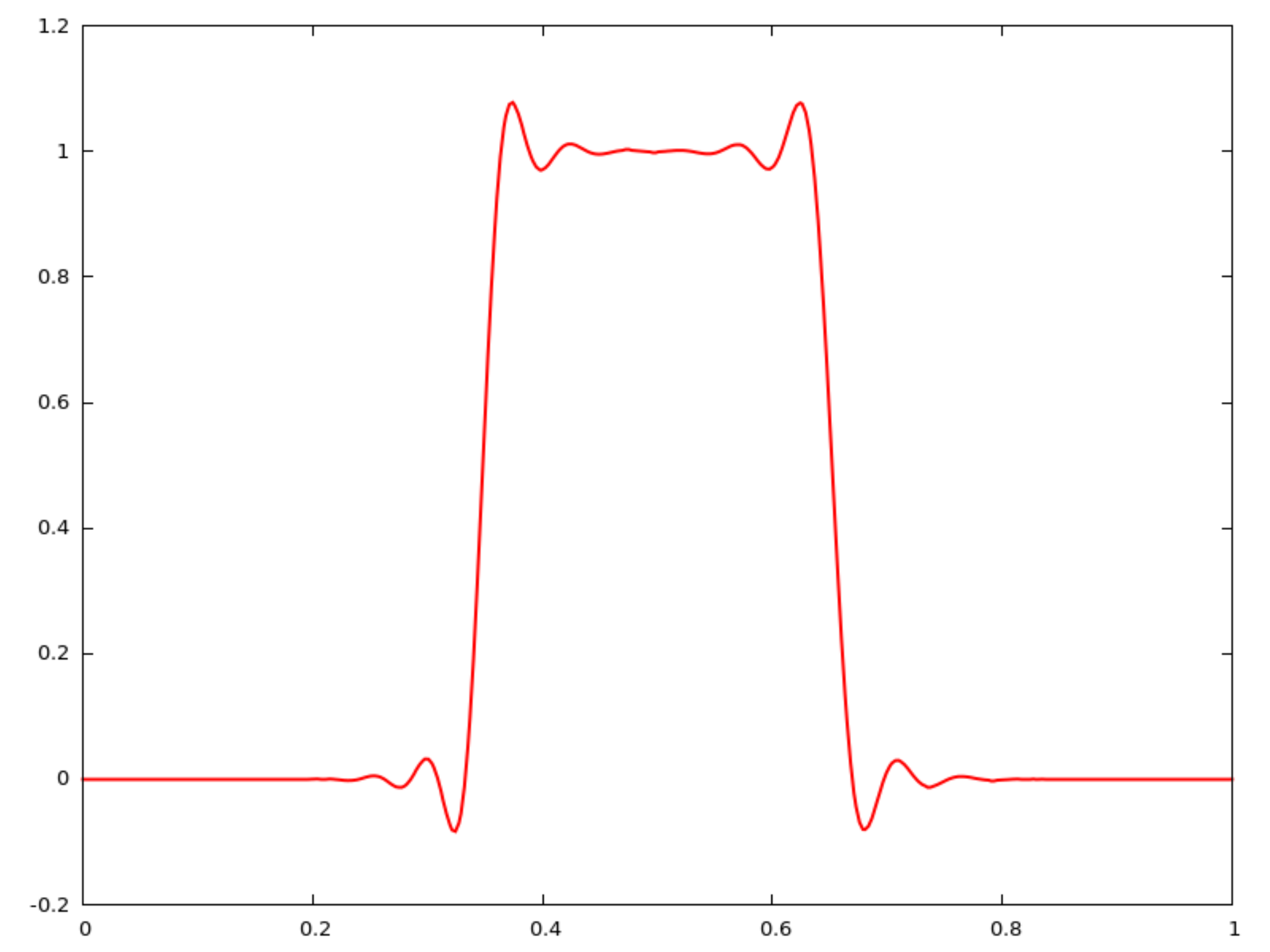}
\caption{Horizontal cross-section at $x=0$ through the solution to the circular wind problem on a mesh of $128\times 128$ elements with trial (2,1), test (2,0), for P\'eclet number $Pe=10^6$, wind force $b=1$.}
\label{fig:Problem4b}
\end{figure}

\begin{figure}
\centering
\includegraphics[scale=0.5]{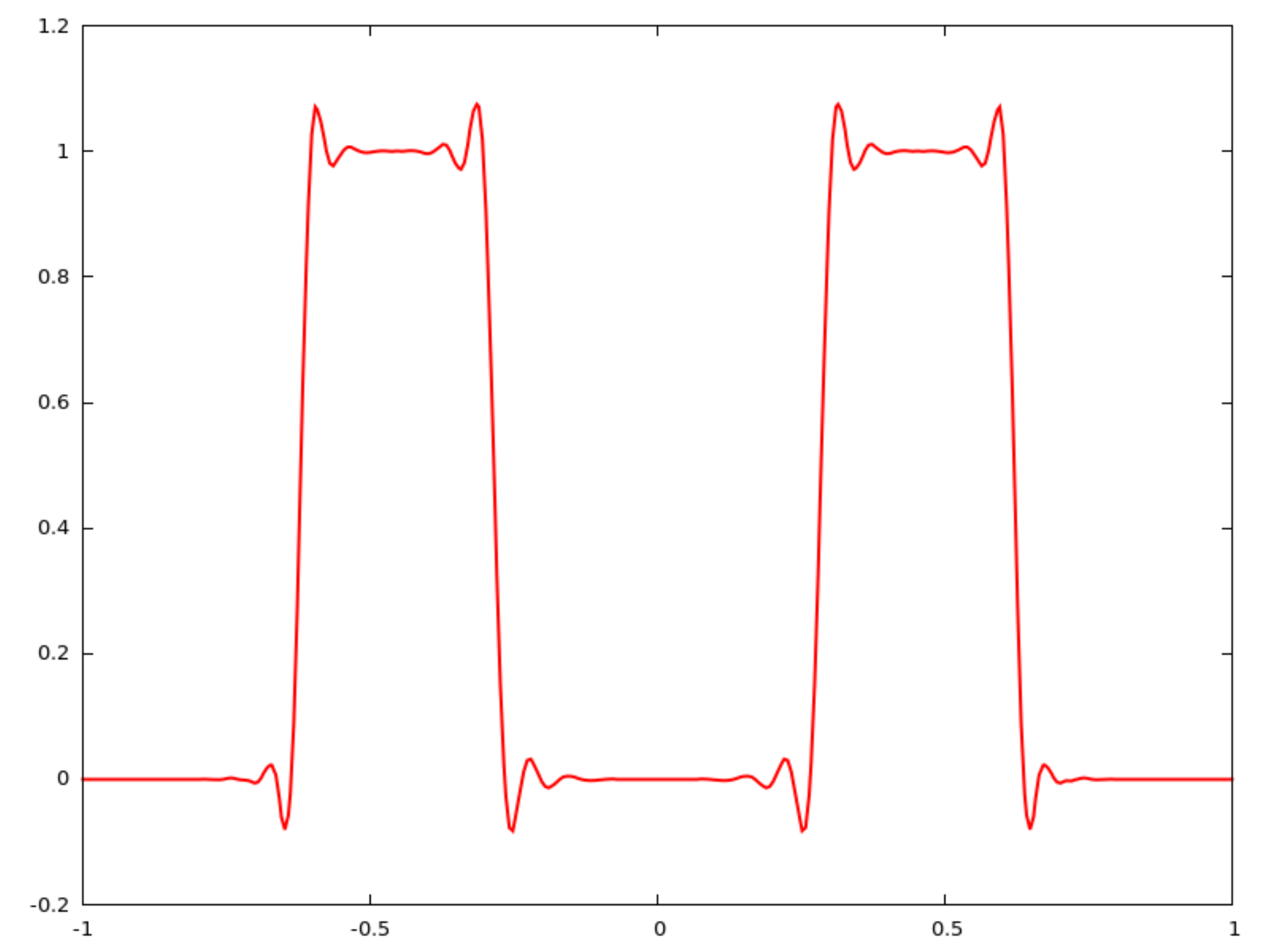}
\caption{Vertical cross-section at $y=0.2$ through the solution to the circular wind problem on a mesh of $128\times 128$ elements with trial (2,1), test (2,0), for P\'eclet number $Pe=10^6$, wind force $b=1$.}
\label{fig:Problem4c}
\end{figure}

We use the iGRM setup (\ref{eq:resmin}) with the preconditioned CG solver described in Section 4. We use $128\times 128$ mesh with solution (2,1) residual (2,0). We use a P\'eclet number $Pe=10^6$. The numerical results are summarized in Figures~\ref{fig:Problem4a}-\ref{fig:Problem4c}.

There are some overshoots and undershoots in the cross-sections of the mesh. They can be removed by performing mesh refinements. They can be also removed by using penalty methods~\cite{ yzbeta1,yzbeta2}.


\begin{remark}
The number of iterations of the conjugate gradient solver grows with the problem size, as presented in Table \ref{tab:CG}. We will explore possible solutions to this issue in future work. 
\revMP{The computational cost of the solver grows like ${\cal O}(Nk)$ where $N$ is the number of degrees of freedom, and $k$ is the number of iterations.}
\end{remark}

\section{Conclusions}
\label{sec:conc}

\revcc{We introduce a stabilized isogeometric method that uses residual minimization in a dual norm.  We design the method to achieve good solution properties. The method exploits the Kronecker product structure of the computational problem to deliver fast solutions. The solution space in our scheme uses maximum continuity B-splines. To accelerate the solution of the algebraic scheme, we introduce a fast solver for the Gramm matrix, but without a Schur complement preconditioner.  We call our method isogeometric residual minimization (iGRM) with direction splitting preconditioner.  We verify the accuracy of the solution on four stationary problems.  In this method, the diffusion and advection coefficient functions can be arbitrary.  Our future work will extend this method to other problems, such as the Stokes~\cite{ evans2013isogeometric} and Maxwell problems~\cite{ maxwell1, maxwell2,maxwell3}, the development of the method for time-dependent problems~\cite{ timedep}, as well as the development of the parallel software dedicated to the simulations of different non-stationary problems with the iGRM method.  We will seek to find a proper preconditioner for the CG problem as well as an alternative method to speed up the solution of the saddle-point problem that the residual minimization delivers, e.g., based on~\cite{ Method1, Method2}.}

\subsection*{Acknowledgments}

The National Science Centre, Poland,  partially supported this work, grant no. 2017/26/M/ ST1/ 00281.  This publication was also made possible in part by the CSIRO Professorial Chair in Computational Geoscience at Curtin University and the Deep Earth Imaging Enterprise Future Science Platforms of the Commonwealth Scientific Industrial Research Organisation, CSIRO, of Australia. Additional support was provided by the European Union's Horizon 2020 Research and Innovation Program of the Marie Skłodowska-Curie grant agreement No. 777778. At Curtin University, support was provided by The Institute for Geoscience Research (TIGeR) and by the Curtin Institute for Computation.


\begin{thebibliography}{99}

\bibitem{Samarskii}
A. A. Samarskij, E. S. Nikolaev, Numerical Methods for Grid Equations: Volume II Iterative Methods, Birkhäuser, Basel, Boston, Berlin (2012)
\bibitem{Quarteroni}
A. Quarteroni, R. Sacco, F. Saleri, Numerical Mathematics (Texts in Applied Mathematics) 2nd Edition, Springer, Berling, Heidelberg, New, York (2006)
\bibitem{ADS1} 
D.W. Peaceman, H.H. Rachford Jr., The numerical solution of parabolic and elliptic differential equations, 
Journal of Society of Industrial and Applied Mathematics 3 (1955) 28–41.
\bibitem{ADS2}
J. Douglas, H. Rachford, On the numerical solution of heat conduction problems in two and three space variables, 
Transactions of American Mathematical Society 82 (1956) 421–439.
\bibitem{ADS3}
E.L. Wachspress, G. Habetler, An alternating-direction-implicit iteration technique, 
Journal of Society of Industrial and Applied Mathematics 8 (1960) 403–423.
\bibitem{ADS4}
G. Birkhoff, R.S. Varga, D. Young, Alternating direction implicit methods, Advanced Computing 3 (1962) 189–273.
\bibitem{Minev1}
J. L. Guermond, P. Minev, {\em A new class of fractional step techniques for the incompressible
Navier-Stokes equations using direction splitting,} Comptes Rendus Mathematique 348(9-10) (2010) 581–585.
\bibitem{Minev2}
J. L. Guermond, P. Minev, J. Shen, {\em An overview of projection methods for incompressible
flows,} Computer Methods in Applied Mechanics and Engineering, 195 (2006) 6011–6054.
\bibitem{IGA} 
J. A. Cottrell, T. J. R. Hughes, Y. Bazilevs, {\em Isogeometric Analysis: Toward Unification of CAD and FEA}
John Wiley and Sons, (2009)
\bibitem{NURBS}
L. Piegl, and W. Tiller, \emph{The NURBS Book (Second Edition)},
Springer-Verlag New York, Inc., (1997).
\bibitem{CG1}
L. Gao, V.M. Calo, Fast Isogeometric Solvers for Explicit Dynamics, Computer Methods in Applied Mechanics and Engineering, 274 (1) (2014) 19-41.
\bibitem{CG2}
L. Gao, V.M. Calo, Preconditioners based on the alternating-direction-implicit algorithm for the 2D steady-state diffusion equation with orthotropic heterogeneous coefficients, 273 (1) (2015) 274-295.
\bibitem{Gao2014}
Longfei Gao, Kronecker Products on Preconditioning, PhD. Thesis, King Abdullah University of Science and Technology (2013).
\bibitem{ED1}
M. \L{}o\'{s}, M. Wo\'{z}niak, M. Paszy\'{n}ski, L. Dalcin, V.M. Calo, Dynamics with Matrices Possessing Kronecker Product Structure,  Procedia Computer Science 51 (2015) 286-295.
\bibitem{ED2}
M. Wo\'{z}niak, M. \L{}o\'{s}, M. Paszy\'{n}ski, L. Dalcin, V. Calo, Parallel fast isogeometric solvers for explicit dynamics, Computing and Informatics, 36(2) (2017) 423-448.
\bibitem{ED3}
M. \L{}o\'{s}, M. Paszy\'{n}ski, A. K\l{}usek, W. Dzwinel, Application of fast isogeometric L2 projection solver for tumor growth simulations, Computer Methods in Applied Mechanics and Engineering, 316 (2017) 1257-1269.
\bibitem{ED4}
M. \L{}o\'{s}, M. Wo\'{z}niak, M. Paszy\'{n}ski, A. Lenharth, K. Pingali,  IGA-ADS : Isogeometric Analysis FEM using ADS solver, Computer \& Physics Communications, 217 (2017) 99-116.
\bibitem{ED5}
G. Gurgul, M. Wo\'{z}niak, M. \L{}o\'{s}, D. Szeliga, M. Paszy\'{n}ski, Open source JAVA implementation of the parallel multi-thread  direction isogeometric L2 projections solver for material science simulations, Computer Methods in Material Science, 17 (2017) 1-11.
\bibitem{LSFEM}
P. Bochev, M. Gunzburger, {\em Least-Squares Finite Element Method}, Springer Applied Mathematical Sciences 166 (2009)
\bibitem{t8}
L. Demkowicz, J. Gopalakrishnan, Recent Developments in Discontinuous Galerkin Finite Element Methods for Partial Differential Equations (eds. X. Feng, O. Karakashian, Y. Xing). In: vol. 157. IMA Volumes in Mathematics and its Applications, (2014). An Overview of the DPG Method, 149–180
\bibitem{VSM}
A. Cohen, W. Dahmen, G. Welper,  {\em Adaptivity and Variational Stabilization for Convection-Diffusion Equations} Mathematical Modelling and Numerical Analysis 46(5) (2012) 1247-1273.
\bibitem{romkes}
V. M. Calo, A. Romkes, E. Valseth, Automatic Variationally Stable Analysis for FE Computations: An Introduction, \emph{https://arxiv.org/abs/1808.01888}
\bibitem{Kopteva}
N. Kopteva, E. O'Riordan, Shishkin meshes in the numerical solution of singularly perturbed differential equations, International Journal of Numerical Analysis and Modeling, 7(3) (2010) 393-415.
\bibitem{Erikkson}
J. Chan, J. A.Evans, A minimal-residual finite element method for the convection-diffusion equations, ICES-REPORT 13-12 (2013)
\bibitem{SUPG}
V. M. Calo, Residual-based multiscale turbulence modeling: Finite volume simulations of bypass transition, Stanford University, Ph.D. Thesis (2005)
\bibitem{SUPG2}
T.J.R. Hughes, L.P. Franca, M. Mallet, A new finite element formulation for fluid
dynamics: VI. Convergence analysis of the generalized SUPG formulation for linear time–
dependent multidimensional advective–diffusive systems, Computer Methods in Applied
Mechanics and Engineering, 6 (1987) 97–112.
\bibitem{evans2013isogeometric}
J. Evans, T. J. R. Hughes, Isogeometric divergence-conforming B-splines for the Darcy--Stokes--Brinkman equations, Mathematical Models and Methods in Applied Sciences, 23(4) (2013) 671-741.
\bibitem{maxwell1}
M. Hochbruck, T. Jahnke, R. Schnaubelt, Convergence of an ADI splitting for Maxwell's equations, Numerishe Mathematik, 129 (2015) 535-561.
\bibitem{maxwell2}
G. Liping, Stability and Super Convergence Analysis of ADI-FDTD for the 2D Maxwell Equations in a Lossy Medium, Acta Mathematica Scientia, 32(6) (2012) 2341-2368.

\bibitem{maxwell3}
M. Paszy\'{n}ski, L. Demkowicz, D. Pardo, Verification of goal-oriented $hp$-adaptivity, Computers \& Mathematics with Applications, 50(8-9) (2005) 1395-1404.

\bibitem{nocedal2006conjugate}
J.~Nocedal, S.~J. Wright, Conjugate gradient methods, Numerical optimization
  (2006) 101-134.
  
\bibitem[{Horn et~al.(1990)Horn, Horn, and Johnson}]{horn1990matrix}
Horn, R.~A., Horn, R.~A., Johnson, C.~R., Matrix analysis. Cambridge university press. (1990)



\bibitem{MUMPS1}
  P. R. Amestoy , I. S. Duff,
  \emph{ Multifrontal parallel distributed symmetric and unsymmetric solvers},
  Computer Methods in Applied Mechanics and Engineering, 184
  (2000)
  501-520.

\bibitem{MUMPS2}
P. R. Amestoy, I. S. Duff, J. Koster, J.Y. L'Excellent,
\emph{ A fully asynchronous multifrontal solver using distributed dynamic scheduling},
SIAM Journal of Matrix Analysis and Applications, 1(23)
  (2001) 15-41.

\bibitem{MUMPS3}
P. R. Amestoy, A. Guermouche, J.-Y. L'Excellent, S. Pralet,
\emph{ Hybrid scheduling for the parallel solution of linear systems},
Computer Methods in Applied Mechanics and Engineering, 2(32)
(2001) 136-156.

\bibitem{weak}
Y. Bazilevs, C. Michler,  V. M. Calo, T.J.R. Hughes, Isogeometric variational multiscale modeling of wall-bounded turbulent flows with weakly enforced boundary conditions on unstretched meshes, Computer Methods in Applied Mechanics and Engineering 199(13-16) (2008) 780-790.

\bibitem{timedep}
M. \L{}o\'{s}, , J. Mu\~{n}oz-Matute, I. Muga, M. Paszy\'{n}ski,  Isogeometric Residual Minimization Method (iGRM) with Direction Splitting for Non-Stationary Advection-Diffusion Problems, Computers and Mathematics with Applications, \revMP{79(2) (2020) 213-229.}
\bibitem{IGA2}
T.J.R. Hughes, J.A. Cottrell, Y. Bazilevs, Isogeometric analysis: CAD, finite elements, NURBS, exact geometry and mesh refinement, Computer Methods in Applied Mechanics and Engineering, (39-41) 4135-4195 (2005)

\bibitem{Sergio}
V. M. Calo, A. Ern, I. Muga, S. Rojas, An adaptive stabilized finite element method based on residual minimization, arXiv:1907.12605 (2019)

\bibitem{Method1}
R. E. Bank, B. D. Welfert, H. Yserentant, A class of iterative methods for solving saddle point problems, Numerische Mathematik, 56(7) \revMP{(1989)} 645-666.


\bibitem{Method2}
M. Arioli, C. Kruse, U. Rude, N. Tardieu, An iterative generalized Golub-Kahan algorithm for problems in structural mechanics,  arXiv:1808.07677v1

\bibitem{yzbeta1}
T. E. Tezduyar, M. Senga, Stabilization and shock capturing
parameters in SUPG formulation of compressible flows, Computer Methods in Applied Mechanics and Engineering, 195 (2006) 1621–1632.

\bibitem{yzbeta2}
T. E. Tezduyar, M. Senga, SUPG finite element computation
of inviscid supersonic flows with YZb shock-capturing, Computers
\& Fluids, 36 (2007) 147–159.

\end{thebibliography}
\end{document}